\documentclass[10pt]{article}
\usepackage{amssymb}
\usepackage{amsthm}
\usepackage{amsmath}
\usepackage{multirow}
\usepackage{longtable}

   \usepackage{epsfig} 
  \usepackage{latexsym}
\usepackage{amsmath,amsfonts}

\usepackage[swedish]{babel}
\usepackage[latin1]{inputenc}

\usepackage{graphicx}

\renewcommand{\theequation}{\arabic{section}.\arabic{equation}}
 
\newcommand{\figura}[2]{%
  \medskip\begin{center}%
  \resizebox!{#2mm}{\includegraphics{#1}}%
  \end{center}\medskip%
}

\newcommand{\figurb}[2]{%
  \ifx\pdfoutput\undefined
    \figura{#1.eps}{#2}\else
    \ifcase\pdfoutput
      \figura{#1.eps}{#2}\else
      \figura{#1.pdf}{#2}\fi
   \fi
}

\newtheorem{theorem}{Theorem}[section]

\newtheorem{prop}[theorem]{Proposition}
\newtheorem{exempel}[theorem]{Example}

\newtheorem{antag}{Assumption}[section]
\newtheorem{lemma}[theorem]{Lemma}

\renewcommand{\theequation}{\arabic{section}.\arabic{equation}}
  
\begin{document}
  
\title{  On a Multilocus  Wright-Fisher    Model with   Mutation and a  Svirezhev-Shahshahani 
Gradient-like Selection Dynamics }
\author{Erik Aurell$^{1,2}$,  Magnus Ekeberg, Timo Koski$^{3, \ast}$ } 
\maketitle 

\noindent $^{1}$ Division  of  Computational  Science and  Technology, Department  of Computer Science, KTH-Royal Institute of Technology, SE-106 91 Stockholm, Sweden \\
 $^{2}$Department  of Applied Physics and Department of Computer Science, Aalto University, FIN-00076 Aalto, Finland \\ 
 $^{3}$Department of Mathematics, KTH-Royal Institute of Technology, SE-100 44  Stockholm, Sweden. \\
$^{\ast}$Corresponding author: tjtkoski@kth.se, Tel: +46-8-790 71 34  

\vspace{1.0cm}

\noindent In this paper  we  introduce  a multilocus diffusion   model of  a population of $N$ haploid, asexually reproducing individuals. The model includes parent-dependent mutation and interlocus selection, the latter limited to pairwise relationships but among  a large number of simultaneous loci. The diffusion  is expressed as a  
system of  stochastic differential equations (SDEs) that are coupled in  the drift functions  through a  Shahshahani 
gradient-like  structure for interlocus selection.     The   system of  SDEs  is derived from a sequence of Markov chains by weak convergence.    We find  the explicit  stationary (invariant) density  by solving the corresponding stationary Fokker-Planck equation 
under parent-independent mutation, i.e., Kingman's house-of-cards mutation. The density  formula enables us  to readily  construct  families of   Wright-Fisher models corresponding  to networks of  loci.  \\  

\vspace{10mm}


\section{Introduction} 
\setcounter{equation}{0}

Recent and future advances of  biotechnology will produce time series data of  allele frequencies, see 
\cite{steinrucken2014novel}, \cite{tataru2017statistical}. 
 The  coupled  Wright-Fisher model below  is a step towards  developing tools for analysis of 
 such data.   

 The   Wright-Fisher model for a single locus has been  extensively studied in population genetics   c.f.  \cite{etheridge2011some},   \cite{hart1997principles} or \cite[pp. 92$-$99]{ewens2012}. 
  An alternative equivalent formulation of the  Wright-Fisher model  is given in \cite{kingman1980mathematics}.
The  overview in \cite{blythe2007}  presents  
applications to    ecology   and linguistics, too.  

One of the present authors and his co-workers   used recently  a Potts model  \cite{ekeberg2013improved}  in structural learning and  analysis of coupled loci for the {\em Pneumococcus}  derived from a whole genome alignment \cite{skwark2017interacting}. 
  We introduce here    a Potts-like  structure of interactions  consisting of  first and second order  interactions as in the Potts model which will produce an explicit selective   fitness term, which is in fact a Shahshahani 
gradient-like  structure,   as is shown below.

 This is related  to the  work    
 on   Quasi-Linkage
Equilibrium   in which the dynamics of the full genotype distribution, including correlations
between alleles at different loci, is given a parametric form by the allele frequencies.      
 In \cite[eqns (41)-(43)]{neher2011statistical}  Neher  and  Shraiman   find  in the Quasi-Linkage Equilibrium   approximation a probability density with   a  broadly similar structure as the one expressed in this paper. The  technical details are, however, completely  different from this work.  

%

 Here  we prove   the weak convergence of a sequence of Markov chains with the allele frequencies of all loci 
 as the state space via  conditionally independent locuswise  multinomial sampling  to the limiting diffusion by  basically applying   the direct techniques of \cite{shimakura1981formulas} and  \cite{sato1978diffusion}. Those techniques were  tailored  especially for Wright-Fisher models.  The  general   methods are by   semigroups  of operators \cite{ethier2009markov} or 
 Stein$^{,}$s method \cite{barbour1990stein}. Diffusion approximation via  duality with a  coalescent, see    \cite{mohle2001forward} is involved are we are   relying on some results  due to  \cite{shiga1981diffusion}    obtained  by   moment  duality.

 The uniqueness of the probability measure induced  by the limiting diffusion (uniqueness of the martingale problem)  is established  by a change of drift  technique and Girsanov$^{,}$s  theorem.  The possibility of  this change of the drift  depends on the additive  appearance  of   Shahshahani 
gradient-like  structure derived from  a fitness  potential in the drift function.

%
%
The paper is organized as follows. 

%

\section{An Outline: The   Wright-Fisher SDEs coupled by a Svirezhev-Shahshahani  gradient form  }  
 \setcounter{equation}{0}

\subsection{The Multiallele Wright-Fisher Process} 
\setcounter{equation}{0}
We start with   the (single locus) M-allele  Wright-Fisher diffusion process by means of a system of stochastic differential  equations.    
   
Denote by  $ {\bf x}= (x_{1}, \ldots, x_{M-1})$ the generic element $ \in {\bf K}$,   the probability simplex    
 \begin{equation}\label{isimplex1} 
 {\bf  K} : = \{  {\bf x} \in  R^{M-1} \mid   x_{l}  \geq 0, \sum_{l=1}^{M -1} x_{l}  \leq 1 \}. 
 \end{equation}
In this     $x_{M} $   serves as    a shorthand notation for $1-\sum_{k=1}^{M -1}x_{k} $.  
Let $\mu({\bf x})$  be a real vector valued function  of ${\bf x} \in {\bf K} \mapsto \mathbf{R}^{M-1}$  satisfying suitable assumptions (specific forms will be  encountered   later)  and 
let for  ${\bf x} \in {\bf K}$ 
\begin{equation}\label{genetic_drift}
d_{kl} ({\bf x})=   \begin{cases} x_k (1-x_k ), k=l \\ -x_k x_l, k \neq l \end{cases}.
\end{equation} 
For ease of writing, $d_{kl} ({\bf x})$ is  sometimes  denoted by $d_{kl}$.
${\bf D}\left({\bf x}\right)$  is the $(M -1) \times (M -1) $  covariance matrix with  arrays given in  (\ref{genetic_drift}).    Let $ X(t)$ be a  random variable
$$ 
 X(t) \in  {\bf  K}, \quad \text{for $t >0$},   
$$ 
 such  that the process $ \mathbf{X}= \{X(t) | t >0 \}$  satisfies (e.g., in the strong sense) the multivariate stochastic differential equation (SDE) 
 \begin{equation}\label{bigd5}  
  dX(t) =  {\bf \mu}(X(t))dt   +  {\bf D}^{1/2}\left(X(t)\right)d {\bf W}(t), 
 \end{equation} 
 where  $ {\bf W}= \{ {\bf W}(t) = \left( W_{1}(t), \ldots, W_{M-1}(t) \right)| t >0\}$ is an $M-1$ dimensional
 Wiener process.  The multivariate diffusion process corresponding  to this SDE  is    called  the  {\em  M-allele  Wright-Fisher diffusion}.    
   $\mu({\bf x})$ is the drift function and  ${\bf D}\left({\bf x}\right)$ is the diffusion (also known as  genetic drift) matrix of the 
 M-allele  Wright-Fisher diffusion.    The presence of  (\ref{genetic_drift})  is the   hallmark of the Wright-Fisher diffusions. 

 The   weak convergence of a sequence  of  Markov chains   to the $M$-allele   Wright-Fisher diffusion  on  ${\bf K} $   is proved  in    \cite[p. 62]{dawson2010introductory},   \cite{sato1976diffusion}, \cite{sato1978diffusion}, \cite{shiga1981diffusion},  and \cite{shimakura1981formulas}.  K.~Sato proved also in  \cite{sato1976diffusion} that  under some weak conditions on $ \mu ({\bf x}) $'s (satisfied in this paper) that the limiting process  stays in    ${\bf K}$. 

The general  methods for  weak convergence of a sequence  of Markov chains  to a multivariate diffusion  as given in      \cite[Chapter 11.2-11.3]{stroock2007multidimensional} are not directly applicable in the cited work. The reasons for this are as follows.  First, the 
diffusion function $d_{kl} ({\bf x})$ yields a degenerate  elliptic operator, see Appendix \ref{diffinverse}. This leads to the question of conditions on  the functions $\mu ({\bf x})$   for the   
martingale problem associated to the generator of   (\ref{bigd5})   to have  a 
unique solution.  These conditions  are  treated  in    \cite{ethier1976class}, see also \cite{shiga1981diffusion}.  The   degenerate elliptic  operators corresponding to   Wright-Fisher diffusion are  studied    in  
\cite{epstein2013degenerate}.     Second, the probability  simplex  
${\bf K} $  has   intricate  geometric properties, c.f.,    \cite{HofrichterJulian2014}.      
Third, the results on  the existence of an invariant measure for a multidimensional diffusion, see   
\cite{bhattacharya1978criteria}, are for these            reasons not applicable to    Wright-Fisher diffusions.

\subsection{The Multilocus  and Multiallele    Wright -Fisher SDE} \label{rewrite}  

Let    $x_{k}^{(i)}$  for $i=1, \ldots, L$   be  the frequency  of  allele type $k$ at locus $i$ in a finite population, and  
 \begin{equation}\label{isimplex} 
 {\bf  K}_{i} = \{ x \in  R^{M_i-1} \mid   x_{l}^{(i)} \geq 0, \sum_{l=1}^{M_{i}-1} x_{l}^{(i)} \leq 1 \}, 
 \end{equation}
and    $x_{M_i}^{(i)}$   equals    $=1-\sum_{k=1}^{M_i-1}x_{k}^{(i)}$.  
 We stack  these    into a single     
column vector,  
 \begin{eqnarray}\label{frekekv} 
{\bf x}&= &(x^{(1)}, x^{(2)},\cdots, x^{(L)}) \\ \nonumber  
& & \\ \nonumber 
& =&  (x_{1}^{(1)},\ldots,x_{M_1-1}^{(1)},x_{1}^{(2)},\ldots,x_{M_2-1}^{(2)},\ldots, x_{1}^{(L)},\ldots,x_{M_L-1}^{(L)}).  \nonumber  
\end{eqnarray}  
Hereafter we  have   ${\bf  K} :=  \times_{i=1}^{L}{\bf  K}_{i}$ instead of (\ref{isimplex1}), so that 
${\bf x} \in {\bf  K}$.   
We define with (\ref{genetic_drift}) for     for all  $i$, $l$ and $k$ $d_{kl}^{(i)}({\bf x})$  by  
\begin{equation}\label{driftlemma44} 
d_{kl}^{(i)}({\bf x}) \stackrel{\rm def}{=}   \begin{cases} x_k^{(i)}(1-x_k^{(i)}), k=l \\ -x_k^{(i)}x_l^{(i)}, k \neq l \end{cases}.
\end{equation}
 ${\bf D}^{(i)}\left({\bf x}\right)$  is the $(M_{i}-1) \times (M_{i}-1) $  covariance matrix with  arrays given in  (\ref{driftlemma44}), which means  that    
${\bf D}^{(i)}\left({\bf x}\right)$ depends only on the  allele frequencies   $x_{k}^{(i)}$ at locus $i$. 
Next, ${\bf O}_{ij}$  is a $(M_{i}-1) \times (M_{j}-1) $ matrix of zeroes.  We introduce now  the  quadratic  $ \sum_{i=1}^{L} (M_{i}-1)  \times  \sum_{j=1}^{L} (M_{j}-1) $ block diagonal diffusion matrix 
\begin{equation}\label{bigd}
{\bf D}\left({\bf x}\right)\stackrel{\rm def}{=} \left(  \begin{array}{ccccc}  
{\bf D}^{(1)}\left({\bf x}\right) & {\bf O}_{12} & \ldots  &  \ldots & {\bf O}_{1L} \\
 {\bf O}_{21} &  {\bf D}^{(2)}\left({\bf x}\right) &  \ldots &  \ldots & {\bf O}_{2L} \\ 
 \vdots      &  \vdots  &  \ddots  & \vdots  &  \vdots    \\  
 {\bf O}_{L1} &  \ldots  &  \ldots  & \ldots      & {\bf D}^{(L)}\left({\bf x} \right) \end{array} \right).  
\end{equation}  
  Next  we consider for $t >0$ 
$$
 X(t) =  \left(   X^{(1)}\left(t\right)   \ldots  X^{(L)}\left(t\right)   
  \right), 
 $$ 
 where each  random variable $X^{(i)}\left(t\right) $ assumes   values in  ${\bf K}_{i}$, respectively. 

In this paper we determine a sequence of Markov chains that after a scaling and interpolation of discrete time will converge weakly to a ${\bf K}$-valued stochastic process 
 $ \{ X(t) \mid t \geq  0   \}$  that satisfies an SDE of the form 
 \begin{equation}\label{bigd55} 
  dX(t) =  {\bf \mu}(X(t))dt + {\bf G} \left(X(t)\right)dt  +  {\bf D}^{1/2}\left(X(t)\right)d {\bf W}(t), 
 \end{equation} 
 where )$^{T}$ is the vector transpose) 
 $$
{\bf W}(t)^{T}= \left( W^{(1)}(t), \ldots, W^{(L)}(t) \right),  
 $$
and   $ W^{(i)}$  are independent    Wiener processes, each with $M_{i}- 1$  independent component Wiener processes.  The drift function  is thus   of the form  ${\bf \mu}({\bf x})  + {\bf G} \left({\bf x}\right)$, where
$$
{\bf \mu} \left({\bf x}\right) = \left({\bf \mu}^{(1)}\left({\bf x}\right), \ldots, {\bf \mu}^{(L)}\left({\bf x}\right)\right),
$$
here each  ${\bf \mu}^{(i)}\left({\bf x}\right)$ depends only on $x^{(i)}$, and 
$$
{\bf G}^{T} \left({\bf x}\right) = \left( {\bf G}^{(1)}\left({\bf x}\right), \ldots, {\bf G}^{(L)}\left({\bf x}\right)\right),
$$
where
\begin{equation}\label{eadrftfact9990}
{\bf G}^{(i)}\left({\bf x}\right) =  {\bf D}^{(i)}\left({\bf x}\right) \nabla_{ x^{(i)}} V \left({\bf x}\right),
  \end{equation} 
with the gradients  
$$
\nabla_{ x^{(i)}} V \left({\bf x}\right)^{T}= \left( 
   V^{'}_{ x_{l}^{(i)}} \left({\bf x}\right),   V^{'}_{ x_{2}^{(i)}} \left({\bf x}\right) ,  
    \ldots,    V^{'}_{ x_{M_{i}-1}^{(i)}} \left({\bf x}\right)   \right). 
$$
for a function $V \left({\bf x}\right)$ (=fitness potential) explicitly constructed below  of at most pairwise interaction between the $L$ loci.   In view of the definition of $ {\bf D}\left({\bf x} \right)$ 
in  (\ref{bigd})  it   holds that  
\begin{equation}\label{eadrftfact9991}
{\bf G} \left({\bf x}\right) =  {\bf D}\left({\bf x} \right)  \nabla_{{\bf x}} V \left({\bf x}\right),     
  \end{equation} 
  where 
$$
  \nabla_{{\bf x}} V \left({\bf x}\right)^{T}=\left( \nabla_{ x^{(1)}} V \left({\bf x}\right) , \nabla_{ x^{(2)}} V \left({\bf x}\right), \ldots, \nabla_{ x^{(L)}} V \left({\bf x}\right)   \right). 
$$  
The result is  that the system of SDEs  in (\ref{bigd55})  is  formally   given by   
\begin{equation}\label{bigdsd} 
  dX(t) =  {\bf \mu}(X(t))dt + {\bf D}(X(t))  \nabla_{{\bf x}} V \left(X(t)\right)dt  +  {\bf D}^{1/2}\left(X(t)\right)d {\bf W}(t).  
 \end{equation}
 The  drift  function (\ref{bigdsd}) will be seen to have the following   intuitive meaning; the  terms  inside any 
$\mu^{(i)}({\bf x})$  are accounting   for 'mutational flow' to and from allele type $k$ respectively at locus $i$, and the second  term ${\bf D}(X(t))  \nabla_{{\bf x}} V \left({\bf x}\right)$ represents selection for the current allele type at the current locus.
 This means  in practice  disregarding the cross-effects of selection and mutation and instead viewing them as independent mechanisms.

The  term   ${\bf D}({\bf x})  \nabla_{{\bf x}} V \left({\bf x}\right)$  will     contain  the (at most pairwise) interactions between    the various loci and their allele types. The quantities in both   ${\bf \mu}\left({\bf x}\right)$   and $ {\bf D}^{1/2}\left({\bf x}\right)$ are  decomposed to act on respective   locus  without interaction.  Hence, were   ${\bf D}({\bf x})  \nabla_{{\bf x}} V \left({\bf x}\right)$  to be removed ($\approx$ weak selection),     $ \{ X(t) \mid t \geq  0   \}$  will consist     
of an  $L$ independent  M-allele  Wright-Fisher diffusion processes  as given in (\ref{bigd55}).  If  ${\bf \mu}({\bf x})$ is the zero vector of appropriate dimensions, then 
\begin{equation}\label{bigdsd2} 
  dX(t) =    {\bf D}(X(t))  \nabla_{{\bf x}} V \left(X(t)\right)dt  +  {\bf D}^{1/2}\left(X(t)\right)d {\bf W}(t).  
 \end{equation}
is a {\em Svirezhev-Shahshahani gradient   SDE}  for the allelic frequency dynamics.  The Svirezhev-Shahshahani gradient is discussed in  \cite[pp.42-43]{burger2000mathematical} and \cite[p.~222-223, p.~303]{svirezhev2012fundamentals},  the paper \cite{huillet2017random} introduces a Wright-Fisher model, where  a  Svirezhev-Shahshahani gradient lies in   the  diffusion matrix.   

In order to make these statements a bit more transparent at this stage,  we look at a two-locus two alleles ($L=2$,  $M_{i}=2$, 
$i=1,2$)   Wright-Fisher model with selection but no  mutation,  \cite[ch. 15, section F, example (b)]{karlin1981second}.  
 
\begin{exempel} \label{karlinsvishahshahani}
{\rm 
The    Wright-Fisher model with selection but without mutation is  in loc.cit. given as the pair of coupled 
SDEs  
\begin{equation} \label{2dimwfdiff11}
\left \{ \begin{array}{cc}
dX^{(1)}_{t} =    
h X^{(1)}_{t}(1-X^{(1)}_{t}) X^{(2)}_{t}dt + \sqrt{ X^{(1)}_{t}(1-X^{(1)}_{t}) }dW^{(1)}_{t} \\
dX^{(2)}_{t} =     
h X^{(2)}_{t}(1-X^{(2)}_{t}) X^{(1)}_{t}dt + \sqrt{ X^{(2)}_{t}(1-X^{(2)}_{t}) }dW^{(2)}_{t}.    
\end{array} \right. 
\end{equation} 
We  find  now $V({\bf x})$ so  that the drift in (\ref{2dimwfdiff11})  is of the form (\ref{bigdsd2}). 
Here ${\bf x}=  (x^{(1)} ,  x^{(2)})$, since $M_{i}=2 $, and $ 0 \leq x^{(i)} \leq 1$, $i=1,2$.   We have here    a special case of the  construction  in section \ref{potshasvhi}, the exact details for  this special case are found in Example \ref{exempel3}. Let  $V({\bf x}) $ be a function of $ \times_{i=1}^{2}{\bf  K}_{i}$ to   $\mathbf{R}$ given by    
$$
V({\bf x}) =    h  x^{(1)}  x^{(2)}.  
$$
Then the gradient is      
$$  
\nabla_{{\bf x}} V \left({\bf x}\right)=  \left(  \begin{array}{c}  
 V^{'}_{x^{(1)}}    \\
V^{'}_{x^{(2)}}  \end{array} \right)  =  \left(  \begin{array}{c}  
 h  x^{(2)}   \\
 h  x^{(1)} \end{array} \right).      
$$
By (\ref{bigd})  
\begin{equation}\label{bigd33}
{\bf D}\left({\bf x}\right)=  \left(  \begin{array}{cc}  
{\bf D}^{(1)}\left({\bf x}\right) & {\bf O}_{12}  \\
 {\bf O}_{21} &  {\bf D}^{(2)}\left({\bf x}\right)  \end{array} \right),   
\end{equation} 
where  
$$
{\bf D}^{(i)}\left({\bf x}\right) =   x^{(i)} (1- x^{(i)}), i=1,2,  {\bf O}_{12}={\bf O}_{21}=0. 
$$
Then 
\begin{equation}\label{svshd33} 
  {\bf D}\left({\bf x} \right)  \nabla_{{\bf x}} V \left({\bf x}\right)=   \left(  \begin{array}{c}  
 h   x^{(1)} (1- x^{(1)})x^{(2)}   \\
 h  x^{(2)} (1- x^{(2)})  x^{(1)} \end{array} \right).      
 \end{equation} 
Hence (\ref{2dimwfdiff11}) is an instance  of (\ref{bigdsd2}). We rewrite (\ref{svshd33}) as 
\begin{equation}\label{svshd335} 
  {\bf D}\left({\bf x} \right)  \nabla_{{\bf x}} V \left({\bf x}\right)=   \left(  \begin{array}{c}  
    x^{(1)} (hx^{(2)}  - h x^{(1)}x^{(2)})   \\
   x^{(2)} (h x^{(1)}  - hx^{(1)} x^{(2)})  \end{array} \right) =  
  \left(  \begin{array}{c}  
    x^{(1)} ( V^{'}_{x^{(1)}}  -V({\bf x}))   \\
   x^{(2)} ( V^{'}_{x^{(2)}}  - V({\bf x}))  \end{array} \right).     
\end{equation}    
At this point  the expressions can be related to  (an underlying)   deterministic dynamics for    $  {\bf x}(t) = (x^{(1)}(t),x^{(2)}(t))$, when (\ref{2dimwfdiff11}) is written  as 
\begin{equation} \label{burgerdet}
\left \{ \begin{array}{cc}
\frac{d}{dt} x^{(1)}(t)=     x^{(1)}(t) ( V^{'}_{x^{(1)}}({\bf x}(t))  -V({\bf x}(t))) 
  \\
\frac{d}{dt} x^{(2)}(t) =    x^{(2)}(t) ( V^{'}_{x^{(2)}}({\bf x}(t))  -V({\bf x}(t))),   
     
\end{array} \right. 
\end{equation}                                 
which is a so-called replicator equation, and in general mathematical terms, this is a gradient system,  
the Svirezhev-Shahshahani gradient system, see, e.g., \cite[p.~103, pp.~349$-$351]{burger2000mathematical}

}\end{exempel}  
%

\section{A Markov Chain of Allele Frequencies: Assumptions and the Transition Probability}\label{transproball}   
 \setcounter{equation}{0}

We consider a very large   population of $N$ haploid, asexually reproducing individuals.  
A new generation   is brought to life as follows: first we sample, independently and with replacement, $N$ new individuals from the previous generation, with the probability of choosing an individual of haplotype $\sigma$ scaled by a selection coefficient $v_{\sigma}$, c.f., the next section. Subsequently, we let a mutation event occur at each locus (independently).

There  are $L \geq 1$ loci. The symbol  ${\sigma}$ specifies  the allele types at these $L$ loci, ${\sigma}=(\sigma_{j_1}, \cdots,\sigma_{j_L})$, where $\sigma_{j_i} \in {\cal S}_{i}= [1,\cdots,M_i]$, i.e., we accept different numbers of possible allele types at different loci.  We write   ${\sigma} \in  \times_{i=1}^{L}{\cal S}_{i} $.    
An individual $r$ is represented by its allelic vector ${\sigma}(r) \in  \times_{i=1}^{L}{\cal S}_{i} $, $r=1, \ldots,N$. 

 Let  $x_{k}^{(i)}$ be  the frequency  of individuals carrying allele type $k$ at locus $i$. Thus  
\begin{equation}\label{iverson}
x_{k}^{(i)} = \frac{1}{N} \sum_{r=1  }^{N}\delta_{\sigma_i(r),k}, \quad   k \in [1,\cdots,M_i],    
\end{equation}   
where  $\delta_{ i,k}$ is the Kronecker delta (or the Iverson bracket). 
Let us set 
\begin{equation}\label{vekt1}
x^{(i)}= \left( x_{1}^{(i)},\cdots,x_{M_i-1}^{(i)}   \right).  
\end{equation} 
In sections \ref{solutionstmk} and \ref{exempel}  $x^{(i)}$ will be treated  as an  $M_{i}-1 \times 1$ vector, but at this stage this interpretation is not operationallly necessary.     Thus  $x^{(i)}$ lies in the simplex 
$ {\bf  K}_{i}$ in (\ref{isimplex}). 

Furthermore,  we have a set of non-negative   integers  or the  occupancy distribution of the $M_{i}$ alleles, 
\begin{equation}\label{frekekv21}
{\bf J}_{(i)}=  \left \{ j^{(i)}=(j_{1}^{(i)}, \ldots, j_{M_{i}}^{(i)})  \in Z_{+ \cup 0}^{M_{i}}  \mid  \sum_{k=1}^{M_{i}}j_{k}^{(i)} = N \right \}. 
\end{equation}
For example,  $j_{k}^{(i)} =  N \cdot x_{k}^{(i)}$.   
Thus  ${\bf J}_{(i)}$ can be  regarded as subset of ${\bf  K}_{i}$ consisting of all  the  lattice points with   mesh  $1/N$. The number of distinquishable  occupancy distributions  in  ${\bf J}_{(i)}$ is  equal to  $ \left( \begin{array}{cc} 
(M_{i} -1)+ N \\ 
N  
\end{array} 
\right)$.   
Corresponding to (\ref{frekekv}) we have 
 \begin{equation}\label{frekekv2} 
{\bf j} = (j^{(1)}, j^{(2)},\cdots, j^{(L)})  
\end{equation} 
and 
 \begin{equation}\label{frekekv3} 
{\bf j}  \in  \times_{i=1}^{L}{\bf J}_{(i)}.  
\end{equation} 
Let now  $n=0,1,2, \ldots $ represent  discrete time and $N>0$.  We consider a Markov chain,  homogeneous 
 in discrete   (scaled) time,     
 $Y^{(N)}  = \{ Y^{(N)}  \left( \frac{n}{N} \right)  \}_{n \in Z_{+}} $  with the state space  $ \times_{i=1}^{L}{\bf J}_{(i)}$.   
The transition probability  is  for any  ${\bf j} \in \times_{i=1}^{L}{\bf J}_{(i)} $ and ${\bf k}  
\in  \times_{i=1}^{L}{\bf J}_{(i)}$ denoted  by  
$$
P_{{\bf j},{\bf k} }= P \left( Y^{(N)} \left( \frac{n+1}{N} \right)= {\bf k}  \mid  Y^{(N)} \left( \frac{n}{N} \right)= {\bf j} \right).    
$$ 
 For each locus $i$ there is  the random process  of occupation numbers    
  $Y^{(N)}_{i}  = \{ Y^{(N)}_{i}\left( \frac{n}{N} \right)  \}_{n \in Z_{+}} $  with the state space  ${\bf J}_{(i)}$ so that 
  $  Y^{(N)} \left( \frac{n}{N} \right) = \left(  Y^{(N)}_{i}\left( \frac{n}{N} \right) \right)_{1 \leq i  \leq L}   $.      
\begin{antag}\label{cindep} 
{ \rm   The locus-wise component  processes  $Y^{(N)}_{i}$ at any time $n+1$ are conditionally independent of  each other given the  process  $Y^{(N)}$ at time $n$: 
 For any ${\bf j} \in \times_{i=1}^{L}{\bf J}_{(i)} $ and ${\bf k}  
\in  \times_{i=1}^{L}{\bf J}_{(i)}$ and any $n \geq 0$  it  holds that  
\begin{equation}\label{condindpe}
P_{{\bf j},{\bf k} }=     \prod_{i=1}^{L} P \left( Y^{(N)}_{i} \left( \frac{n+1}{N} \right) =     
 k^{(i)} \mid Y^{(N)} \left( \frac{n}{N} \right)= {\bf j}   \right).   
\end{equation}}
\end{antag}  \qed \\  
   There is  clearly   for any locus $i$   a vector process of random occupation numbers, i.e.,  
\begin{equation}\label{onelinterpol}
Y^{(N)}_{i}(n) = \left(Y_{ i,1}^{(N)} \left( \frac{n}{N} \right),Y_{i,2}^{(N)} \left( \frac{n}{N} \right), \ldots, Y_{i, M_{i}}^{(N)}\left( \frac{n}{N} \right)    \right),  
\end{equation}
 where, if we sum over the allele types at any  locus $i$,    $ \sum_{r=1}^{M_{i}} Y_{ i,r}^{(N)} \left( \frac{n}{N} \right) =   N$.

 The  transition  probability       
 $P_{{\bf j},{\bf k}}$ is given by  specification of  the  conditional probabilities $P \left( Y^{(N)}_{i} \left(\frac{n+1}{N} \right) =     
 k^{(i)} \mid Y \left(\frac{n}{N} \right)= {\bf j}   \right) $.  If  $  {\bf j} \in  \times_{i=1}^{L}{\bf J}_{(i)}  $ and  $k \in  {\bf J}_{(i)}$, we take  
\begin{equation}\label{transitionmech}
P \left( Y^{(N)}_{i} \left(\frac{n+1}{N} \right)= k^{(i)}] \mid  Y^{(N)} \left(\frac{n}{N} \right)= {\bf j}  \right) = \frac{N!}{k^{(i)}_{1}! \cdots k^{(i)}_{M_{i}}!} 
 p^{(i)}_1({\bf j})^{k^{(i)}_{1}} \cdots  p^{(i)}_{M_{i}}({\bf j})^{k^{(i)}_{M_{i}} },    
\end{equation}  
which is  a  multinomial distribution, where  $p^{(i)}_k({\bf j})$ is the probability  of  the  allele type $k$ at locus $i$.  Our goal is  to express the dependence of      
$$
[p^{(i)}_1( {\bf j}),\cdots,p^{(i)}_{M_{i}-1}( {\bf j})] \in  {\bf  K}_{i} 
$$           
on ${\bf j}$.   For any  ${\bf j}$ we may compute the corresponding  relative frequency vector ${\bf x}\in   {\bf  K} $,  of the form (\ref{frekekv}) for  the current population. 
We drop, for simplicity  of  expression,  the   dependence on the  $L$ occupancy distributions  in   $  {\bf  j }$ in the  formulas that  in the rest of this section.

 We start by  the fraction of individuals with   haplotype 
 ${\sigma} = \left( \sigma_{j_1}, \cdots,\sigma_{j_L}   \right) $ denoted  by   $f(\sigma)$. This is  simply the product of the (relative)  population frequencies for   an  allele $\sigma_{j_{i}} \in {\cal S}_{i}= [1,\cdots,M_i]$  for  each locus,
\begin{equation}\label{frac_of_haplo}
 f(\sigma) = \prod_{i=1}^L x_{\sigma_{j_{i}}}^{(i)}, 
\end{equation} 
  i.e. this is the  product of  fractions   picked from  (\ref{frekekv})  according  to $ (\sigma_{j_1}, \cdots,\sigma_{j_L})$ and (\ref{iverson}). 

Let us next  define $f^{(i)}_k(\sigma)$ as the conditional frequency of the haplotype $\sigma$ given $\sigma_i=k \in {\cal S}_{i}$, which is simply the same expression as above but with Kronecker delta  $\delta_{\sigma_i,k}$ substituted for $x_{\sigma_i}^{(i)}$,
\begin{equation}\label{condallelfrekv}  
 f^{(i)}_k(\sigma) = \delta_{\sigma_i,k}\prod_{\underset{ j \neq i }{j=1}}^L x_{\sigma_{l_j}}^{(j)} =\frac{\delta_{\sigma_i,k}}{x_k^{(i)}}f(\sigma).
\end{equation}
If  $x_k^{(i)}=0$  for the  population at a time, then $f(\sigma)=0$ by (\ref{frac_of_haplo}),  we can take 
by convention $ f^{(i)}_k(\sigma)=0 $. 
The function $v_{\sigma}$ gives the viability of an individual with the allelic vector $\sigma$. We set 
\begin{equation}\label{selecstrength}
\bar{v} = \sum_{\sigma    \in  \times_{i=1}^{L}{\cal S}_{i}}  f(\sigma)v_{\sigma}, \bar{v}_k^{(i)} = \sum_{\sigma  \in  \times_{i=1}^{L}{\cal S}_{i}} f_k^{(i)}(\sigma)v_{\sigma}.  
\end{equation}  
$\bar{v}$ can be interpreted as the average selection strength for the population as a whole (in the current state ${\bf x}\in  \times_{i=1}^{L}{\bf  K}_{i}$, and $\bar{v}_k^{(i)}$ is   the average selection strength for allele type $k$ at locus $i$.

Without mutation, the probability of drawing an individual with allele type $k$ at locus $i$ is 
 \begin{equation}
q^{(i)}_k =\frac{\sum_{\sigma  \in  \times_{i=1}^{L}{\cal S}_{i}} \delta_{\sigma_i,k} f(\sigma)v_{\sigma}}{\sum_{\sigma  \in  \times_{i=1}^{L}{\cal S}_{i}} f(\sigma)v_{\sigma}},
\end{equation}
which can be in view of (\ref{selecstrength})  written as
\begin{equation}\label{q_def}
q^{(i)}_k =x^{(i)}_k\frac{\bar{v}_k^{(i)}}{\bar{v}}.
\end{equation}
Let next  $\upsilon_{lk}^{(i)}$  be the probability  that an $l$-allele at locus $i$ mutates   to 
 an  $k$-allele at locus  $i$ after the selection event. In this  $\upsilon_{lk}^{(i)}$  does not  depend  on  ${\bf j}$.   
Any other allele type can mutate into a $k$-allele, as governed by the probabilities $\upsilon_{lk}^{(i)}$, so the final probability of ending up with a $k$-allele at locus $i$ is $p^{(i)}_k = \sum_{l=1}^{M_i} \upsilon_{lk}^{(i)}q_l^{(i)}$.
As  the probability of no mutation can be written as $\upsilon^{(i)}_{kk} = 1- \sum_{\underset{l \neq k}{l=1}}^{M}\upsilon_{lk}^{(i)}$, we get
\begin{equation} \label{finalprob} 
\begin{aligned}
 p_k^{(i)} = \sum_{l=1}^{M_i} \upsilon_{lk}^{(i)}q^{(i)}_l = \sum_{\underset{l \neq k}{l=1}}^{M_i} \upsilon_{lk}^{(i)}q^{(i)}_l + \upsilon^{(i)}_{kk}q^{(i)}_k \\ = \sum_{\underset{l \neq k}{l=1}}^{M_i} \upsilon_{lk}^{(i)}q^{(i)}_l + \left (1- \sum_{\underset{l \neq k}{l=1}}^{M_i}\upsilon_{kl}^{(i)}\right )q^{(i)}_k
 =\sum_{\underset{l \neq k}{l=1}}^{M_i} \left [ \upsilon_{lk}^{(i)}q^{(i)}_l  - \upsilon_{lk}^{(i)}q^{(i)}_k  \right ] + q_k.\\
\end{aligned}
\end{equation}
By inserting (\ref{q_def}) we obtain
\begin{equation} 
\label{not_approximated1}
 p_k^{(i)}\left({\bf j} \right) =\sum_{\underset{l \neq k}{l=1}}^{M_i} \left [ \upsilon_{lk}^{(i)}x_l^{(i)} \left ( \frac{\bar{v}_l^{(i)}}{\bar{v}} \right ) - \upsilon_{kl}^{(i)}x_k^{(i)} \left ( \frac{\bar{v}_k^{(i)}}{\bar{v}} \right ) \right ] + x_k^{(i)} \left (\frac{\bar{v}_k^{(i)}}{\bar{v}}\right ).\\
\end{equation}
This completes the  description of the transition probability in (\ref{transitionmech}). 
We proceed  by re-scalings and translations of the  quantities in $ p_k^{(i)}\left({\bf j} \right)$  to get over to a  continuous time  SDE.

\section{The Drift Function  }\label{driftikin}  
 \setcounter{equation}{0}
\subsection{Scaling of the Transition Probability; The limiting  Drift function   } 
    
With    $L$ loci,  and   $M_{i}$ alleles at locus $i$, let us   consider   for all $l,k$ and all $i$ 
the parameters   $u_{lk}^{(i)}$, $\bar{m}$ and   $\bar{m}_k^{(i)}$ obtained by scaling and shifting with the inverse population size
the parameters of (\ref{not_approximated1}) as follows:      
\begin{equation}\label{scale1}
u_{lk}^{(i)} =  \frac{\upsilon_{lk}^{(i)}}{ \frac{1}{N}}, \quad   \bar{m} =  \frac{\bar{v}}{ \frac{1}{N}} - \frac{1}{ \frac{1}{N}}, \quad  
 \bar{m}_k^{(i)} =   \frac{ \bar{v}_k^{(i)}}{ \frac{1}{N}}- \frac{1}{ \frac{1}{N}}.    
\end{equation}  
 Then  it follows  by straightforward substitutions  in (\ref{not_approximated1}),   lemma \ref{driftlemma1} gives the detailed limiting  argument,  that  as $ N \rightarrow +\infty$ 
$$
 \frac{\left( p^{(i)}_k({\bf j}) -x^{(i)}_k \right)}{ \frac{1}{N} }\rightarrow  p_{k}^{(i)}({\bf x}), 
    $$ 
where  
\begin{equation}\label{scale12} 
 p_{k}^{(i)}({\bf x}) \stackrel{\rm def}{=} \sum_{\underset{l \neq k}{l=1}}^{M_i} \left [ u_{lk}^{(i)}x_l^{(i)} - u_{kl}^{(i)}x_k^{(i)} \right ]+  x_k^{(i)} \left ( \bar{m}_k^{(i)}-\bar{m}\right ).
\end{equation}    
For things to make the desired  sense,  $ p_{k}^{(i)}({\bf x})$ in (\ref{scale12}) should be the $k$th component   of the  vector     for locus $i$ in  
$ {\bf \mu}({\bf x})  + {\bf D}({\bf x})  \nabla_{{\bf x}} V \left({\bf x}\right)$   for some suitable $V \left({\bf x}\right)$.   
In order to establish  this   we start by making    an extra assumption,   
the   parent-independent mutation.  This assumption is also known as  {\em Kingman's house of cards assumption}, see  \cite{huillet2017random} and  \cite{kingman1980mathematics}.   
 
\begin{antag}\label{parent}
{\rm 
\begin{equation}\label{indmut223}  
 u_{lk}^{(i)} = u_{k}^{(i)}
 \end{equation}  
 for all $l$, $k$ and $i$. 
In addition we  assume that  
\begin{equation}\label{pos223}  
  u_{k}^{(i)}> 0 
 \end{equation}  
 for all   $k$ and $i$. 
 }
 \end{antag}  \qed \\

\begin{lemma}\label{ealemma4} 
Assume that  (\ref{indmut223}) holds.   
  Let  
\begin{equation}\label{eafact54}
g^{(i)}_{k} \stackrel{\rm def}{=} \sum_{\underset{l \neq k}{l=1}}^{M_i} \left [ u_{lk}^{(i)}x_l^{(i)} - u_{kl}^{(i)}x_k^{(i)} \right ] 
\end{equation}
    Then  
\begin{equation}\label{eafact55}  
g^{(i)}_{k}   = u_{k}^{(i)}-  \bar{u}     x_k^{(i)}, 
\end{equation}  
where  
\begin{equation}\label{eabar}
\bar{u}^{(i)} =  \sum_{ l=1}^{M_i}  u_{l}^{(i)}  
 \end{equation}
and  in  (\ref{scale12})  
\begin{equation}\label{driftlimit2222}    
  p_{k}^{(i)}({\bf x}) = g^{(i)}_{k}  +  x_k^{(i)} \left ( \bar{m}_k^{(i)}-\bar{m}\right ).     
\end{equation} 
 \end{lemma} 
 \noindent {\em Proof:}    
  By   (\ref{eafact54}) and  (\ref{indmut223}) we get 
$$
g^{(i)}_{k} =  \sum_{\underset{l \neq k}{l=1}}^{M_i} \left [ u_{k}^{(i)}x_l^{(i)} - u_{l}^{(i)}x_k^{(i)} \right ] 
 =  u_{k}^{(i)}(1-x_k^{(i)}) -  x_k^{(i)}\sum_{\underset{l \neq k}{l=1}}^{M_i}  u_{l}^{(i)}. 
$$
Now we evoke  $  \bar{u}^{(i)}  =  \sum_{ l=1}^{M_i}  u_{l}^{(i)}$ and get 
$$ 
= u_{k}^{(i)}(1-x_k^{(i)}) - \left( \bar{u}^{(i)} -u_{k}^{(i)} \right)  x_k^{(i)}  
 =  u_{k}^{(i)}  -   \bar{u}^{(i)}   x_k^{(i)}. 
$$
  \qed \\  
  
\noindent Herewith  we   set 
\begin{equation}\label{bigd2}
 \mu^{(i)}({\bf x}) \stackrel{\rm def}{=}  \left( \begin{array}{c} 
 
u^{(i)}_{1}  - \bar{u}^{(i)}  x^{(i)}_{1}    \\
 
  \vdots \\
  
 u^{(i)}_{M_{i}-1}  - \bar{u}^{(i)} x^{(i)}_{M_{i}-1}  \end{array} \right).   
\end{equation} 
 
\subsection{Population fitnesses}  

 Next we study the second term in the right hand side of  $ p_{k}^{(i)}({\bf x})$  in (\ref{scale12}), i.e., 
  $  x_k^{(i)} \left ( \bar{m}_k^{(i)}-\bar{m}\right )$.
Here $\bar{m}_k^{(i)}$ and $\bar{m}$ are expressing a   population   fitness, which  is a quantitative trait of a population,   thought of  as   mapping   the genotype to the expected reproductive success of an organism.  Here we have   
\begin{equation}
 \bar{m}= \sum_{\sigma} f(\sigma) m_{\sigma},   \bar{m}_k^{(i)}= \sum_{\sigma} f_k^{(i)}(\sigma) m_{\sigma},
 \end{equation}
 where  $ m_{\sigma}$ is the Potts-type interaction  map  
\begin{equation} 
\label{pairwise_selection}
 m_{\sigma}\stackrel{\rm def}{=}1+\sum_{r=1}^L h_r(\sigma_r)+\sum_{1 \leq r < s \leq L} J_{rs}(\sigma_r,\sigma_s),
\end{equation}
i.e.   selective interaction  between loci is limited to pairwise links. 
We assume  a double symmetry in the sense that 
\begin{antag}\label{symmm}
{\rm 
 \begin{equation}\label{symmetri} 
J_{rs}(k,l)=J_{sr}(l,k).
\end{equation} }
\end{antag} \qed \\
Then we get  (see  Appendix A) 
\begin{equation}
 \bar{m}= \sum_{\sigma    \in  \times_{i=1}^{L}{\cal S}_{i}}  f(\sigma)m_{\sigma}= 1+\sum_{r=1}^L\sum_{t=1}^{M_r}h_r(t)x_t^{(r)}+\sum_{1 \leq r < s \leq L}\sum_{t=1}^{M_r}\sum_{n=1}^{M_s}J_{rs}(t,n)x_t^{(r)}x_n^{(s)},
\end{equation} 
 and  (see  Appendix A), 
\begin{equation}\label{pairwise_selection2}
\begin{aligned}
 \bar{m}_k^{(i)}= 1&+h_i(k)+\sum_{\underset{r \neq i}{r=1}}^L\sum_{t=1}^{M_r}h_r(t)x_t^{(r)}\\
&+\sum_{\underset{r \neq i}{r=1}}^L \sum_{t=1}^{M_r}J_{ir}(k,t)x_t^{(r)}+\sum_{\underset{r,s \neq i}{1 \leq r < s \leq L}}\sum_{t=1}^{M_r}\sum_{n=1}^{M_s}J_{rs}(t,n)x_t^{(r)}x_n^{(s)},
\end{aligned}
\end{equation}
which yields 
\begin{eqnarray}\label{pairwise_selection3}
x_k^{(i)} \left ( \bar{m}_k^{(i)}-\bar{m} \right)& = & \\ \nonumber 
&  &  \\
& & x_k^{(i)} \left ( h_i(k)-\sum_{k'=1}^{M_i}h_i(k')x_{k'}^{(i)} + \sum_{\underset{r \neq i}{r=1}}^L\sum_{t=1}^{M_r} \left [J_{ir}(k,t)-\sum_{k'=1}^{M_i}J_{ir}(k',t)x_{k'}^{(i)} \right ]x_{t}^{(r)} \right). \nonumber 
\end{eqnarray}
  Let us set for simplicity of writing
 \begin{equation}\label{hejp2} 
\widetilde{h}_{i}(k) \stackrel{\rm def}{=} \left ( h_i(k)+\sum_{\underset{r \neq i}{r=1}}^{L}\sum_{m=1}^{M_r}  J_{ir}(k,m)x_{m}^{(r)}  \right ). 
\end{equation} 
  Then we have the following lemma.     
\begin{lemma}\label{ealemma3} 
For  $k=1, \ldots, M_{i}-1$, 
 \begin{equation}\label{eadrftfact22}
 x^{(i)} _{k} \left( \bar{m}_k^{(i)}-\bar{m} \right)=      \sum_{l=1}^{M_{i}} d_{kl}^{(i)} \widetilde{h}_{i}(l)   
\end{equation}
\end{lemma}     
\noindent   The proof   is a lengthier  technical  exercise  recapitulated  in Appendix \ref{matrix}. \qed  
  \begin{lemma}\label{driftaddition} 
\begin{equation} \label{slutsumma2}
 \sum_{l=1}^{M_{i}} d_{kl}^{(i)} \widetilde{h}_{i}(l)= \sum_{l=1}^{M_{i}-1} d_{kl}^{(i)} \widetilde{h}_{i}(l) -  x^{(i)}_{k}   
 x^{(i)}_{M_{i}}\widetilde{h}_{i}\left(M_{i}\right). 
\end{equation} 
 \end{lemma}   
\noindent  {\em Proof:}
\begin{eqnarray}  
 \sum_{l=1}^{M_{i}} d_{kl}^{(i)} \widetilde{h}_{i}(l) & = & 
 \sum_{l=1}^{M_{i}-1} d_{kl}^{(i)} \widetilde{h}_{i}(l) +  d_{kM_{i}}^{(i)}  \widetilde{h}_{i}\left(M_{i}\right) \nonumber \\ 
 & &  \\ 
 & = &  \sum_{l=1}^{M_{i}-1} d_{kl}^{(i)} \widetilde{h}_{i}(l) -  x^{(i)}_{k}   
 x^{(i)}_{M_{i}}\widetilde{h}_{i}\left(M_{i}\right). \nonumber 
\end{eqnarray}  
 \qed \\
The  following identity holds for any  function $V({\bf x})$ that has the required partial derivatives.  
\begin{lemma}\label{driftaddition2}
\begin{equation}\label{driftalt}
\sum_{l=1}^{M_i-1}d_{kl}^{(i)} ({\bf x})    
   V^{'}_{ x_{l}^{(i)}} \left({\bf x}\right) = x^{(i)}_{k} \left[  V^{'}_{ x_{k}^{(i)}} \left({\bf x}\right)
- \sum_{l=1 }^{M_i-1} x_{l}^{(i)}      
   V^{'}_{ x_{l}^{(i)}} \left({\bf x}\right)  \right].  
\end{equation} 
\end{lemma}  
{\em Proof}:$$
\sum_{l=1}^{M_i-1}d_{kl}^{(i)} ({\bf x})    
   V^{'}_{ x_{l}^{(i)}} \left({\bf x}\right) = 
  -x^{(i)}_{k} \sum_{l=1, l\neq k }^{M_{i}-1} x_{l}^{(i)}      
   V^{'}_{ x_{l}^{(i)}} \left({\bf x}\right)  + (1-x^{(i)}_{k}) x^{(i)}_{k}  V^{'}_{ x_{k}^{(i)}} \left({\bf x}\right)  
$$
$$
=  -x^{(i)}_{k} \sum_{l=1,   }^{M_i-1} x_{l}^{(i)}      
   V^{'}_{ x_{l}^{(i)}} \left({\bf x}\right)  +  x^{(i)}_{k}  V^{'}_{ x_{k}^{(i)}} \left({\bf x}\right).  
$$
\qed \\  
We construct explicitly the  potential $V({\bf x})$ in Svirezhev-Shahshahani gradient form  ${\bf D}({\bf x})  \nabla_{{\bf x}} V \left({\bf x}\right)$.   
  

 \subsection{The Potential $V({\bf x})$}\label{potshasvhi}

Let us first  expand  our basic  notation in (\ref{vekt1}) (now   a column vector)  as follows. We introduce the $ M_i  \times 1$ vector    
 \begin{equation}\label{frekekvloc} 
{\bf x}^{(i)} =\left( \begin{array}{c}  x^{(i)} \\    x_{M_i}^{(i)} \end{array} \right)  = 
\left(  \begin{array}{c} x_{1}^{(i)} \\ \vdots  \\ x_{M_i-1}^{(i)}  \\  x_{M_i}^{(i)} \end{array}  \right).  
\end{equation} 
  We stack these vectors to   the $ \sum_{i=1}^{L} M_{i}  \times 1 $   vector   (written in the  transposed ($^{T}$) form  for economy  of space)  
 \begin{eqnarray}\label{frekekvv} 
{\bf \underline{x}}^{T}&= &({\bf x}^{(1)}, {\bf x}^{(2)},\cdots, {\bf x}^{(L)}) \\ \nonumber  
& & \\ \nonumber 
& =&  (x_{1}^{(1)},\ldots,x_{M_1}^{(1)},x_{1}^{(2)},\ldots,x_{M_2}^{(2)},\ldots, x_{1}^{(L)},\ldots,x_{M_L}^{(L)}). \nonumber  
\end{eqnarray}  
Let ${\bf h}$ be the $ \sum_{i=1}^{L} M_{i}   \times 1 $   vector  of one locus   selection parameters 
$$
{\bf h}^{T} = \left(   h_1(1) ,  h_1(2) ,  \cdots ,  h_1(M_{1})    \cdots   h_L(1),   h_L(2),   \cdots   h_L(M_{i})  \right).   
$$ 
The next goal is to define a  $ \sum_{i=1}^{L} M_{i}  \times \sum_{i=1}^{L} M_{i}$ matrix   $A$ so that the quadratic form 
$$  
{\bf \underline{x}}^{T}A {\bf \underline{x}} 
$$
can be used to define the desired potential $V({\bf x}) $.   
  
Let first ${\bf 0}^{(i)}$ denote the $M_{i} \times M_{i}$  matrix of zeroes for $i=1,2,\ldots, L$ (not to be confused with the zero matrices  of  other dimensions in (\ref{bigd})).  
These  matrices  ${\bf 0}^{(i)}$ are inserted as   block matrices  in the main diagonal of $A$, i.e.,  
\begin{equation}\label{potmatrix}  
A \stackrel{\rm def}{=} \left( \begin{array}{cccccccc}
{\bf 0}^{(1)} &  {\bf J}_{1}(M_{2}) &  {\bf J}_{1}(M_{3}) & \ldots & \ldots  & \ldots  &\ldots & {\bf J}_{1}(M_{L})   \\ 
{\bf J}_{2}(M_{1}) &    {\bf 0}^{(2)} &   {\bf J}_{2}(M_{3}) & \ldots & \ldots  &  \ldots  & \ldots &  {\bf J}_{2}(M_{L})  \\   
\vdots &  \vdots &  \vdots  &  \ddots  &\vdots  & \ldots  & \ldots &  \vdots      \\ 
{\bf J}_{i}(M_{1}) &  \ldots & \ldots &  {\bf J}_{i}(M_{i-1}) & {\bf 0}^{(i)} &   {\bf J}_{i}(M_{i+1})& \ldots 
&  {\bf J}_{i}(M_{L}) \\  

\vdots &  \vdots &  \vdots  &  \vdots  &\vdots & \ddots  &    \ldots    &  \vdots      \\ 

{\bf J}_{L}(M_{1})   & {\bf J}_{L}(M_{2}) &   {\bf J}_{L}(M_{3}) & \ldots & \ldots  & \ldots  & \ldots    & {\bf 0}^{(L)} \end{array} \right).  
\end{equation}   
 Here   ${\bf J}_{i}(M_{l})$ is  a   block matrix of dimension  $M_{i} \times M_{l}$. 
 It is given as 
 \begin{equation}\label{potmatrix1}  
{\bf J}_{i}(M_{l}) = \left( \begin{array}{cccc}

 J_{il}(1,1) &  J_{il}(1,2) &  \ldots  &    J_{il}(1,M_{l}) \\
  J_{il}(2,1), &  J_{il}(2,2) &  \ldots  & J_{il}(2,M_{l}) \\
  
  \vdots & \vdots &  \vdots &   \vdots \\

 J_{il}(M_{i},1) &  J_{il}(M_{i},2) &  \ldots &  J_{il}(M_{i},M_{l}) \end{array} \right) 
   \end{equation}  
 by means of the two locus selection parameters at locus $i$.  By the  symmetry assumption (\ref{symmetri})  
 ${\bf J}_{i}(M_{l})$ is a symmetric matrix.  
 Thus,   a  generic   $1 \times \sum_{l=1}^{L}M_{i}$  row  in   $A$  looks like 
\begin{equation}\label{potmatrix2}                                                                                                                                                                                                                                                                                                                                                                                                                                                                                                                                                                                                                                                                                                                                                                                                                                                                                                                                                                                                                                                                                                                                                                                                                                                                                                                                                                                                                                                                                                                                                                                                                                                                                  
{\small  J_{i1}(k,1), \ldots   J_{i1}(k,M_{1}),  J_{i2}(k,1), \ldots   J_{i2}(k,M_{2}), \ldots,  {\bf 0}_{k}^{(i)},
 \ldots,  J_{iL}(k,1), \ldots   J_{iL}(k,M_{L}) )  },   
\end{equation}  
where now  ${\bf 0}_{k}^{(i)}$  is the $k$th row of in ${\bf 0}^{(i)}$ with $M_{i}$ zeroes.    
We note also that  ${\bf J}_{i}(M_{l})$ and ${\bf x}^{(l)}$ in (\ref{frekekvloc}) are compatible for the matrix  multiplication 
$ {\bf J}_{i}(M_{l}){\bf x}^{(l)}$.   

We observe  that  the  elements   in the   $ \sum_{i=1}^{L} M_{i}   \times 1 $   vector 
$A {\bf \underline{x}}$ are by the construction above for all cases of  $(i,k)$  nothing but the expressions  
\begin{equation}\label{potmatrix25}                                                                                                                                                                                                                                                                                                                                                                                                                                                                                                                                                                                                                                                                                                                                                                                                                                                                                                                                                                                                                                                                                                                                                                                                                                                                                                                                                                                                                                                                                                                                                                                                                                                                                  
\sum_{\underset{r \neq i}{r=1}}^{L}\sum_{m=1}^{M_r}  J_{ir}(k,m)x_{m}^{(r)}.    
\end{equation} 
By the  symmetry  of   ${\bf J}_{i}(M_{l})$s     the matrix $A$ is 
a symmetric matrix.  We set 
\begin{equation}\label{potmatrix3} 
   W({\bf \underline{x}}) \stackrel{\rm def}{=}  {\bf \underline{x}}^{T} {\bf h} + \frac{1}{2} {\bf \underline{x}}^{T}A {\bf \underline{x}}.  
\end{equation} 
We note that  $W({\bf \underline{x}})$ is a function of $\sum_{l=1}^{L}M_{i}$ variables, and that 
a partial derivative  like $ \frac{\partial}{\partial  x_{M_i}^{(i)} } W({\bf \underline{x}})$ means differentiation w.r.t. the appropriate position in   ${\bf \underline{x}}$.    
     It turns out that 
   \begin{equation}\label{Vdef} 
  V({\bf x}):=     W({\bf \underline{x}})
    \end{equation}         
is the function sought for the  Svirezhev-Shahshahani gradient form.

 \subsection{The Svirezhev-Shahshahani gradient }  
\begin{lemma} \label{lemmasvishas} 
\begin{equation}\label{flowlocusallellepart303} 
   \sum_{l=1}^{M_{i}} d_{kl}^{(i)} \widetilde{h}_{i}(l) =    \sum_{l=1}^{M_i-1}d_{kl}^{(i)}      
   V^{'}_{ x_{l}^{(i)}} \left({\bf x}\right)   
 \end{equation}
 \end{lemma}    
\noindent  {\em Proof}:  We prove the assertion of the lemma by expanding the right hand side of 
(\ref{flowlocusallellepart303}).  The pertinent partial derivatives  are  for $l=1, \ldots, M_{i}-1$ 
$$
 V^{'}_{ x_{l}^{(i)}} \left({\bf x}\right)= 
\frac{\partial}{\partial  x_{l}^{(i)} } W({\bf \underline{x}}) -  \frac{\partial}{\partial  x_{M_i}^{(i)} } W({\bf \underline{x}}), 
$$
since the derivative of the  inner function is  $x_{M_i}^{(i)}= 1-  \sum_{k=1}^{M_{i}-1}x_{l}^{(i)}$ in 
$W$  w.r.t 
$ x_{l}^{(i)}$ equals $-1$. 
We note that 
$$
\frac{\partial}{\partial  x_{l}^{(i)} }  {\bf \underline{x}}^{T} {\bf h} - \frac{\partial}{\partial   x_{M_i}^{(i)} }  {\bf \underline{x}}^{T} {\bf h}
 =  h_{i}(l) - h_{i}(M_{i}).  
$$
Next,   since $A$ is a symmetric matrix,  
$$ 
\nabla_{\bf \underline{x}}  \left[ \frac{1}{2} {\bf \underline{x}}^{T}A {\bf \underline{x}} \right]=  
 A {\bf \underline{x}}.  
$$ 
Hence we obtain   by (\ref{potmatrix25})  that 
$$
\frac{\partial}{\partial  x_{l}^{(i)} } W({\bf \underline{x}}) =  \sum_{\underset{r \neq i}{r=1}}^{L}\sum_{m=1}^{M_r}  J_{ir}(l,m)x_{m}^{(r)} 
$$ 
and 
$$
\frac{\partial}{\partial  x_{M_{i}}^{(i)} } W({\bf \underline{x}})=  
\sum_{\underset{r \neq i}{r=1}}^{L}\sum_{m=1}^{M_r}  J_{ir}(M_{i},m)x_{m}^{(r)}. 
$$
Hence 
$$
 V^{'}_{ x_{l}^{(i)}} \left({\bf x}\right) = \frac{\partial}{\partial  x_{l}^{(i)} } W({\bf \underline{x}}) -  \frac{\partial}{\partial  x_{M_i}^{(i)} } W({\bf \underline{x}})  
$$ 
$$ 
= h_{i}(l) - h_{i}(M_{i}) +  \sum_{\underset{r \neq i}{r=1}}^{L}\sum_{m=1}^{M_r}  J_{ir}(k,m)x_{m}^{(r)}  - 
 \sum_{\underset{r \neq i}{r=1}}^{L}\sum_{m=1}^{M_r}  J_{ir}(M_{i},m)x_{m}^{(r)}. 
$$
In view of  (\ref{hejp2})  we have thus  shown that       
\begin{equation}\label{flowlocusallellepart505} 
  V^{'}_{ x_{l}^{(i)}} \left({\bf x}\right)=  \widetilde{h}_{i}(l)   -  \widetilde{h}_{i}\left(M_{i}\right).   
\end{equation}  
 Hence   
\begin{equation}\label{flowlocusallellepart506} 
 \sum_{l=1}^{M_i-1}d_{kl}^{(i)}      
   V^{'}_{ x_{l}^{(i)}} \left({\bf x}\right)  =   
   \sum_{l=1}^{M_{i}-1 } d_{kl}^{(i)} \widetilde{h}_{i}(l) -   \sum_{l=1}^{M_{i}-1} d_{kl}^{(i)}  \widetilde{h}_{i}\left(M_{i}\right).  
 \end{equation}
 The last term is   
 $$
 \sum_{l=1}^{M_{i}-1} d_{kl}^{(i)}  \widetilde{h}_{i}\left(M_{i}\right) = 
 \widetilde{h}_{i}\left(M_{i}\right)  \sum_{l=1}^{M_{i}-1} d_{kl}^{(i)}      
 $$ 
 $$
=  \widetilde{h}_{i}\left(M_{i}\right)\left[ 
   \sum_{l=1, l \neq k }^{M_{i}-1 } (-x^{(i)}_{k}x^{(i)}_{l})   +   x^{(i)}_{k} - (x^{(i)}_{k})^{2}\right]  
  $$
$$
=  \widetilde{h}_{i}\left(M_{i}\right)  x^{(i)}_{k} \left[ 
  (-1)\sum_{l=1, l \neq k }^{M_{i}-1}x^{(i)}_{l}   +  1  - x^{(i)}_{k}\right] 
$$
$$
= \widetilde{h}_{i}\left(M_{i}\right)  x^{(i)}_{k} \left[ 
  (-1)( 1 -x^{(i)}_{k} -  x^{(i)}_{M_{i}})   +  1  - x^{(i)}_{k}\right] 
$$
$$
= \widetilde{h}_{i}\left(M_{i}\right)  x^{(i)}_{k} \left[  -1 +  x^{(i)}_{k} +  x^{(i)}_{M_{i}} +1  - x^{(i)}_{k}\right] =  x^{(i)}_{k}   x^{(i)}_{M_{i}} \widetilde{h}_{i}\left(M_{i}\right). 
$$     
Thus we have   in (\ref{flowlocusallellepart506}) that 
\begin{eqnarray} \label{slutsumma} 
 \sum_{l=1}^{M_i-1}d_{kl}^{(i)}      
   V^{'}_{ x_{l}^{(i)}} \left({\bf x}\right) & = &   
   \sum_{l=1}^{M_{i}-1} d_{kl}^{(i)} \widetilde{h}_{i}(l) -  x^{(i)}_{k}   x^{(i)}_{M_{i}} \widetilde{h}_{i}\left(M_{i}\right)  \nonumber \\
   & & \\ 
   & = &    \sum_{l=1}^{M_{i}-1} d_{kl}^{(i)} \widetilde{h}_{i}(l) +   d^{(i)}_{k M_{i}}  \widetilde{h}_{i}\left(M_{i}\right), \nonumber  
 \end{eqnarray} 
which is the left hand side of (\ref{flowlocusallellepart303}), as claimed. \qed \\
\noindent In view of the preceding lemma \ref{ealemma4} 
\begin{equation}\label{bigd46}      
  p_{k}^{(i)}({\bf x})   = u_{k}^{(i)}-  \bar{u}      +  x_k^{(i)} \left ( \bar{m}_k^{(i)}-\bar{m}\right )      
\end{equation}
and by  lemma  \ref{ealemma3} 
\begin{equation}\label{bigd47} 
p_{k}^{(i)}({\bf x})   =   u_{k}^{(i)}-  \bar{u}   +     \sum_{l=1}^{M_{i}} d_{kl}^{(i)} \widetilde{h}_{i}(l)   
\end{equation}
 and by lemma \ref{lemmasvishas}   for the potential in (\ref{Vdef})  
\begin{equation}\label{bigd48} 
p_{k}^{(i)}({\bf x})   =   u_{k}^{(i)}-  \bar{u}   +    \sum_{l=1}^{M_i-1}d_{kl}^{(i)}      
   V^{'}_{ x_{l}^{(i)}} \left({\bf x}\right).  
\end{equation} 
By rules of matrix calculus we observe that  $\sum_{l=1}^{M_i-1}d_{kl}^{(i)}      
   V^{'}_{ x_{l}^{(i)}} \left({\bf x}\right)$ is the $k$th component in the $(M_{i}-1)\times 1$ - vector 
$$
{\bf D}^{(i)}\left({\bf x}\right)  V^{'}_{ x^{(i)}}\left({\bf x}\right).  
$$     
When  define        an $(M_{i}-1) \times 1 $  vector  
\begin{equation}\label{bigd4}
{\bf G}^{(i)}\left({\bf x}\right) \stackrel{\rm def}{=} 
 \left( \begin{array}{c} 
x^{(i)} _{1} \left( \bar{m}_1^{(i)}-\bar{m} \right) \\ 
\vdots \\ 
x^{(i)} _{M_{i}-1} \left( \bar{m}_{M_{i}-1}^{(i)}-\bar{m} \right)     
 \end{array} \right), 
\end{equation}  
we have shown  that 
$$
{\bf G}^{(i)}\left({\bf x}\right) =  {\bf D}^{(i)}\left({\bf x}\right)  V^{'}_{ x^{(i)}}\left({\bf x}\right).
$$ 
Thus  adding    
$$
{\bf G} \left({\bf x}\right) = {\bf D} \left({\bf x}\right)  V^{'}_{{\bf x}}\left({\bf x}\right) = \left( \begin{array}{c} 
 
{\bf D}^{(1)}\left({\bf x}\right) \nabla_{ x^{(1)}} V \left({\bf x}\right)   \\
 
  \vdots \\
  
{\bf D}^{(L)}\left({\bf x}\right) \nabla_{ x^{(L)}} V \left({\bf x}\right)  \end{array} \right) 
$$ 
to ${\bf \mu} \left({\bf x}\right)$ in (\ref{bigd2}), we have  an explicit expression for the drift function  in (\ref{bigdsd}).  
Next we prove the weak convergence of the sequence of interpolated Markov chains  to   (\ref{bigdsd}).

\section{The  Diffusion Approximation    } \label{limitmgle}   
 \setcounter{equation}{0} 
The studies  summarized  in this section   prove  the  
  weak convergence of a sequence of the Markov chains defined in section  \ref{transproball}  to a diffusion process and consist of straigthforward  verifications  the conditions  for weak convergence found  in   \cite[Theorem 7.1]{Durrett1996}, \cite[Lemma 4.1.]{sato1976diffusion}, see also  \cite[ch.10 thm 3.5]{ethier2009markov}. 
  
These conditions correspond one-to-one  to   the technical lemmas in Appendix D, and are given as  (\ref{driftlimit}), (\ref{difflimit}),  (\ref{gdr_finished1})  and (\ref{mixed}).  The level of mathematical effort  herewith  is merely to verify that these  conditions 
are valid in the current  situation, not to contribute to a general advancement   of diffusion approximation.   
 
  We    recall  variables (\ref{scale1})  scaled by  the  inverse population size, or,  more conveniently:    
\begin{antag}\label{scaling}
{\rm For all $l,k$ and all $i$   
\begin{equation}\label{popscale1} 
\upsilon_{lk}^{(i)} = \frac{u_{lk}^{(i)} }{N},  
\end{equation}
  \begin{equation}\label{pairwise_selection11}
 \bar{v}  = 1 +  \frac{\bar{m} }{N},  
\end{equation}
and 
\begin{equation}\label{pairwise_selection12}
 \bar{v}_k^{(i)}  = 1 +  \frac{ \bar{m}_k^{(i)}}{N},   
\end{equation}  
where  $\bar{m} $  is    given in (\ref{pairwise_selection}) and
 $\bar{m}_k^{(i)}$ is   in  (\ref{pairwise_selection3}). }
 \end{antag} \qed \\ 
  We normalize in (\ref{onelinterpol}) to  get  
\begin{equation}\label{normvariab2}
X^{(N)}_{i,k}(n) =  \frac{1}{N}Y_{i,k}^{(N)} \left( \frac{n}{N} \right)   
\end{equation}
and 
\begin{equation}\label{onelinterpol22}
X^{(N)}_{(i)}(n) = \left( X_{i,1}^{(N)}(n), X_{i,2}^{(N)}(n), \ldots, X_{i,M_{i}-1}^{(N)}(n)   \right). 
\end{equation}
We have   thus also  
\begin{equation}\label{onelinterpol2223}
X^{(N)}(n) = \left\{ X^{(N)}_{(i)}(n) \right \}_{ 1 \leq i  \leq M_{i}},
\end{equation}
\begin{equation}\label{onelinterpol2224} 
 X^{(N)}  = \left\{ X^{(N)} (n)   \right \}_{ n \geq 0}.
\end{equation} 
 \begin{equation}\label{normvariab21}
X^{(N)}_{i,k}(n) =  \frac{1}{N}Y_{i,k}^{(N)} \left( \frac{n}{N} \right)   
\end{equation}
 and 
\begin{equation}\label{onelinterpol22245} 
 X^{(N)}  = \left\{ X^{(N)} (n)   \right \}_{ n \geq 0}.
\end{equation}
 Let us  next define   the process $ \{ X^{(N)}(t) \mid   t \geq  0    \}$   as   the piece-wise constant continuous time   interpolation of the process $X^{(N)}$   in (\ref{onelinterpol2224})  
$$
X^{(N)}(t)  =  X^{(N)} \left( \frac{1}{N}\lfloor  tN  \rfloor    \right), 0 \leq t   < +\infty,   
$$
where $\lfloor x \rfloor $  is the integer part of a real number $x$.  
Let now for every $i \in \{ 1, \ldots, L\}$
\begin{equation}\label{satodomain2} 
{\cal D}^{(i)} \stackrel{\rm def }{=} \left \{ f|_{{\bf K}_{i}} \mid f \in C^{4}(\mathbf{R}^{\sum_{i=1}^{L} M_{i}}) \quad  \text{ and} \quad  \frac{\partial}{\partial x_{i}} \left( f|_{{\bf K}_{i}} \right) =  \left(\frac{\partial}{\partial x_{i}} f \right)|_{{\bf K}_{i}} \right\}
\end{equation}
 ${\cal D}_{i}$ contains the restrictions of four times  differentiable  real valued  functions $f$ on  $\mathbf{R}^{\sum_{i=1}^{L} M_{i}})$ to ${\bf K}$  such that 
the restriction of  a partial derivative of $f$  to  ${\bf K}_{i}$, $\left(\frac{\partial}{\partial x_{i}} f \right)|_{{\bf K}_{i}} $,  equals the same partial derivative of the restriction   $f|_{{\bf K}_{i}}$.     
 The requirement of  $C^{4}$ -functions in the domain  is implicitly  needed in the uniqueness 
part of the proof of the next proposition and is due  to  \cite{ethier1976class}. 
  The differential generator ${\cal L}$ of the desired limiting 
process is  defined  by    
\begin{equation}\label{mgproblem}
\begin{aligned}
{\cal L} =  \sum_{i=1}^{L}\sum_{k=1}^{M_i-1} p_{k}^{(i)}({\bf x}) \frac{\partial }{\partial x_{k}^{(i)}}  \\
&+\frac{1}{2}\sum_{i=1}^{L} \sum_{k=1}^{M_i-1}\sum_{l=1}^{M_i-1} \left[d_{kl}^{(i)}({\bf x}) \right]\frac{\partial^2}{\partial x_{k}^{(i)}\partial x_{l}^{(i)}}  
\end{aligned}
\end{equation} 
with the domain  given from   (\ref{satodomain2}) by    
$$
{\cal D} =  \times_{i=1}^{L} {\cal D}^{(i)}. 
$$

\begin{prop}\label{diffapprthm}  
  If  the assumptions   \ref{cindep}, (\ref{indmut223})  and (\ref{pos223})  hold,   and  if $X^{(N)}(0) \rightarrow x_{0}$, as   $N  \rightarrow +\infty$,   
 then  the process 
$ \{ X^{(N)}(t) \mid   t \geq  0    \}$ converges under  the scalings in Assumption \ref{scaling} weakly as   $N  \rightarrow +\infty$ to $ \{ X(t) \mid t \geq  0   \}$, which is the unique solution of  the  martingale problem for 
${\cal L}$  (\ref{mgproblem}). 
\end{prop}   
 \noindent{\em Proof}:

The proof is in three steps. The first, {\bf Step a) }  identifies the drift $\mu ({\bf x})+ {\bf D} \left({\bf x}\right)  V^{'}_{{\bf x}}\left({\bf x}\right)  $  and diffusion ${\bf D} \left({\bf x}\right) $ of a  limiting 
process.  It remains to  prove the uniqueness of the martingale problem defined by this  drift
 and diffusion. If  the Svirezhev-Shahshahani gradient form is removed, we are faced with   independent and uncoupled multiallelic  processes, where  the uniqueness result of \cite{ethier1976class} can be applied,as shown in {\bf Step b)}.  
In step   {\bf Step c) } we use the finding in  {\bf Step b)} by an adding  the Svirezhev-Shahshahani gradient form. We  can then  apply the Girsanov theorem on trnafomraions of drift and measure  to prove the desired 
uniqueness.   
   
\begin{description}  
 
 \item[Step a)]    The  convergences, which are uniform in     
${\bf x}\in \times_{i=1}^{L}{\bf  K}_{i}$,  in   (\ref{driftlimit}), (\ref{difflimit}),  (\ref{gdr_finished1})  and (\ref{mixed})  with    the continuous  limiting drift  $\mu ({\bf x}) + {\bf D} \left({\bf x}\right)  V^{'}_{{\bf x}}\left({\bf x}\right)$ and  
and  the diffusion  (genetic drift)   (\ref{genetic_drift})  as well as (\ref{onelinterpol28}),   for every locus  and for every allele type at  every locus, have been shown  in the Appendix D.  

By these facts the assertion in the proposition follows by, see e.g.,   \cite[Theorem 7.1]{Durrett1996},  and      \cite[Lemma 4.1.]{sato1976diffusion}, if   the  martingale problem associated   to the differential operator ${\cal L}$, or ${\rm MP}( {\bf \mu} + {\bf D}  \nabla  V , D)$ (c.f.  (\ref{bigdsd})) is unique. We shall now establish  the desired uniqueness by first  proving uniqueness for an uncoupled system of SDE's.

\item[Step b) No selection]  Now we consider at  every locus $i$ the multivariate diffusion $Y^{i} = \{ Y^{(i)}(t) =  \left(  Y^{(i)}_{1}(t), \ldots,  Y^{(i)}_{M_{i}-1}(t)  \right) \mid t \geq 0\}$ with values in   $  {\bf K}_{i}$ 
and with the notation (\ref{bigd2})  satisfying  
\begin{equation}\label{ethierone}  
  d Y^{(i)}(t) = \mu^{(i)}(Y^{(i)}(t) )dt   + D^{(i)}(Y^{(i)}(t))^{1/2}dW^{(i)}(t).  
 \end{equation}   
  For diffusions with  values in $  {\bf K}_{i}$  and with diffusion matrix  $D^{(i)}({\bf x})$, it holds, here we need $C^{4}$ in the  domain ${\cal  D}$,  by   \cite{ethier1976class}, or,  \cite[p.991]{campiti2019binomial},   \cite[p.134 and p. 135 Corollary 2.1]{shiga1981diffusion}, that 
   if  for every  ${\bf y} \in {\bf K}_{i}$    
\begin{equation}\label{ethierentydig1} 
 \sum_{k=1}^{M_i-1} \mu_{k}^{(i)}({\bf y})  \leq 0,     \mbox{ with  $\sum_{k=1}^{M_{i}-1} y^{(i)}_{k}=1$} 
\end{equation} 
and for every $k$ 
\begin{equation}\label{ethierentydig2}  
 \mu_{k}^{(i)}({\bf y}) \geq 0  \mbox{ if  $  y^{(i)}_{k}=0$,}  
\end{equation} 
then the martingale problem corresponding to  $\left(Y^{(i)},   \mu ^{(i)}({\bf y}), D^{i} \right)$  has a unique solution.  We check next that these conditions hold.   First, if  $\sum_{k=1}^{M_{i}-1} y^{(i)}_{k}=1$, then  
$$
 \sum_{k=1}^{M_i-1} \mu_{k}^{(Y_i)}({\bf y}) =  \sum_{k=1}^{M_i-1} \left( u^{(i)}_{k}  - \bar{u}^{(i)}  y^{(i)}_{k} \right)  
$$ 
$$
= \sum_{k=1}^{M_i-1}  u^{(i)}_{k}  - \bar{u}^{(i)} \sum_{k=1}^{M_i-1}  y^{(i)}_{k} = \sum_{k=1}^{M_i-1}  u^{(i)}_{k}  - \bar{u}^{(i)} = - u^{(i)}_{M_{i}} \leq 0, 
$$ 
and (\ref{ethierentydig1}) is checked.   For  (\ref{ethierentydig2}) we observe that if  $y^{(i)}_{k} =0$, then  
$$
\mu_{k}^{(i)}({\bf y}) =  u^{(i)}_{k}  - \bar{u}^{(i)}  y^{(i)}_{k} =  u^{(i)}_{k} > 0. 
$$  
 Hence, the martingale problem  ${\rm MP} \left(  \mu ^{(i)}({\bf y}), {\bf D}^{(i)} \right)$  has a unique solution.  By (\ref{bigd2})  we have $
{\bf \mu}({\bf y})= \left(    \mu^{(1)}({\bf y}), \ldots,  \mu^{(L)}({\bf y})\right)$.    
As the processes  $Y^{(i)}$  are independent and uncoupled, it follows that  the martingale problem corresponding to the diffusion 
$$
 {\rm MP} \left( {\bf \mu}({\bf y}),  {\bf D} \right)  
$$  
and represented by  the system of locus-wise decoupled stochastic differential equations  
\begin{equation}\label{bigd44} 
  dY(t) =  {\bf \mu}(Y(t))dt +    {\bf D}^{1/2}\left(Y(t)\right)d {\bf W}(t), 
 \end{equation} 
has a unique solution ${\bf P}$. 
\item[Step c)Change of drift and change of measure]  
  The equation (\ref{bigd4}) 
 is clearly obtained from (\ref{bigd44}) by a   change of drift.  Let now  
 ${\bf P}$ the unique probability measure such that $X$ satisfies 
 \begin{equation}\label{bigd445} 
  dX(t) =  {\bf \mu}(X(t))dt +    {\bf D}^{1/2}\left(X(t)\right)d {\bf B}(t), 
 \end{equation}   
w.r.t. a ${\bf P}$ -Wiener process ${\bf B}=\{{\bf B}^{(i)} \}_{i=1}^{L}$.  
Then we recall (\ref{bigdsd}), or, 
$$ 
  dX(t) =  {\bf \mu}(X(t))dt + {\bf D}(X(t))  \nabla_{{\bf x}} V \left(X(t)\right)dt  +  {\bf D}^{1/2}\left(X(t)\right)d {\bf W}(t).  
$$ 
  We set for any $i$  and $t>0$   
\begin{equation}\label{bigd356} 
H (t)  = (-1) \cdot   {\bf D}^{1/2}\left(X (t) \right) \nabla_{ {\bf x}} V \left(X(t)\right) .  
\end{equation}    
The expression  ${\bf D}^{1/2}\left({\bf x}\right) \nabla_{{\bf x}} V \left({\bf x}\right)$ is  well defined 
on all of  ${\bf K}$
and its vector norm is uniformly bounded  there. Then,  by \cite[Theorem 6.4.7, p. 153]{kallianpur2014stochastic}           
 we can thus define ${\bf Q}$ by the Radon-Nikodym derivative 
 $$
 \frac{d{\bf Q}}{d{\bf P}} = e^{  \int_{0}^{S} H (t)^{T} d {\bf B}(t) -\frac{1}{2}\int_{0}^{S} H (t)^{T}H (t)dt }. 
 $$
Then the process     $\left\{{\bf W}(t)  \right\}$  given by  
$$
{\bf W} (t)= {\bf B} (t)  + \int_{0}^{t} H (s)ds  
$$ 
is,  by Girsanov$^{'}$s theorem,    Wiener processes w.r.t. the measure  ${\bf Q}$, with the same 
covariances (quadratic variations) as ${\bf B} (t)$. 
Then,  by  matrix multiplication 
$$
 {\bf D} ^{1/2}\left(X(t)\right)  d{\bf W} (t)=  {\bf D} ^{1/2}\left(X(t)\right)d{\bf B} (t)  +   {\bf D} ^{1/2}\left(X(t)\right) H (t) dt  
$$
$$
=  {\bf D}^{1/2}\left(X(t)\right)d{\bf B} (t)  -  {\bf D}^{1/2}\left(X(t)\right)   {\bf D}^{1/2}\left(X (t)\right) \nabla_{ {\bf x}} V \left(X(t)\right)dt  
$$
$$
= {\bf D}^{1/2}\left(X(t)\right)d{\bf B}(t)  -  {\bf D}\left(X(t)\right) \nabla_{{\bf x} } V \left(X(t)\right)dt.   
$$
When we rearrange this, we get    
$$
{\bf \mu} (X(t))dt +   {\bf D} \left(X(t)\right) \nabla_{ {\bf x}} V \left(X(t)\right)dt  + {\bf D}^{1/2}\left(X(t)\right)d{\bf W}(t)  
$$
$$
=  {\bf \mu} (X(t))dt +   {\bf D}^{1/2}\left(X(t)\right)d{\bf B} (t). 
$$
Since  $dX (t)= {\bf \mu} (X(t))dt +   {\bf D}^{1/2}\left(X(t)\right)d{\bf B}(t)$ via  the canonical 
path space, we have   
$$
dX (t)= {\bf \mu} (X(t))dt +   {\bf D} \left(X(t)\right) \nabla_{ {\bf x}} V \left(X(t)\right)dt  + {\bf D}^{1/2}\left(X(t)\right)d{\bf W}(t). 
$$ 
 
Hence, by theorem (5.2) of \cite[pp.204$-$205]{Durrett1996}  there is a 1-1 correspondence  between 
the solution ${\bf P}$ of the  martingale problems $ {\rm MP} \left(  \mu, {\bf D}\right)$ and 
the solution ${\bf Q}$ of   $  {\rm MP} \left(   {\bf \mu}  +  {\bf D}  \nabla_{ {\bf x }} V ,  {\bf D} \right)$.  Hence    ${\rm MP} \left(   {\bf \mu}  + {\bf D}  \nabla_{ {\bf x }} V ,  {\bf D} \right)$ has a unique solution.    \qed

\end{description}

 \section{The  Stationary Probability Distribution: The Explicit Solution}\label{solutionstmk}
\setcounter{equation}{0} 
\subsection{The  Fokker -Planck Equation and  the Probability Flow }
  Given pre-specified one- and two-locus selection parameters $h_i(k)$ and $J_{ij}(k,l)$ and mutation intensities   $u_{lk}^{(i)}$, the probability density function   $P({\bf x},t)$ of the    diffusion  $ \{ X(t) \mid  t  \geq 0   \}$ in (\ref{bigd55})     occupying state ${\bf x}$ at time $t$ is  governed by the multidimensional Fokker-Planck (or Kolmogorov Forward) equation
\begin{equation}\label{final_Fokker_Planck} 
\begin{aligned}
\frac{\partial P({\bf x},t)}{\partial t} &= -\sum_{i=1}^{L}\sum_{k=1}^{M_i-1}\frac{\partial }{\partial x_{k}^{(i)}} \left[p_{k}^{(i)}({\bf x})P({\bf x},t) \right ]\\
&+\frac{1}{2}\sum_{i=1}^{L}\sum_{k=1}^{M_i-1}\sum_{l=1}^{M_i-1}\frac{\partial^2}{\partial x_{k}^{(i)}\partial x_{l}^{(i)}}\left[d_{kl}^{(i)}({\bf x})P({\bf x},t)\right]. 
\end{aligned}
\end{equation} 

In order to solve  this   stationary Fokker-Planck equation   we are going to study  the  probability flow \cite[pp. 133$-$134]{risken1989fokker}. The reference cited   is, however, not the source of the technical details below.
$$
J =  \left(J^{(1)}, \ldots , J^{(L)}    \right),  
$$
where 
$$
J^{(i)} = \left(J_{1}^{(i)}, \ldots , J_{M_{i}-1}^{(i)}    \right).    
$$  
   Here  
\begin{equation}\label{flowlocusallelle} 
J_{k}^{(i)}  \stackrel{\rm def}{=}  p_{k}^{(i)}({\bf x})P - \frac{1}{2}\sum_{l=1}^{M_i-1}\frac{\partial }{ \partial x_{l}^{(i)}}\left[d_{kl}^{(i)}({\bf x})P \right].
\end{equation} 
Let now 
\begin{equation}\label{flownabla1} 
\nabla J  =  \left(     \nabla_{x^{(1)}} J^{(1)} , 
\ldots, 
  \nabla_{x^{(L)}} J^{(L)} \right),
\end{equation}
and 
\begin{equation}\label{flownabla2}
  \nabla_{x^{(i)}}     J^{(i)} =  \left(     \frac{\partial} {\partial x_{1}^{(i)}} J^{(i)} ,  
\ldots  , 
 \frac{\partial} {\partial x_{M_{i}-1}^{(i)}} J^{(i)}   \right).  
\end{equation} 
Then    (\ref{final_Fokker_Planck}) can be written as 
$$ 
\frac{\partial }{\partial t} P  = -  \nabla \cdot  J.   
$$ 
Then  the   solution to the equation  (\ref{final_Fokker_Planck})   that satisfies, when  existing,    
\begin{equation}\label{nablaflow} 
 -  \nabla \cdot  J=0 
\end{equation}
is denoted by $P\left( {\bf x} \right)$ and is called the stationary or invariant probability density   w.r.t. to the measure $d {\bf x}$ on ${\bf K}$.   
We shall next solve  the stationary  Fokker-Planck equation.  In fact, we shall find a solution under  a  stronger   condition on the probability flow, namely, 
\begin{equation}\label{noflow}
J = {\bf 0}  \quad \text{ on $ \times_{i=1}^{L} {\bf K}_i $}. 
\end{equation}
The possibility of explicit solution rests   upon     the presence of the Svirezhev-Shahshahani gradient form.

\subsection{The Explicit Solution  } 
 \subsubsection{Auxiliaries} 
 
We set 
 \begin{equation}
\label{stationary_dist11}
\pi_{i} \left(  x^{(i)}  \right) = \left[ \left( x_{M_{i}}^{(i)}\right)^{2u_{M_{i}}^{(i)}-1}   \prod_{l=1}^{M_{i} -1} (x_l^{(i)})^{2u_l^{(i)}-1} \right]  d {\bf x}.  
   \end{equation}
This is the non-normalized  density of a Dirichlet distribution on ${\bf K}_{i}$, sometimes denoted as 
${\cal D}\left( u_1^{(i)}, \ldots, u_{M_{i} -1}^{(i)}; u_{M_{i}}^{(i)} \right)$.      
Due to the   functional relationship between the $M_{i}$  variables, their joint
probability distribution is degenerate, and the density is only for the $M_{i}-1$ variables inside $ x^{(i)}$.  Of course, the numbering  of the $M_{i}$ alleles is  arbitrary,  but so is  $\pi_{i} \left(  x^{(i)}  \right)$ invariant w.r.t. permutations.
\begin{antag}\label{mutassump}
{\rm  We assume that for all $i$ and $l$
\begin{equation}\label{dircond}
u_l^{(i)} > 0.    
\end{equation}  }
\end{antag}\qed \\
Then we set 
 \begin{equation}
\label{stationary_dist1}
\pi({\bf x}) \stackrel{\rm def}{=}   \prod_{i=1}^{L} \pi_{i} \left(  x^{(i)}  \right).
\end{equation} 
  Let  next      $V({\bf x})$ be any sufficiently  differentiable  fitness potential. 
We  set
\begin{equation}\label{stationary_dist23}
P\left({\bf x} \right) \stackrel{\rm def}{=}  \pi({\bf x})e^{2 V({\bf x})}. 
\end{equation}     
 \begin{lemma}\label{bassats1}
 Assume  (\ref{indmut223})  and   (\ref{dircond}).    
  Then,    for all  ${\bf x} \in \times_{i=1}^{L}{\bf K}_{i}$, $k=1, \ldots, M_{i}-1$,  and $i=1, \ldots, L$, 
 we have  
 \begin{equation}\label{flowlocusallellepart13}  
 \frac{1}{2}\sum_{l=1}^{M_i-1}\frac{\partial }{ \partial x_{l}^{(i)}}\left[d_{kl}^{(i)}({\bf x})P\left({\bf x}\right) \right] =  P  \left({\bf x}\right)  \left[    u_{k}^{(i)}-  \bar{u}  +   \sum_{l=1}^{M_i-1}d_{kl}^{(i)} ({\bf x})    
   V^{'}_{ x_{l}^{(i)}} \left({\bf x}\right) \right], 
 \end{equation}  
 where  $V^{'}_{ x_{l}^{(i)}} \left({\bf x}\right) $ is the partial derivative of  $V \left({\bf x}\right)$ 
w.r.t. to  $x_{l}^{(i)}$.
\end{lemma}   
\noindent The proof is a computational exercise   recapitulated   in Appendix D.   \qed \\

\subsubsection{The potential and  the probability flow}

 \begin{prop}\label{bassats2}
 Assume   (\ref{indmut223}) and (\ref{dircond}).   Let  the function $V({\bf x})$ be defined by  (\ref{potmatrix3})  as  
\begin{equation} \label{choiceofv}
V({\bf x}) \stackrel{\rm def}{=}  W({\bf \underline{x}}).          
\end{equation}  
 Then,    for all  ${\bf x} \in \times_{i=1}^{L}{\bf K}_{i}$,  
 \begin{equation}\label{stationary_dist24}
P\left({\bf x} \right) = \pi({\bf x})e^{2 V({\bf x})}  
\end{equation}  
 solves the equation of zero probability flow 
\begin{equation}\label{flowlocusallelle22} 
J_{k}^{(i)}  =   p_{k}^{(i)}({\bf x})P - \frac{1}{2}\sum_{l=1}^{M_i-1}\frac{\partial }{ \partial x_{l}^{(i)}}\left[d_{kl}^{(i)}({\bf x})P \right]=0, 
\end{equation} 
and    the stationary Fokker-Planck equation  (\ref{nablaflow}). 
\end{prop} 
{\em Proof}: We study $p_{k}^{(i)}({\bf x})P $. By virtue of   (\ref{bigd46}) - (\ref{bigd48}) it holds under (\ref{indmut223})  that  
$$  
  p_{k}^{(i)}({\bf x}) P  = \left(u_{k}^{(i)}-  \bar{u}      +  \sum_{l=1}^{M_i-1}d_{kl}^{(i)} ({\bf x})    
   V^{'}_{ x_{l}^{(i)}} \left({\bf x}\right)                \right)P.       
$$
By (\ref{flowlocusallellepart13})  we now see that    (\ref{flowlocusallelle22}) is satisfied, i.e.,    
       $J_{k}^{(i)}  = 0$  and therefore   $ -  \nabla \cdot  J=0$. \qed \\  
The normalized stationary density is again denoted by $P\left( {\bf x} \right)$, i.e., 
\begin{equation}\label{stationary_dist44} 
P\left( {\bf x} \right)=  \frac{1}{Z} \pi({\bf x})e^{2 V({\bf x})}. 
\end{equation}                  
Here  
$$
Z \stackrel{\rm def}{=} \int_{\times_{i=1}^{L} {\bf K}_i} \pi({\bf x})e^{2 V({\bf x})} d{\bf x}  
$$     
is required to exist when integrated  w.r.t. the  Lebesgue measure restricted to ${\bf K}$.

For   a single locus $L=1$, ${\bf x} = ( x_1,  x_2, \ldots, x_{M})$ with  multiple allele types as well as  mutation and   selection, Watterson \cite{watterson1977heterosis}  finds  (with a sketch of the explicit  calculations)  the stationary density   as 
\begin{equation}\label{dirichletexp} 
P \left( x_1,  x_2, \ldots, x_{M-1}, x_{M} \right)  =  \frac{1}{Z}  
  x_1^{2u_1 -1}    \cdots   x_{M }^{2u_{M }-1}  
  e^{ 2  U({\bf x}) }  dx_{1} \ldots   dx_{M-1},      
\end{equation}
where $u_{l } > 0$ for every $l$.  Watterson states  also a method of computation of       
   the constant  $Z$ in  one special case.   In spite of the obvious similarity with (\ref{stationary_dist44}), we cannot in any  straightforward  manner    regard this as  a special case (\ref{stationary_dist44}), since  the matrix $A$ in (\ref{potmatrix})   becomes for $L=1$  the $M \times M$ matrix of zeroes.  
In \cite{fearnhead2006} the loci are unlinked.

We present next an example of the computation of a stationary density by  the techniques   above. 
    
 \begin{exempel}\label{exempel3}[{\bf Two  loci, two  alleles with   selection   and mutation} ]    
{\rm  We  have $L=2$, $M_{1}=M_{2}=2$, ${\bf x} =  (x^{(1)},  x^{(2)})$. Then  we identify 
$x^{(1)}= x^{(1)}_{1},  x^{(2)}= x^{(2)}_{1}$ and have $x^{(1)}_{2}= 1-x^{(1)}_{1}$ and 
$x^{(2)}_{2}= 1-x^{(2)}_{1}$. The augmented state vector is thus 
$$
{\bf \underline{x}}^{T}= \left(  
 x^{(1)}_{1}, 
 x^{(1)}_{2},  
 x^{(2)}_{1}, 
 x^{(2)}_{2}  \right). t).  
$$  
and 
$$
{\bf \underline{x}}^{T} A {\bf \underline{x}}=  2 h  x^{(1)}_{1} x^{(2)}_{1}.  
$$
When we return to the variables $x^{(i)}$, this yields  
$$
V({\bf x}) =  \frac{1}{2} {\bf \underline{x}}^{T} A {\bf \underline{x}} =  h  x^{(1)}  x^{(2)}.  
$$
Hereafter we obtain  as in Example  \ref{karlinsvishahshahani} the following system of stochastic  differential equations 
\begin{equation} \label{2dimwfdiff}
\left \{ \begin{array}{cc}
dX^{(1)}_{t} =  u_{1}^{(1)}dt -( u_{1}^{(1)} + u_{2}^{(1)}) X^{(1)}_{t} dt  +  
h X^{(1)}_{t}(1-X^{(1)}_{t}) X^{(2)}_{t}dt + \sqrt{ X^{(1)}_{t}(1-X^{(1)}_{t}) }dW^{(1)}_{t} \\
dX^{(2)}_{t} =  u_{1}^{(2)}dt -( u_{1}^{(2)} + u_{2}^{(2)}) X^{(2)}_{t} dt  +  
h X^{(2)}_{t}(1-X^{(2)}_{t}) X^{(1)}_{t}dt + \sqrt{ X^{(2)}_{t}(1-X^{(2)}_{t}) }dW^{(2)}_{t},    
\end{array} \right. 
\end{equation} 
where  $W^{(1)}$ and $W^{(2)}$ are independent  Wiener processes. This is a system of  two coupled  Wright-Fisher 
stochastic differential equations with mutation and selection. If $h=0$,  the processes are obviously independent.      

The normalized stationary distribution (density)     is   by  
(\ref{stationary_dist24}) equal to    
\begin{equation}\label{aekdistn1} 
P \left( x^{(1)},  x^{(2)} \right) = \frac{1}{Z} \pi_{1} \left( x^{(1)}  \right) 
 \pi_{2} \left(  x^{(2)} \right) e^{2 h x^{(1)} x^{(2)} } d x^{(1)} d x^{(2)},    
\end{equation} 
 where 
$$
\pi_{1} \left( x^{(1)}  \right) =  (x^{(1)})^{2u^{(1)}_1 -1} (1-x^{(1)})^{2u^{(1)}_{2}-1} 
 $$  
 and 
   $$  
\pi_{2} \left(   x^{(2)} \right)  =  (x^{(2)})^{2u^{(2)}_1 -1} (1-x^{(2)})^{2u^{(2)}_{2}-1}.   
$$  
We can in this example determine  the normalization constant $Z$ explicitly.  First,  
\begin{equation}\label{forkummer1}
 \int_{0}^{1} P \left( x^{(1)},  x^{(2)} \right)  dx^{(2)}  =  
 (x^{(1)})^{2u^{(1)}_1 -1} (1-x^{(1)})^{2u^{(1)}_{2}-1}   \int_{0}^{1} 
 \pi_{2} \left(  x^{(2)} \right)  e^{ 2 h  x^{(2)}  x^{(1)} }  dx^{(2)}.  
\end{equation}
Here
$$  
\int_{0}^{1} 
 \pi_{2} \left(   x^{(2)} \right)  e^{ 2 h  x^{(2)}  x^{(1)} }  dx^{(2)}   
  $$
\begin{equation}\label{forkummer22}  
=\int_{0}^{1} 
(x^{(2)})^{2u^{(2)}_1 -1} (1-x^{(2)})^{2u^{(2)}_{2}-1}  
  e^{ 2  h x^{(1)} x^{(2)} }     dx^{(2)},   
  \end{equation}
Here    Kummer$^{,}$s (confluent hypergeometric)  function  $M(a,b,z)$ \cite[section 13.1]{abramowitz1964handbook}  contributes to computing the normalization  constant  in view of     
the  integral representation  \cite[eqn. 13.2.1]{abramowitz1964handbook} 
\begin{equation}\label{forkummer23} 
M(a,b,z)= \frac{\Gamma(b)}{\Gamma(a)\Gamma(b-a)}\int_0^1 e^{zu}u^{a-1}(1-u)^{b-a-1}\,du.         
 \end{equation}  
 Kummer$^{,}$s    function  has    the expansion \cite[eqn. 13.1.2]{abramowitz1964handbook} 
\begin{equation}\label{kummerexp}
    M(a,b,x)=\sum_{n=0}^\infty \frac {a^{(n)} x^n} {b^{(n)} n!}, -\infty < x <  \infty,  
\end{equation}
where $a^{(0)}=1$, $a^{(n)}=a(a+1)(a+2)\cdots(a+n-1)$.

 The integral representation of the  Kummer function gives by  (\ref{forkummer23})  in 
 (\ref{forkummer22}) with $a= 2u^{(2)}_1$, $b-a= 2u^{(2)}_2$ so that 
 $b= 2 \bar{u}^{(2)}$ 
$$
\int_{0}^{1} 
(x^{(2)})^{2u^{(2)}_1 -1} (1-x^{(2)})^{2u^{(2)}_{2}-1}  
  e^{ 2  h x^{(1)} x^{(2)} }     dx^{(2)}=    \frac{\Gamma(2u^{(2)}_1 )\Gamma(2u^{(2)}_2)  }{\Gamma(2(u^{(2)}_1+ u^{(2)}_2 ) )} M(2u^{(2)}_1,2 \bar{u}^{(2)} ,2      h x^{(1)}  ). 
$$
Then  we get in view of (\ref{forkummer1}) that 
$$
Z=\int_{0}^{1}(x^{(1)})^{2u^{(1)}_1 -1} (1-x^{(1)})^{2u^{(1)}_{2}-1}    \frac{\Gamma(2u^{(2)}_1 )\Gamma(2u^{(2)}_2)  }{\Gamma(2(u^{(2)}_1+ u^{(2)}_2 ))} M(2u^{(2)}_1,2 \bar{u}^{(2)},2      h x^{(1)}  ) dx^{(1)} 
$$
$$
=   \frac{\Gamma(2u^{(2)}_1 )\Gamma(2u^{(2)}_2)  } {\Gamma(2(u^{(2)}_1 +  u^{(2)}_2 ))}\int_{0}^{1} 
 (x^{(1)})^{2u^{(1)}_1 -1} (1-x^{(1)})^{2u^{(1)}_{2}-1}
 M \left(2u^{(2)}_1,2 \bar{u}^{(2)},  2h x^{(1)} \right )     dx^{(1)}.    
$$ 
With the expansion in  (\ref{kummerexp}) we get  
$$
  \int_{0}^{1} 
 (x^{(1)})^{2u^{(1)}_1 -1} (1-x^{(1)})^{2u^{(1)}_{2}-1}
 M \left(2u^{(2)}_1,2 \bar{u}^{(2)},  2h x^{(1)} \right )     dx^{(1)}  
 $$
 $$=  
\sum_{n=0}^\infty \frac {(2u^{(2)}_1)^{(n)}  (2h)^n} {(2 \bar{u}^{(2)})^{(n)} n!}  \int_{0}^{1} 
 (x^{(1)})^{2u^{(1)}_1 +n  -1} (1-x^{(1)})^{2u^{(1)}_{2}-1}
     dx^{(1)}
$$
$$
= \sum_{n=0}^\infty \frac {(2u^{(2)}_1)^{(n)}  (2h)^n} {(2 \bar{u}^{(2)})^{(n)} n!}  
 \frac{\Gamma(2u^{(1)}_1 +n)\Gamma(2u^{(1)} )  }{\Gamma(2(u^{(1)}_1 +u^{(1)}_{2}) +n)}.
$$
In summary, we have  found 
$$
Z=    \frac{\Gamma(2u^{(2)}_1 ) \Gamma(2u^{(2)}_2 )   } {\Gamma(2(u^{(2)}_1 +  u^{(2)}_2 ))}\sum_{n=0}^\infty \frac {(2u^{(2)}_1)^{(n)}  (2h)^n} {(2 \bar{u}^{(2)})^{(n)} n!}  
 \frac{\Gamma(2u^{(1)}_{2} )  }{\Gamma(2(u^{(1)}_1 +u^{(1)}_{2}) +n)}{\Gamma(2u^{(1)}_1 +n)}.
$$   
In other words,   
{\small 
\begin{equation}\label{aek2}
P \left( x^{(1)},  x^{(2)} \right) = \frac{(x^{(1)})^{2u^{(1)}_1 -1} (1-x^{(1)})^{2u^{(1)}_{2}-1}(x^{(2)})^{2u^{(2)}_1 -1} (1-x^{(2)})^{2u^{(2)}_{2}-1}  
  e^{ 2 h  x^{(1)} x^{(2)} } }
  { \frac{\Gamma(2u^{(2)}_1 )\Gamma(2u^{(2)}_2)  }{\Gamma(2(u^{(2)}_1+ u^{(2)}_2 ))}\sum_{n=0}^\infty \frac {(2u^{(2)}_1)^{(n)}  (2h)^n} {(2 \bar{u}^{(2)})^{(n)} n!}  
 \frac{\Gamma(2u^{(1)}_1 +n)\Gamma(2u^{(1)}_{2} )  }{\Gamma( 2(u^{(1)}_1 +u^{(1)}_{2}) +n))}}. 
\end{equation} }
If $  h =0$, then  only  the term with $n=0$  ($0^{0}=1$, $0! =1$)  in the summation in the numerator gives a non-zero contribution and  the  density in (\ref{aek2}) becomes  a product of two Beta densities, or describes two independent loci with 
two alleles and mutation,  as it should.  Or, we are for $h=0$ dealing with a pair of  independent  Wright-Fisher models  with mutation.  

The probability density function (\ref{aek2}) might     be called a (non-centralized) {\em bivariate Beta density}. However, the nome bivariate Beta density is already assigned to a different  bivariate  density, see   \cite{gupta2011non} and its references.   

In \cite {ethier1989}  the two-locus Wright-Fisher model for mutation, selection, and
random genetic drift in a panmictic, monoecious, diploid population of $N$ 
individuals is given a  diffusion approximation under various forms of  selection.  The resulting 
diffusion processes do  not seem to include explicitly   the  Svirezhev-Shahshahani  selection term   of (\ref{2dimwfdiff}).  


 }\end{exempel} \qed

\section{   Svirezhev-Shahshahani gradients   and  undirected Graphs  }\label{exempel} 
\setcounter{equation}{0}

We    derive    some  instances  of the multilocus and multiallele  model by     choices  of the structure   of $A$ above. 
This turns out to be a very flexible and effective  way to derive Wright-Fisher  diffusions of the form (\ref{bigdsd}). These examples are computational desktop constructions   and do not necessarily  emulate  any known  real-life biological situations.  
 
 In each of these examples  the normalization constant is denoted generically  as  $Z$, 
but has to be  computed anew in each example. 
In each of these examples we  take  also  single locus selection  parameters as zeroes, i.e.,  
$$
{\bf h}= {\bf 0}. 
$$

 \begin{exempel}\label{exempel311}[{\bf Two  loci, two  alleles with   selection   and mutation} ]    
{\rm  We  have $L=2$, $M_{1}=M_{2}=2$, ${\bf x} =  (x^{(1)},  x^{(2)})$. Then  we identify 
$x^{(1)}= x^{(1)}_{1},  x^{(2)}= x^{(2)}_{1}$ and have $x^{(1)}_{2}= 1-x^{(1)}_{1}$ and 
$x^{(2)}_{2}= 1-x^{(2)}_{1}$. The augmented state vector is thus 
$$
{\bf \underline{x}}^{T}= \left(  
 x^{(1)}_{1}, 
 x^{(1)}_{2},  
 x^{(2)}_{1}, 
 x^{(2)}_{2}  \right). 
 $$ 
We take take 
$$
J_{2}(M_{1}) = \left( \begin{array}{cc}  h_1  &  0 \\ 0 & h_2  \end{array} \right) = J_{1}(M_{2})
$$
so that in   (\ref{potmatrix1}) we  obtain the  symmetric $4 \times 4 $ matrix
$$
A= \left( \begin{array}{cccc} 
0 & 0 &  h_1 & 0 \\
0 & 0 & 0 & h_2 \\
h_1  &  0 & 0 & 0 \\
0  & h_2 & 0 & 0 \end{array} \right),  
$$
  Then 
 $$
{\bf \underline{x}}^{T} A {\bf \underline{x}}=  2 (h_{1} + h_{2})  x^{(1)}_{1} x^{(2)}_{1} -2h_{2} 
(x^{(1)}_{1} + x^{(2)}_{1}) + 2h_{2}.     
$$
or with  the variables $x^{(i)}$, this yields  
$$
V({\bf x}) =      \left(h_{1} + h_{2}\right)  x^{(1)}  x^{(2)}  - 2h_{2} 
\left(x^{(1)}  + x^{(2)} \right)  + h_{2}.  
$$

The  in the drift functions are  then found by  lemma \ref{driftaddition2}  as  
$$  
d_{11}^{(1)}({\bf x}) V^{'}_{x^{(1)}} = x^{(1)}\left(1-x^{(1)} \right) \cdot \left[ \left(h_{1} + h_{2} \right.)  x^{(2)}   
- h_{2} \right],      
$$
 and  
$$
d_{11}^{(2)}({\bf x}) V^{'}_{x^{(2)}} = x^{(2)}\left(1-x^{(2)} \right) \cdot \left[ \left(h_{1} + h_{2}\right)  x^{(1)}  
- h_{2} \right].      
$$.     
 Hence we obtain   the following system of stochastic  differential equations 
{\small $$
\left \{ \begin{array}{cc}
dX^{(1)}_{t} =  \mu^{(1)}( X^{(1)}_{t}) dt  +  
\left(h_{1} + h_{2} \right)   X^{(1)}_{t}(1-X^{(1)}_{t}) X^{(2)}_{t}dt -   h_{2}  X^{(1)}_{t}(1-X^{(1)}_{t})dt + \sqrt{ X^{(1)}_{t}(1-X^{(1)}_{t}) }dW^{(1)}_{t} \\
dX^{(2)}_{t} =  \mu^{(2)}( X^{(2)}_{t})  dt  +  
\left(h_{1} + h_{2}\right)   X^{(2)}_{t}(1-X^{(2)}_{t}) X^{(1)}_{t}dt -   h_{2}  X^{(2)}_{t}(1-X^{(2)}_{t})dt + \sqrt{ X^{(2)}_{t}(1-X^{(2)}_{t}) }dW^{(2)}_{t},    
\end{array} \right. 
$$} 
 
 }

\end{exempel} \qed

 \begin{exempel}\label{exempel4}[{\bf Four   loci, two  alleles with   selection   and mutation:  The General Case } ]  
{\rm  We take $L=4$, $M_{i}=2$ for $i=1,2,3,4$. We set ($^{T}$ is the vector transpose)
\begin{eqnarray}\label{frekekvvcase1} 
{\bf \underline{x}}^{T}&= &({\bf x}^{(1)}, {\bf x}^{(2)}, {\bf x}^{(3)}, {\bf x}^{(4)})  \\ \nonumber  
& & \\ \nonumber 
& =&  (x_{1}^{(1)},  x_{2}^{(1)},x_{1}^{(2)},x_{2}^{(2)}, x_{1}^{(3)},x_{2}^{(3)},x_{1}^{(4)},x_{2}^{(4)})  \nonumber  
\end{eqnarray} 
where $ x_{2}^{(i)}=1- x_{1}^{(i)}$ for $i=1,2,3,4$.  Let us furthermore set  
\begin{equation}\label{fournull}
{\bf 0}= \left( \begin{array}{cc}  0 &  0 \\ 0 & 0  \end{array} \right)  
\end{equation} 
and 
$$
{\bf J}_{2}(M_{1}) = \left( \begin{array}{cc}  h_{2} &  0 \\ 0 & 0  \end{array} \right) = {\bf J}_{1}(M_{2})
$$
$$
{\bf J}_{3}(M_{1}) = \left( \begin{array}{cc}  h_{5} &  0 \\ 0 & 0  \end{array} \right) = {\bf J}_{1}(M_{3}) 
$$
$$
{\bf J}_{3}(M_{2}) = \left( \begin{array}{cc}  h_{3} &  0 \\ 0 & 0  \end{array} \right) = {\bf J}_{2}(M_{3}) 
$$
$$
{\bf J}_{4}(M_{1}) = \left( \begin{array}{cc}  h_{7} &  0 \\ 0 & 0  \end{array} \right) = {\bf J}_{1}(M_{4})  
$$
 $$
{\bf J}_{4}(M_{2}) = \left( \begin{array}{cc}  h_{6} &  0 \\ 0 & 0  \end{array} \right) = {\bf J}_{2}(M_{4})  
$$
and 
$$
{\bf J}_{4}(M_{3}) = \left( \begin{array}{cc}  h_{4} &  0 \\ 0 & 0  \end{array} \right) = {\bf J}_{3}(M_{4}), 
$$
where the right most  inequalities are  enforced  by the symmetry required. These $2 \times 2$ matrices are symmetric, but in general the matrices ${\bf J}_{i}(M_{k})$ cannot always be symmetric, since these matrices
are not necessarily square.  We are obviously taking here  all 
$h_{i}(M_{k}) = 0$ for ease of work.  
 Then the matrix A      in    (\ref{potmatrix1}) boils down to the  symmetric $8 \times 8 $ matrix, again designated by $A$,  
    \begin{equation}\label{potmatrix15} 
A= \left( \begin{array}{cccc}
{\bf 0}  &  {\bf J}_{2}(M_{1}) &  {\bf J}_{3}(M_{1}) &    {\bf J}_{4}(M_{1})   \\ 
{\bf J}_{2}(M_{1}) &    {\bf 0}  &    {\bf J}_{3}(M_{2})    &    {\bf J}_{4}(M_{2})   \\   
{\bf J}_{3}(M_{1}) &    {\bf J}_{3}(M_{2})  &   {\bf 0} &    {\bf J}_{4}(M_{3})   \\ 
{\bf J}_{4}(M_{1}) &    {\bf J}_{4}(M_{2})  &   {\bf J}_{4}(M_{3})  &  {\bf 0}    \end{array} \right).  
\end{equation} 
We get by a simple piece of  algebra 
\begin{equation}\label{gradient0}
 A{\bf \underline{x}} =  \left(  \begin{array}{c}   h_{2} x^{(2)}_{1} +  h_{5} x^{(3)}_{1} +  h_{7} x^{(4)}_{1} \\
 0 \\
 h_{2} x^{(1)}_{1} +   h_{3} x^{(3)}_{1} +  h_{6} x^{(4)}_{1}\\
 0 \\
 h_{5} x^{(1)}_{1} +   h_{3} x^{(2)}_{1} +  h_{4} x^{(4)}_{1} \\
0 \\ 
 h_{7} x^{(1)}_{1} +   h_{6} x^{(2)}_{1} +  h_{4} x^{(3)}_{1} \\
0 
 \end{array}    \right).  
 \end{equation} 
When we revert  to  $x^{(i)}_{1}=  x^{(i)}$, 
\begin{equation}\label{wpotent}
{\bf \underline{x}}^{T} A{\bf \underline{x}} = 2 \left( 
      h_{2} x^{(1)} x^{(2)}  +  h_{5} x^{(1)}  x^{(3)}  +  h_{7} x^{(1)}  x^{(4)}      + 
h_{3} x^{(2)} x^{(3)} + h_{6} x^{(2)} x^{(4)}  +  h_{4} x^{(3)} x^{(4)}\right).  
\end{equation}
The    Svirezhev-Shahshahani gradient  $ {\bf D}\left({\bf x} \right)  \nabla_{{\bf x}} V \left({\bf x}\right)$, a $4 \times 1$ vector     in this case,   comprising the selective interactions between the four  loci in the  corresponding   stochastic differential equation   
(\ref{bigd55}), is  
\begin{equation}\label{gdriftb}
{\bf D}({\bf x})  \nabla_{{\bf x}} V \left({\bf x}\right)   = \left( \begin{array}{c}  
x^{(1)} ( 1- x^{(1)}) \left( h_{2} x^{(2)}  +  h_{5} x^{(3)}  +  h_{7} x^{(4)}  \right)    \\ 
  x^{(2)} ( 1- x^{(2)}) \left( h_{2} x^{(1)}  +   h_{3} x^{(3)}  +  h_{6} x^{(4)}   \right)  \\
x^{(3)} ( 1- x^{(3)}) \left(  h_{5} x^{(1)}  +   h_{3} x^{(2)}  +  h_{4} x^{(4)} \right)  \\
 x^{(4)} ( 1- x^{(4)}) \left( h_{7} x^{(1)}  +   h_{6} x^{(2)}  +  h_{4} x^{(3)}  \right)   \end{array} \right). 
\end{equation} 
Hence we see that every  locus  interacts  with  every other in a symmetric manner. Let us  agree to regard the four loci as nodes and to draw an undirected edge between two loci, as soon as  these  appear simultaneously   in the same product term  in $V({\bf x})$, or are jointly  in a    component of the vector ${\bf G} \left({\bf x}\right)$.  Then the current 
Wright-Fisher  model with Svirezhev-Shahshahani  selection is represented   by the complete  graph below.

\includegraphics[width=5.0cm,height=3.0cm]{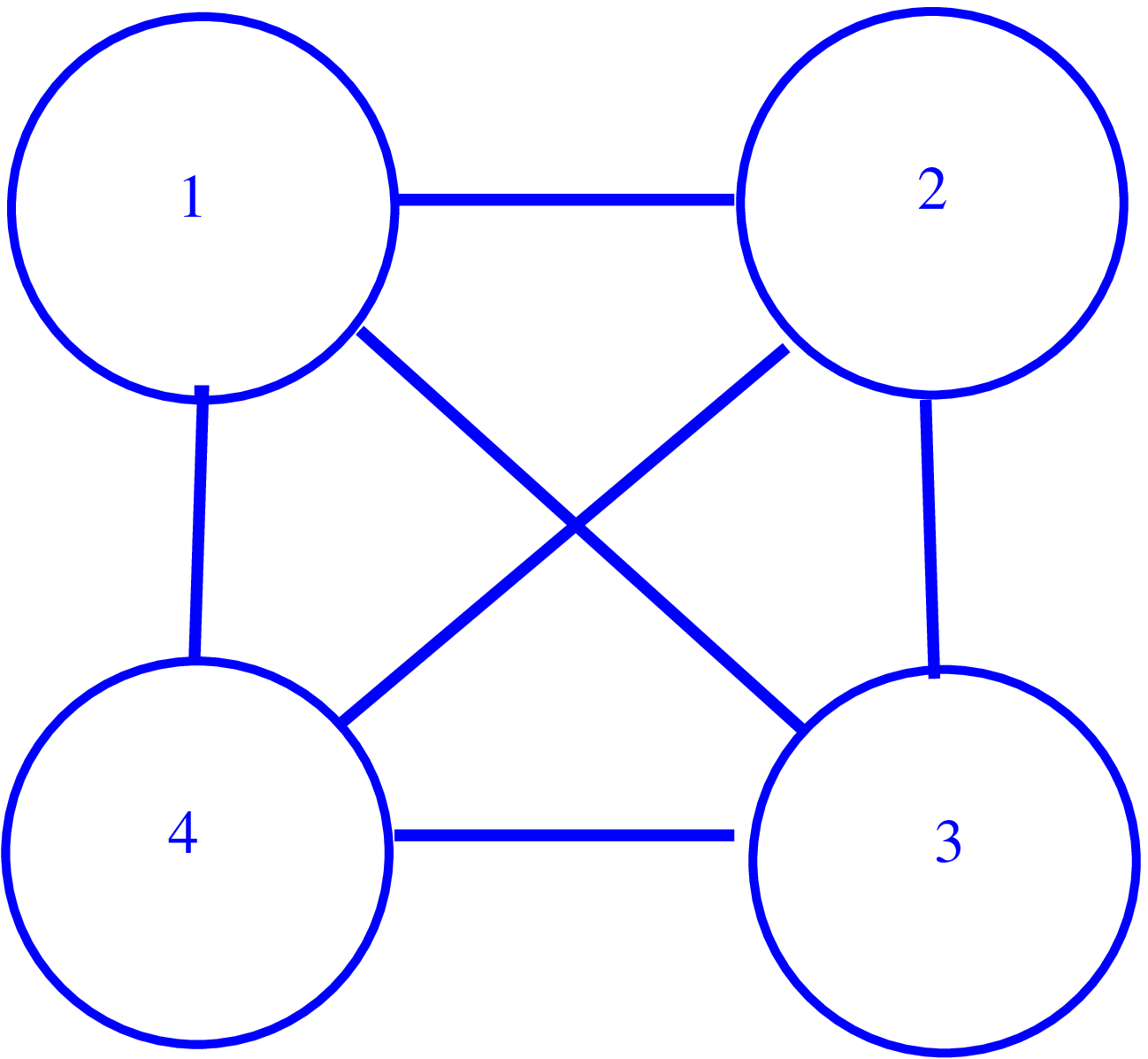}%

\noindent The normalized stationary distribution (density)     is   by  
(\ref{stationary_dist24}) equal to    
\begin{equation}\label{aekdistn2b} 
P \left({\bf x} \right) = \frac{1}{Z} \pi({\bf x}) e^{2 \left( 
      h_{2} x^{(1)} x^{(2)}  +  h_{5} x^{(1)}  x^{(3)}  +  h_{7} x^{(1)}  x^{(4)}      + 
h_{3} x^{(2)} x^{(3)} + h_{6} x^{(2)} x^{(4)}  +  h_{4} x^{(3)} x^{(4)}\right) },      
\end{equation} 
 where 
 $$
 \pi({\bf x}) =  \prod_{i=1}^{4}   (x^{(i)})^{2u^{(i)}_1 -1} (1-x^{(i)})^{2u^{(i)}_{2}-1}  
$$ 
is the product of  the non-normalized marginal Beta densities.  The standardization constant $Z$  can, at least in special cases,   again be developed using the Kummer function along the lines of  example  \ref{exempel4} but seems to produce  a  rather  messy  final formula. 

}\end{exempel}  \qed 
 
\begin{exempel}\label{exempel5}[{\bf Four   loci, two  alleles with   selection   and mutation Case I} ]  
{\rm      If   $h_{3}=h_{4}=h_{6}=0$    in    (\ref{potmatrix15}) we  get  the    matrix again denoted by $A$, 
\begin{equation} 
A= \left( \begin{array}{cccc}
{\bf 0}  &  {\bf J}_{2}(M_{1}) &  {\bf J}_{3}(M_{1}) &    {\bf J}_{4}(M_{1})   \\ 
{\bf J}_{2}(M_{1}) &    {\bf 0}  &   {\bf 0} &  {\bf 0}   \\   
{\bf J}_{3}(M_{1}) &    {\bf 0}  &   {\bf 0} &  {\bf 0}  \\ 
{\bf J}_{4}(M_{1}) &   {\bf 0}  &   {\bf 0} &  {\bf 0}    \end{array} \right).  
\end{equation} 
 This entails   by (\ref{gradient0}),    when we revert  to  $x^{(i)}_{1}=  x^{(i)}$,
\begin{equation}\label{gradient1}
 A{\bf \underline{x}} =  \left(  \begin{array}{c}   h_{2} x^{(2)}_{1} +  h_{5} x^{(3)}_{1} +  h_{7} x^{(4)}_{1} \\
 0 \\
 h_{2} x^{(1)} \\
 0 \\
  h_{5} x^{(1)}  \\
0 \\ 
 h_{7} x^{(1)}  \\
0 
 \end{array}    \right)    
\end{equation} 
 and  from (\ref{wpotent}) 
$$
{\bf \underline{x}}^{T} A{\bf \underline{x}} =  
2  x^{(1)} \left(    h_{2} x^{(2)}  +  h_{5} x^{(3)}  +  h_{7} x^{(4)}   \right)   
$$
Thus, by (\ref{gdriftb}),  
\begin{equation}\label{gdriftboo}
{\bf D}(X(t))  \nabla_{{\bf x}} V \left( {\bf x} \right)   = \left( \begin{array}{c}  
x^{(1)} ( 1- x^{(1)})   \left(h_{2} x^{(2)}  +  h_{5} x^{(3)}  +  h_{7} x^{(4)}\right) \\ 
 h_{2} x^{(2)} ( 1- x^{(2)})  x^{(1)}  \\
 h_{5} x^{(3)} ( 1- x^{(3)})  x^{(1)} \\
 h_{7} x^{(4)} ( 1- x^{(4)})  x^{(1)}  \end{array} \right). 
\end{equation}  
Hence, in this case the loci $2,3,4$ do not interact with each other except mediated by the locus $1$. 
This can be representing  a possible biologically interesting feature.  This is illustrated by    
  graph  drawn by the same principle as in the general case, i.e. in the example \ref{exempel3}.     

\includegraphics[width=5.0cm,height=3.0cm]{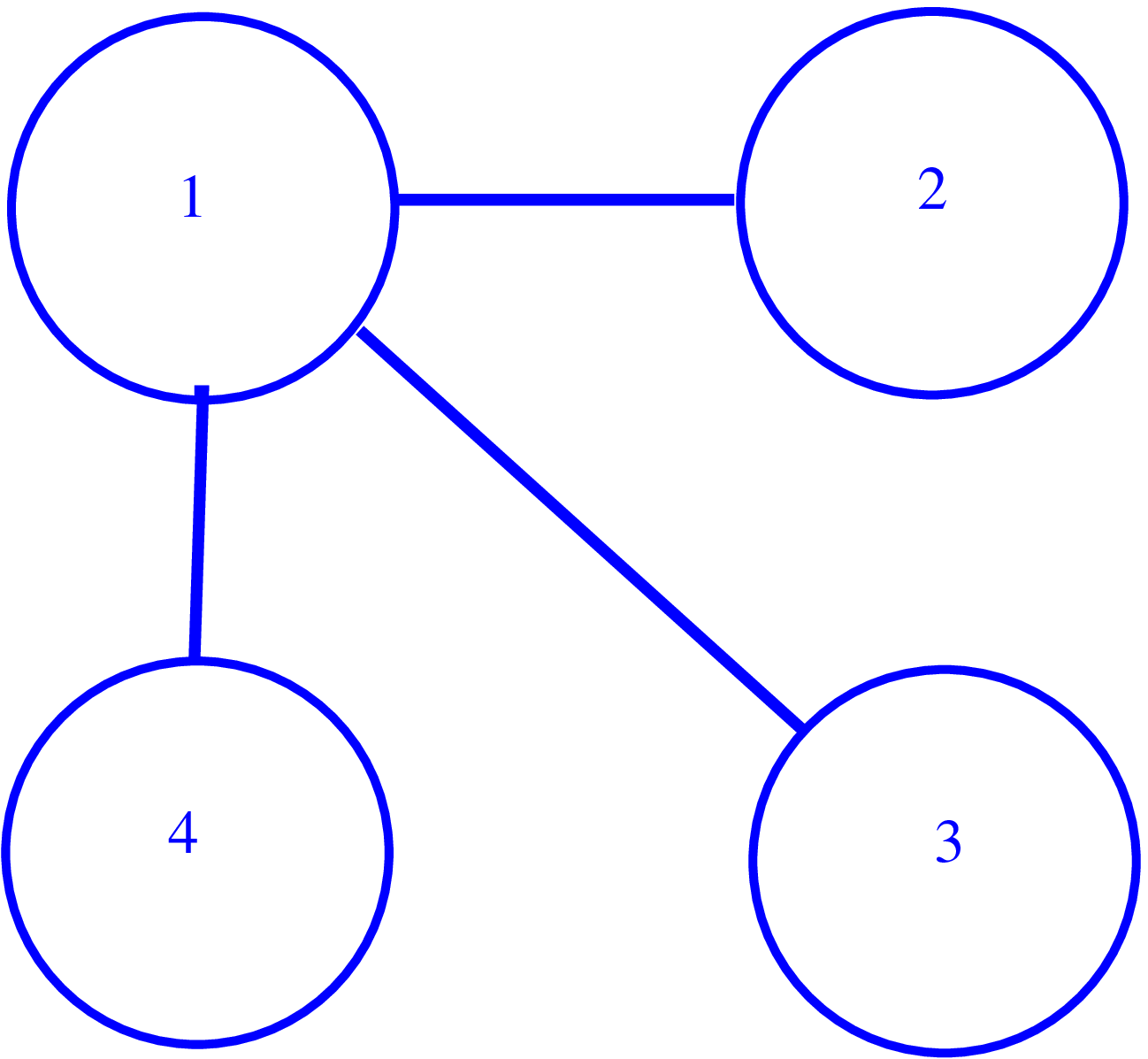}%

When  compared to the complete graph in the preceding example, we see clearly the correspondence between the  deleted  
edges and the  zero matrices ${\bf J}_{i}\left(M_{k} \right)$ of interaction  imposed.  }\end{exempel} \qed   
 
\begin{exempel}\label{exempel6}[{\bf Four   loci, two  alleles with   selection   and mutation Case II} ]  
{\rm    We  consider  in (\ref{potmatrix15})  the  the following special case, i.e.,  $h_{5}=h_{6}=h_{7}=0$, so that  
\begin{equation} 
A= \left( \begin{array}{cccc}
{\bf 0}  &  {\bf J}_{2}(M_{1}) & {\bf 0}   &   {\bf 0}   \\ 
{\bf J}_{2}(M_{1}) &    {\bf 0}  &  {\bf J}_{3}(M_{2})  &  {\bf 0}   \\   
{\bf 0} &  {\bf J}_{3}(M_{2})    &   {\bf 0} & {\bf J}_{4}(M_{3})  \\ 
  {\bf 0}  &   {\bf 0}  &  {\bf J}_{4}(M_{3})  &  {\bf 0}    \end{array} \right).  
\end{equation} 
 This gives,   by specialization of the general case,   
\begin{equation}\label{gradient2}
 A{\bf \underline{x}} =  \left(  \begin{array}{c}   h_{2}x_{1}^{(2)}  \\
 0 \\ 
     h_{2}x_{1}^{(1)} +   h_{3}x_{1}^{(3)}  \\
     0 \\
 h_{3}x_{1}^{(2)} +   h_{4}x_{1}^{(4)} \\
 0 \\
  h_{4}x_{1}^{(3)}  \\
  0    
 \end{array}    \right).
 \end{equation}  
The  return to  $x^{(i)}_{1}=  x^{(i)}$ yields  by (\ref{wpotent}) 
$$
{\bf \underline{x}}^{T} A{\bf \underline{x}} =  
   2 h_{2}x ^{(1)}x^{(2)} +  2  h_{3}x^{(2)}x^{(3)}    
      +  2 h_{4}x^{(3)}x^{(4)}.   
$$   
We get by   (\ref{gdriftb}) 
\begin{equation}\label{gdrift2}
{\bf D}\left({\bf x} \right)  \nabla_{{\bf x}} V \left({\bf x}\right)= \left( \begin{array}{c}  
 x^{(1)} ( 1- x^{(1)})  h_{2} x^{(2)}   \\ 
   x^{(2)} ( 1- x^{(2)})  \left(h_{2}x^{(1)} +   h_{3}^{(3)} x^{3}\right)     \\
   x^{(3)} ( 1- x^{(3)})  \left( h_{3}x^{(2)} +   h_{4}x ^{(4)}  \right)  \\
   x^{(4)} ( 1- x^{(4)})   h_{4}x^{(3)}   \end{array} \right). 
\end{equation}  
In this case  the loci 1 and 4 have no direct interaction  with each other, but  interact through the loci 2 and 3. 
 
The normalized stationary distribution (density)     is  equal to    
\begin{equation}\label{aekdistn23} 
P \left({\bf x} \right) = \frac{1}{Z} \pi({\bf x}) e^{2 \left( h_{2}x ^{(1)}x^{(2)} +     h_{3}x^{(2)}x^{(3)}  +  h_{4}x^{(3)}x^{(4)}  \right) },      
\end{equation} 
 where 
 $$
 \pi({\bf x}) =  \prod_{i=1}^{4}   (x^{(i)})^{2u^{(i)}_1 -1} (1-x^{(i)})^{2u^{(i)}_{2}-1}  
$$ 
is as in the preceding example.  The standardization constant $Z$  can perhaps  be developed using the Kummer function along the lines of  example  \ref{exempel4}. 

\includegraphics[width=5.0cm,height=3.0cm]{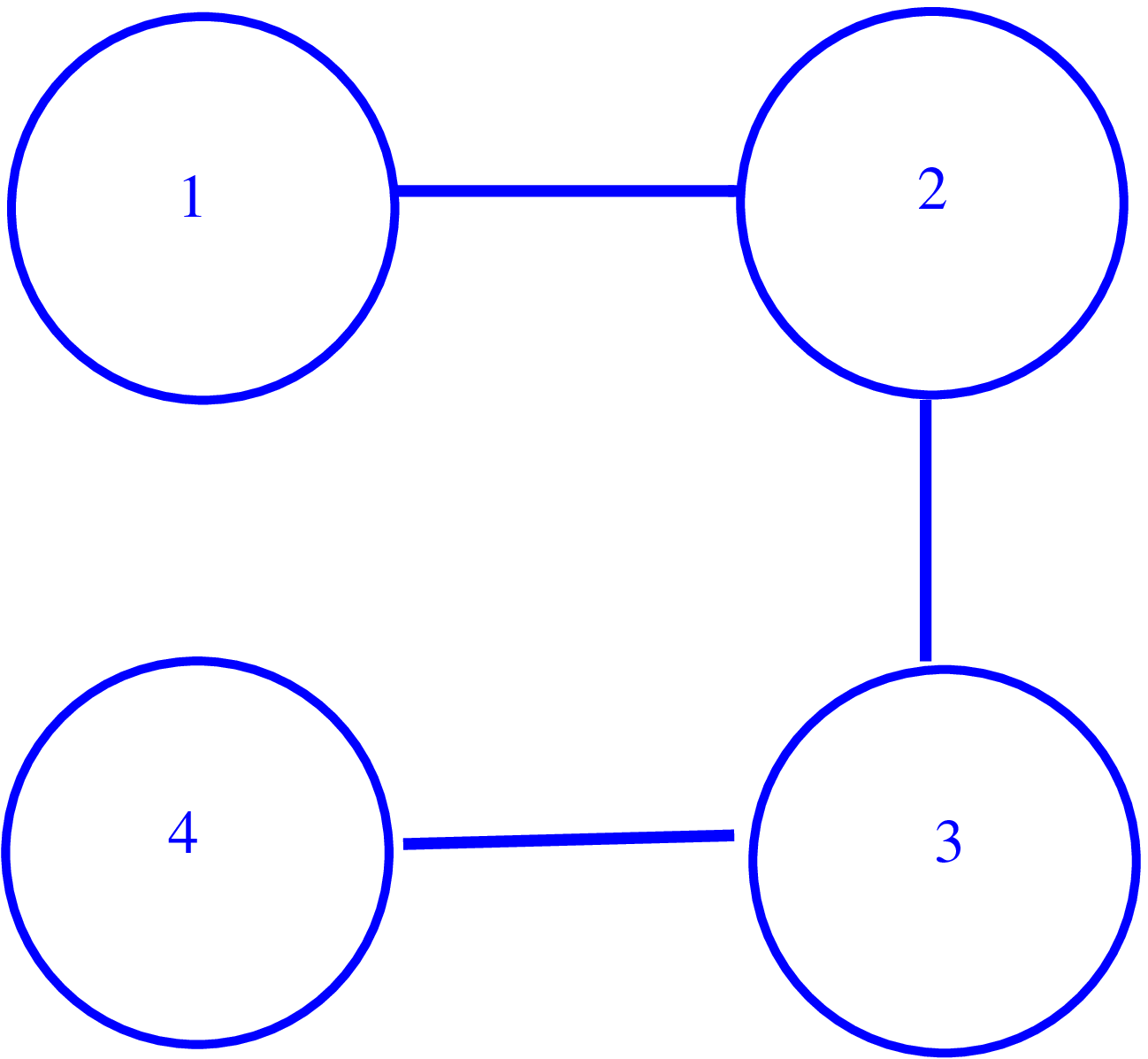}%

The graph above is drawn by the same principle as in the preceding examples. 

}\end{exempel} \qed    
    
\begin{exempel}\label{exempel61}[{\bf Four   loci, two  alleles with   selection   and mutation Case III} ]  
{\rm     We take in (\ref{potmatrix15}) $h_{2}=h_{3}=h_{4} =0$  to  get  the  symmetric $8 \times 8 $ matrix  
\begin{equation} 
A= \left( \begin{array}{cccc}
{\bf 0}  &  {\bf 0}  &  {\bf J}_{4}(M_{1})   &    {\bf J}_{4}(M_{1})    \\ 
{\bf 0}  &    {\bf 0}  &    {\bf 0} &  J_{4}(M_{4})    \\   
 {\bf J}_{3}(M_{1}) & {\bf 0}      &   {\bf 0} &  {\bf 0}   \\ 
  {\bf J}_{4}(M_{1})  &     {\bf J}_{4}(M_{2}) & {\bf 0}     &  {\bf 0}    \end{array} \right).  
\end{equation} 
 This gives  with   $x^{(i)}_{1}=  x^{(i)}$ and by (\ref{gradient0}) 
\begin{equation}\label{gradient3}
 A{\bf \underline{x}} =  \left(  \begin{array}{c}   h_{5}x^{(3)} +  h_{7}x^{(4)} \\
 0 \\ 
         h_{6}x^{(4)}  \\
     0 \\
 h_{5}x^{(1)}   \\
 0 \\
 h_{7} x^{(1)} +   h_{6}x^{(2)}  \\
  0    
 \end{array}    \right).
 \end{equation}  
Thus 
$$
{\bf \underline{x}}^{T} A{\bf \underline{x}} =  2 \left(
    h_{5}x^{(1)}x^{(3)} +    h_{6} x^{(2)}x^{(4)} +  h_{7} x^{(1)} x^{(4)} \right). 
$$    
Hence we get 
\begin{equation}\label{gdrift3}
  {\bf D}\left({\bf x} \right)  \nabla_{{\bf x}} V \left({\bf x}\right)= \left( \begin{array}{c}  
 x^{(1)} ( 1- x^{(1)})  (h_{5}x^{(3)} + h_{7}x^{(4)} )   \\ 
   x^{(2)} ( 1- x^{(2)})  h_{6}x^{(4)}         \\
   x^{(3)} ( 1- x^{(3)}) h_{5}x^{(1)}     \\
   x^{(4)} ( 1- x^{(4)})  (h_{7} x^{(1)} +   h_{6}x^{(2)})   \end{array} \right). 
\end{equation} 
 In this case the locus 1 interacts with the locus 2 only through the locus 4, and  the  
 locus 4 interacts with the locus 3 only through the locus 1. The graph is obvious.  
\includegraphics[width=5.0cm,height=3.0cm]{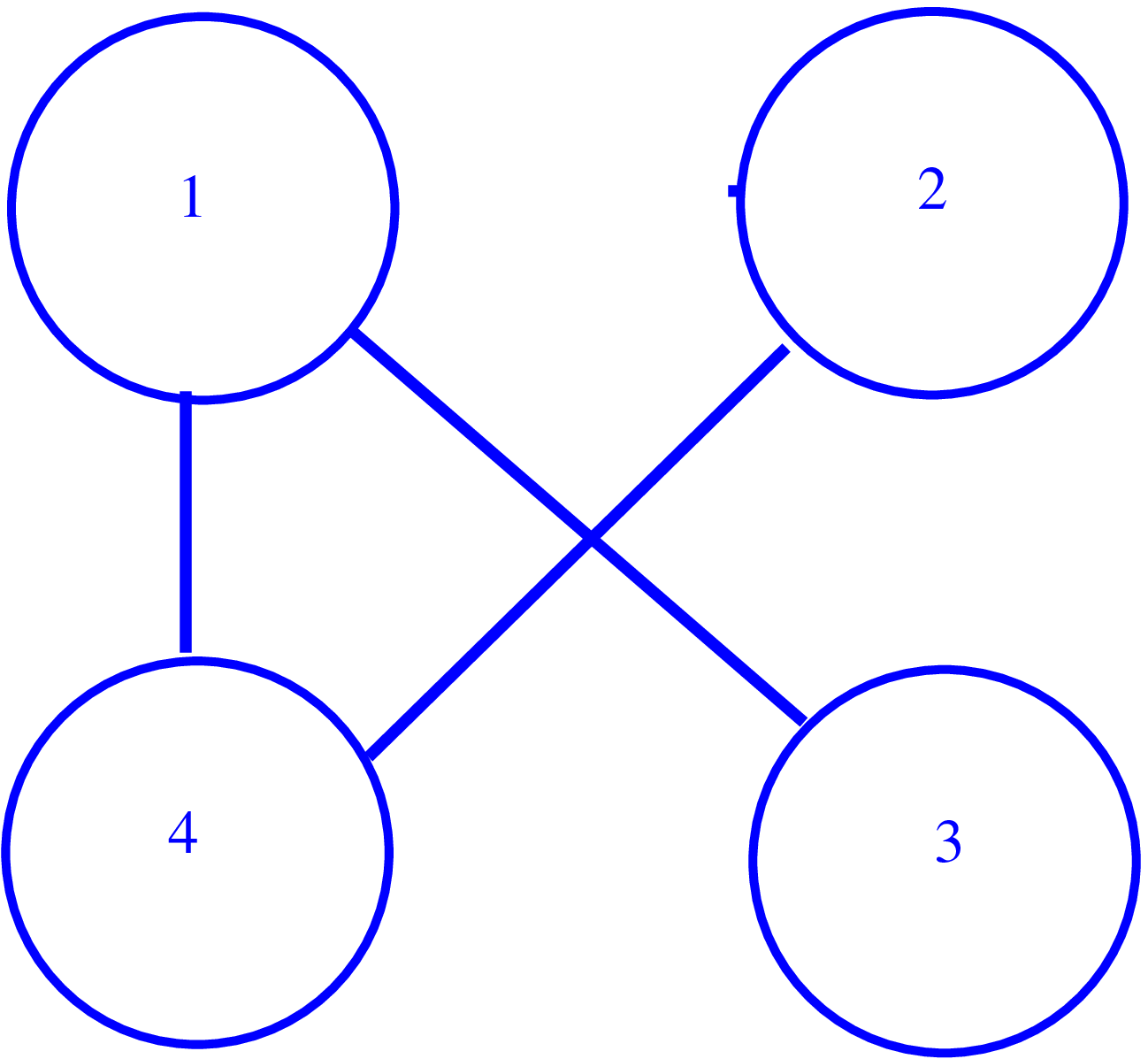}%

}\end{exempel} \qed  

 \begin{exempel}\label{exempel7}[{\bf Three   loci, two  alleles with   selection   and mutation.  The General Case } ] 
{\rm If we take $h_{4}=h_{6}=h_{7}=0$  in  (\ref{potmatrix15}), we eliminate $x^{(4)}$ and what remains 
or appears in  a Wright-Fisher model with $L=2$ and $M_{i}=2$, for $i=1,2,3$ and the $6 \times  6$ matrix, again denored by $A$,         
    \begin{equation}\label{potmatrix33} 
A= \left( \begin{array}{ccc}
{\bf 0}  &  {\bf J}_{2}(M_{1}) &  {\bf J}_{3}(M_{1})     \\ 
{\bf J}_{2}(M_{1}) &    {\bf 0}  &    {\bf J}_{3}(M_{2})      \\   
{\bf J}_{3}(M_{1}) &    {\bf J}_{3}(M_{2})  &   {\bf 0}    \\ 
    \end{array} \right).  
\end{equation} 
Here we see clearly that the model  with three loci and two nodes is nested inside the model with 
four  loci and two nodes.  To write down the  invariant density here and to inspect various special cases, one only needs to invoke the formulas  (\ref{wpotent}) and/or  (\ref{gdriftb}) with  $h_{4}=h_{6}=h_{7}=0$ 
and other specifications of zeroes.

}\end{exempel} \qed

 \begin{exempel}\label{exempel8}[{\bf Six    loci, two  alleles with   selection   and mutation.  The General Case } ] 
{\rm  But    $L=5$ and  and $M_{i}=2$, for $i=1,2,3,4$  is, of course  nested inside any  model  with a higher  number of 
loci and two alleles.     The for $L=6$ we have  the $12 \times 12$ matrix  
    \begin{equation}\label{potmatrix66} 
A= \left( \begin{array}{cccccc}
{\bf 0}  &  {\bf J}_{2}(M_{1}) &  {\bf J}_{3}(M_{1})   &  {\bf J}_{4}(M_{1})&  {\bf J}_{5}(M_{1})&  {\bf J}_{6}(M_{1})  \\ 
{\bf J}_{2}(M_{1}) &    {\bf 0}  &    {\bf J}_{3}(M_{2})  &  {\bf J}_{4}(M_{2})&  {\bf J}_{5}(M_{2})&  {\bf J}_{6}(M_{3})      \\   
{\bf J}_{3}(M_{1}) &    {\bf J}_{3}(M_{2})  &   {\bf 0} & {\bf J}_{4}(M_{3}) &    {\bf J}_{5}(M_{3}) & {\bf J}_{6}(M_{3})  \\
   
{\bf J}_{4}(M_{1}) &    {\bf J}_{4}(M_{2})  &    {\bf J}_{4}(M_{3}) & {\bf 0} &    {\bf J}_{5}(M_{4}) & {\bf J}_{6}(M_{4}) \\
{\bf J}_{5}(M_{1}) &    {\bf J}_{5}(M_{2})  &    {\bf J}_{5}(M_{3}) &     {\bf J}_{5}(M_{4}) & {\bf 0} & {\bf J}_{6}(M_{5}) \\ 
{\bf J}_{6}(M_{1}) &    {\bf J}_{6}(M_{2})  &    {\bf J}_{6}(M_{3}) &     {\bf J}_{6}(M_{4}) &   {\bf J}_{6}(M_{5}) & {\bf 0}
    \end{array} \right).  
\end{equation} 
If we want to write down the formulas like in the above examples here,  we  need to amend the matrix (\ref{potmatrix15}) with  the nine additional  matrices 
$$
{\bf J}_{i}(M_{k}) = \left( \begin{array}{cc}  h_{ik} &  0 \\ 0 & 0  \end{array} \right)  (= {\bf J}_{k}(M_{i}))
$$  
for    $k=1,2,3,4$ if $i=5$ and  $k=1,2,3,4,5$  if  $i=6$ and to perform the necessary  matrix multiplications. 
But even without any numbers  we can state  something.  Suppose we have in (\ref{potmatrix66}) 
 \begin{equation}\label{potmatrix67} 
A=  \left( \begin{array}{cccccc}
{\bf 0}  &  {\bf J}_{2}(M_{1}) &  {\bf J}_{3}(M_{1})   &  {\bf J}_{4}(M_{1})&  {\bf 0}&   {\bf 0}   \\ 
{\bf J}_{2}(M_{1}) &    {\bf 0}  &    {\bf J}_{3}(M_{2})  &  {\bf J}_{4}(M_{2})&   {\bf 0}&   {\bf 0}       \\   
{\bf J}_{3}(M_{1}) &    {\bf J}_{3}(M_{2})  &   {\bf 0} & {\bf J}_{4}(M_{3}) &    {\bf 0} &  {\bf 0}  \\
   
{\bf J}_{4}(M_{1}) &    {\bf J}_{4}(M_{2})  &    {\bf J}_{4}(M_{3}) & {\bf 0} &     {\bf J}_{5}(M_{4})   &  {\bf J}_{6}(M_{4})  \\
 {\bf 0}  &    {\bf 0}   &   {\bf 0}  &     {\bf J}_{5}(M_{4})  & {\bf 0} & {\bf J}_{6}(M_{5}) \\ 
 {\bf 0}  &    {\bf 0}   &     {\bf 0}  &      {\bf J}_{6}(M_{4})   &   {\bf J}_{6}(M_{5}) & {\bf 0}
    \end{array} \right).  
\end{equation}
But then  it clearly  holds, assuming that there  are no further zero matrices in (\ref{potmatrix67}),  that  the loci  1 to 3  interact with 5 and 6    only through the locus 4, and vice versa by the symmetries assumed.  
The interaction graph is below.

\includegraphics[width=5.0cm,height=3.0cm]{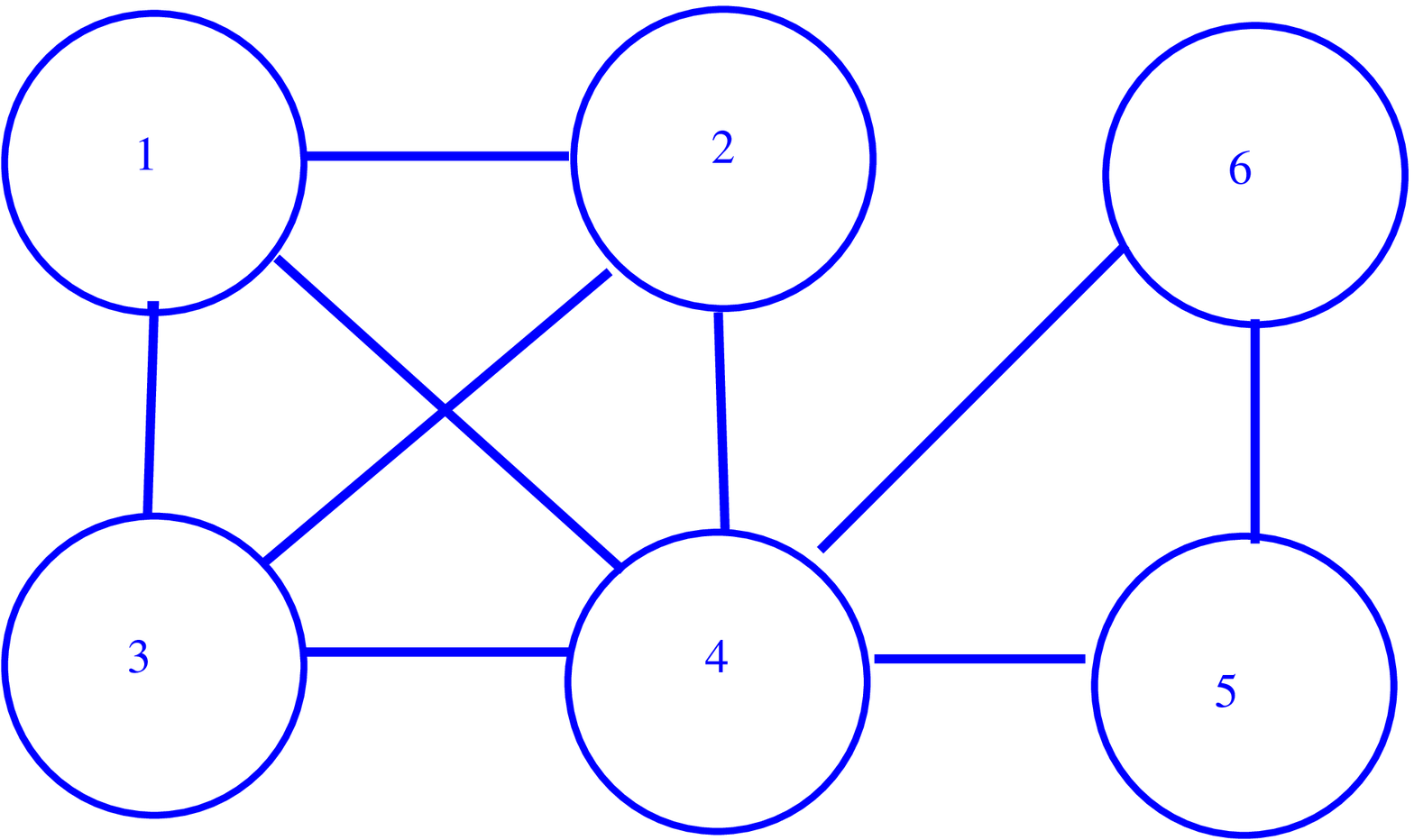}%

}\end{exempel} \qed

 \begin{exempel}\label{exempel9}[{\bf Eight   loci, two  alleles with   selection   and mutation. The converse } ] 
{\rm  Continuing with  Example \ref{exempel8} in  this manner  it is  easy to find, by extension of  the matrix $A$ in (\ref{potmatrix67}), and hence ${\bf G}({\bf x})$ and 
$V({\bf x})$ and then draw the interaction graph.   

However, even the converse is true.  If we are given  the graph in the next figure and are told that there are two alleles at every locus, we can  find   the corresponding $A$.    

\includegraphics[width=5.0cm,height=3.0cm]{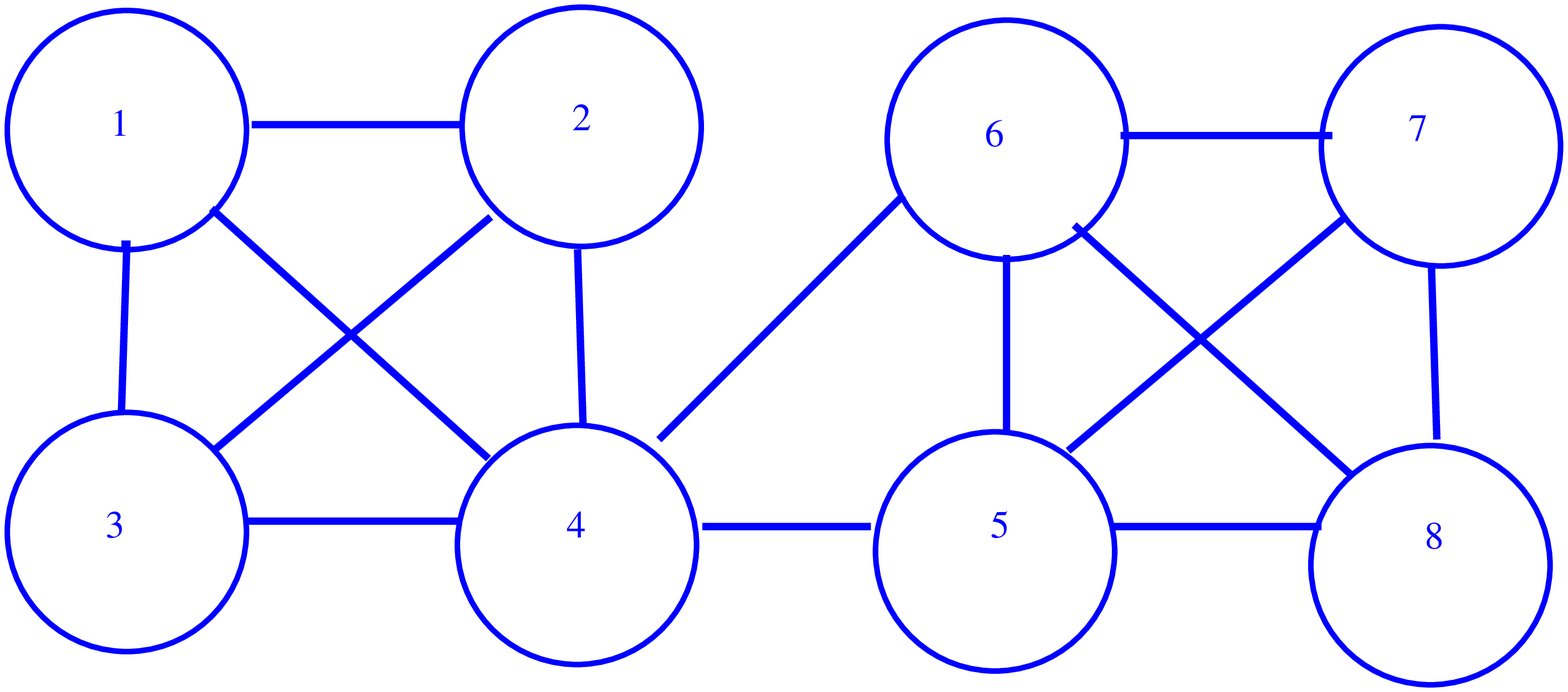}%
 
}\end{exempel} \qed

\begin{exempel}\label{exempel10}[{\bf Two  loci, one with  three   alleles  and  the other locus with  two alleles at two loci } ] 
{\rm If there are, e.g. three alleles at locus $2$, $M_{1}=2$, $M_{2}=3$  the corresponding two locus interaction matrices are  taken as 
 $$
{\bf J}_{1}(M_{2}) = \left( \begin{array}{ccc}  h_{1} &    h_{2}   & 0 \\ 
0&   0  & 0                                                                                                                                                                                                                        \end{array} \right)   
$$ 
and 
 $$
{\bf J}_{2}(M_{1}) = \left( \begin{array}{cc }   h_{1}  &       0 \\ 
 h_{2} &  0  \\
0  & 0  \end{array} \right).
$$ 
and matrix $A$ is    
$$
A= \left( \begin{array}{ccccc}  
0 & 0    &   h_{1} &    h_{2} & 0    \\ 
0 & 0  &  0 &   0   & 0 \\ 
  h_{1}  &  0 & 0 & 0 & 0 \\ 
 h_{2}  &  0 & 0 & 0 & 0 \\
 0 & 0 &   0 & 0 & 0 \\ 
    \end{array} \right). 
$$
With $$
{\bf \underline{x}}^{T}= \left(  
 x^{(1)}_{1}, 
 x^{(1)}_{2},  
 x^{(2)}_{1}, 
 x^{(2)}_{2},
 x^{(2)}_{3}
  \right)  
 $$
we get 
$$
A{\bf \underline{x}}= \left( \begin{array}{c}  
    h_{1} x^{(2)}_{1} +  h_{2}  x^{(2)}_{2}     \\ 
0  \\ 
  h_{1} x^{(1)}_{1}  \\ 
 h_{2} x^{(1)}_{1}    \\
 0  \\ 
    \end{array} \right).
$$
This entails
$$
V \left({\bf x}\right)=  V \left(  x^{(1)}_{1}, x^{(2)}_{1},  x^{(2)}_{2} \right) =  x^{(1)}_{1} \left(  h_{1} x^{(2)}_{1} +  h_{2}  x^{(2)}_{2}\right).
$$ 
From   (\ref{eadrftfact9991}) we get  here  
$$
  {\bf D}\left({\bf x} \right)  \nabla_{{\bf x}} V \left({\bf x}\right)=      
$$
\begin{equation} 
 = \left(  \begin{array}{ccc}
 x^{(1)}(1-x^{(1)}) & 0 &  0  \\
 0  &   x_{1}^{(2)}(1-x_{1}^{(2)})  &  - x_{1}^{(2)}x_{2}^{(2)} \\ 
 0 &    - x_{1}^{(2)}x_{2}^{(2)}   &   x_{2}^{(2)}(1-x_{2}^{(2)})  \end{array} \right) 
 \left(  \begin{array}{c} 
    h_{1} x^{(2)}_{1} +  h_{2}  x^{(2)}_{2}  \\
   h_{1} x^{(1)}_{1}  \\ 
    h_{2} x^{(1)}_{1}  \end{array} \right)
\end{equation}  
\begin{equation} 
 = \left(  \begin{array}{c} 
 (x^{(1)}(1-x^{(1)}) (  h_{1} x^{(2)}_{1} +  h_{2}  x^{(2)}_{2} ) \\
  h_{1} x^{(1)}_{1} x_{1}^{(2)}(1-x_{1}^{(2)})   -  h_{1} x^{(1)}_{1}  x_{1}^{(2)}x_{2}^{(2)} \\
   - h_{1} x^{(1)}_{1}  x_{1}^{(2)}x_{2}^{(2)}  +    h_{2} x^{(1)}_{1}  x_{2}^{(2)}(1-x_{2}^{(2)}) 
   \end{array} \right).
 \end{equation} 
We cannot give the corresponding  SDE in  fully  explicit form  without treating  $ {\bf D}\left({\bf x} \right)^{1/2}$, 
which  is not meaningful  for the current purpose.     
 
}\end{exempel} \qed

\section{Final Comments} 

An ultimate goal  to learn the structure of  the matrix $A$ from  data, which  will be done elsewhere, 
once successful  learning algorithms have been established.  Here the structure  of the SDE 
\begin{equation}\label{shesav} 
  dX(t) =  {\bf \mu}(X(t))dt + {\bf D}(X(t))  \nabla_{{\bf x}} V \left(X(t)\right)dt  +  {\bf D}^{1/2}\left(X(t)\right)d {\bf W}(t),   
\end{equation}  
will  be crucial, as  $V \left(X(t)\right)$ is easily  determined by  the structures  to be learned.      
The studies in  \cite{favero2019dual} and  \cite{garcia2019exact}  are exploring some of the steps  
required  for  statistical inference with  (\ref{shesav}).

There is a certain degree of skepticism  on diffusion approximations in population genetics  voiced    by  John F. Kingman in \cite[p.39]{kingman1980mathematics}.  Amongst other things, 
the question of  justifying  the  stationary density by a diffusion approximation is not simple. Techniques   for  this are  given in \cite{norman1975limit} for   one dimensional  Wright-Fisher models. 

The problem  of ergodicity and existence of an invariant measure is  studied   in \cite{kbasak1992stability} for a class of degenerate    multidimensional diffusions, but does not discuss  explicitly the  Wright-Fisher models. This analysis  deals with  the properties of transition probability function of the Wright-Fisher  diffusion, see, e.g.,   \cite{barbour2000transition}.  

In \cite{tran2013introduction}  the Fokker-Planck equation associated
with Wright-Fisher model with  two alleles at a diploid locus under
random genetic drift in a population of fixed size   without mutation or selection is shown to  possess  a unique (global)  solution. The expression   for the invariant density is, roughly stated,    a series expansion  in terms of  
Gegenbauer polynomials and eigenvectors  of the Wright-Fisher generator.  By aid of this density formula 
 these  authors can find,  .e.g.,   the
expectation and the second moment of the absorption time,
 fixation probabilities, the probability of  coexistence, or the probability of heterogeneity. 
Similar  applications  using our formula have not been attempted.

\section{Acknowledgements} 
Prof. Jukka Corander, Faculty of Medicine,   University of  Oslo, is thanked for  communicating to  the  authors  about the emergence  of  new  data in the form of  time-series of allele frequencies.  

This research has been supported by the Swedish Science
Council through grant 621-2012-2982 (EA), by the Academy of Finland through its Center of Excellence COIN (EA),  and by the Chinese Academy of Sciences
through its CAS$^{'}$  President$^{'}$s International Fellowship Initiative (PIFI) GRANT No. 2016VMA002 (EA) and    the Swedish Science  Council through grant  40-2012-5952  (TK).         
 
The funding agencies have had no involvement  in study design and in the writing of the report; and in the decision to
submit the article for publication.  There are no conflicts of interest. 
 

\newpage
\section*{Appendix A: Fitness }\label{Viability}   
\renewcommand{\theequation}{\rm{A}.\arabic{equation}}
\setcounter{equation}{0} 
 
For $\bar{m}$, we have
\begin{equation}
\begin{aligned} 
 \bar{m}= \sum_{\sigma} f(\sigma) m_{\sigma}= \sum_{\sigma} f(\sigma) \left (1+\sum_{r=1} h_r(\sigma_r)+\sum_{1 \leq r < s \leq L} J_{rs}(\sigma_r,\sigma_s) \right )\\
=1+\sum_{r=1}^L\underbrace{\sum_{\sigma}f(\sigma)h_r(\sigma_r)}_{\sum_{m=1}^{M_r}h_r(t)x_t^{(r)}}+\sum_{1 \leq r < s \leq L}\underbrace{\sum_{\sigma}f(\sigma)J_{rs}(\sigma_r,\sigma_s)}_{\sum_{t=1}^{M_r}\sum_{n=1}^{M_s}J_{rs}(t,n)x_t^{(r)}x_n^{(s)}}\\
=1+\sum_{r=1}^L\sum_{t=1}^{M_r}h_r(t)x_t^{(r)}+\sum_{1 \leq r < s \leq L}\sum_{t=1}^{M_r}\sum_{n=1}^{M_s}J_{rs}(t,n)x_t^{(r)}x_n^{(s)}.
\end{aligned}
\end{equation}
The underbraces  are due to the calculations
\begin{equation}
\sum_{\sigma}f(\sigma)h_r(\sigma_r)=\sum_{m=1}^{M_r} \sum_{\underset{\sigma_r=t}{\sigma:}}f(\sigma)h_r(\sigma_r)=\sum_{t=1}^{M_r}h_r(m) \underbrace{\sum_{\underset{\sigma_r=t}{\sigma:}}f(\sigma)}_{x_t^{(r)}}=\sum_{t=1}^{M_r}h_r(t)x_t^{(r)},
\end{equation}
and 
\begin{equation}
\sum_{\sigma}f(\sigma)J_{rs}(\sigma_r,\sigma_s)=\sum_{t=1}^{M_r}\sum_{n=1}^{M_s}J_{rs}(t,n) \sum_{\underset{\substack{\sigma_r=t \\ \sigma_s=n}}{\sigma:}}f(\sigma)=\sum_{t=1}^{M_r}\sum_{n=1}^{M_s}J_{rs}(t,n)x_t^{(r)}x_n^{(s)}.
\end{equation}
For $\bar{m}_k^{(i)}$, we have
\begin{equation}
\begin{aligned} 
 \bar{m}_k^{(i)}= \sum_{\sigma} f_k^{(i)}(\sigma) m_{\sigma}= \sum_{\sigma} f_k^{(i)}(\sigma) \left (1+\sum_{r=1}^L h_r(\sigma_r)+\sum_{1 \leq r < s \leq L} J_{rs}(\sigma_r,\sigma_s) \right )\\
=1+\sum_{r=1}^L\sum_{\sigma}f_k^{(i)}(\sigma)h_r(\sigma_r)+\sum_{1 \leq r < s \leq L}\sum_{\sigma}f_k^{(i)}(\sigma)J_{rs}(\sigma_r,\sigma_s).
\end{aligned}
\end{equation}
Per definition it holds that $f_k^{(i)}(\sigma)=\frac{\delta_{\sigma_i,k}}{x_k^{(i)}}f(\sigma)$, which gives for the second term
\begin{equation}
\begin{aligned} 
\sum_{r=1}^L\sum_{\sigma}f_k^{(i)}(\sigma)h_r(\sigma_r)&=\sum_{\sigma}f_k^{(i)}(\sigma)h_i(\sigma_i)+\sum_{\underset{r \neq i}{r=1}}^L\sum_{\sigma}f_k^{(i)}(\sigma)h_r(\sigma_r)\\
&=\frac{1}{x_k^{(i)}}\sum_{\underset{\sigma_i=k}{\sigma:}}f(\sigma)h_i(\sigma_i)+\frac{1}{x_k^{(i)}}\sum_{\underset{r \neq i}{r=1}}^L\sum_{\underset{\sigma_i=k}{\sigma:}}f(\sigma)h_r(\sigma_r)\\
&=\frac{1}{x_k^{(i)}}h_i(k)\underbrace{\sum_{\underset{\sigma_i=k}{\sigma:}}f(\sigma)}_{x_k^{(i)}}+\frac{1}{x_k^{(i)}}\sum_{\underset{r \neq i}{r=1}}^L\sum_{t=1}^{M_r}h_r(t)\underbrace{\sum_{\underset{\substack{\sigma_i=k \\ \sigma_r=t}}{\sigma:}}f(\sigma)}_{x_k^{(i)}x_m^{(r)}}\\
&=h_i(k)+\sum_{\underset{r \neq i}{r=1}}^L\sum_{t=1}^{M_r}h_r(t)x_t^{(r)}.
\end{aligned}
\end{equation}
Now   we assume (\ref{symmetri}), $J_{rs}(k,l)=J_{sr}(l,k)$.   
This yields 
\begin{equation}
\begin{aligned} 
\sum_{1 \leq r < s \leq L}\sum_{\sigma}f_k^{(i)}(\sigma)J_{rs}(\sigma_r,\sigma_s)=\sum_{\underset{r \neq i}{r=1}}^L\sum_{\sigma}f_k^{(i)}(\sigma)J_{ir}(\sigma_i,\sigma_r)+\sum_{\underset{r,s \neq i}{1 \leq r < s \leq L}}\sum_{\sigma}f_k^{(i)}(\sigma)J_{rs}(\sigma_r,\sigma_s)\\
=\frac{1}{x_k^{(i)}}\sum_{\underset{r \neq i}{r=1}}^L\sum_{\underset{\sigma_i=k}{\sigma:}}f(\sigma)J_{ir}(\sigma_i,\sigma_r)+\frac{1}{x_k^{(i)}}\sum_{\underset{r,s \neq i}{1 \leq r < s \leq L}}\sum_{\underset{\sigma_i=k}{\sigma:}}f(\sigma)J_{rs}(\sigma_r,\sigma_s)\\
=\frac{1}{x_k^{(i)}}\sum_{\underset{r \neq i}{r=1}}^L\sum_{t=1}^{M_r}J_{ir}(k,t)\underbrace{\sum_{\underset{\substack{\sigma_i=k\\\sigma_r=t}}{\sigma:}}f(\sigma)}_{x_k^{(i)}x_t^{(r)}}+\frac{1}{x_k^{(i)}}\sum_{\underset{r,s \neq i}{1 \leq r < s \leq L}}\sum_{t=1}^{M_r}\sum_{n=1}^{M_s}J_{rs}(t,n)\underbrace{\sum_{\underset{\substack{\sigma_i=k\\\sigma_r=t\\\sigma_s=n}}{\sigma:}}f(\sigma)}_{x_k^{(i)}x_t^{(r)}x_n^{(s)}}\\
=\sum_{\underset{r \neq i}{r=1}}^L \sum_{t=1}^{M_r}J_{ir}(k,t)x_t^{(r)}+\sum_{\underset{r,s \neq i}{1 \leq r < s \leq L}}\sum_{t=1}^{M_r}\sum_{n=1}^{M_s}J_{rs}(t,n)x_t^{(r)}x_n^{(s)},
\end{aligned}
\end{equation}
thus completing the expression for $\bar{m}_k^{(i)}$ as
\begin{equation}
\begin{aligned}
 \bar{m}_k^{(i)}= 1&+h_i(k)+\sum_{\underset{r \neq i}{r=1}}^L\sum_{t=1}^{M_r}h_r(\sigma_r)x_t^{(r)}\\
&+\sum_{\underset{r \neq i}{r=1}}^L \sum_{t=1}^{M_r}J_{ir}(k,t)x_t^{(r)}+\sum_{\underset{r,s \neq i}{1 \leq r < s \leq L}}\sum_{t=1}^{M_r}\sum_{n=1}^{M_s}J_{rs}(t,n)x_t^{(r)}x_n^{(s)}.
\end{aligned}
\end{equation}
\section*{Appendix B: Proof of Lemma \ref{ealemma3}   } \label{matrix}   
\renewcommand{\theequation}{\rm{B}.\arabic{equation}}
\setcounter{equation}{0} 
 {\em Proof:}  We  evaluate   
$\sum_{l=1}^{M_{i}} d_{kl}^{(i)} \widetilde{h}_{i}(l) $. We have from (\ref{hejp2}) that 
$$
\sum_{l=1}^{M_{i}} d_{kl}^{(i)} \widetilde{h}_{i}(l) =  
\sum_{l=1}^{M_{i}} d_{kl}^{(i)} h_i(l)  + \sum_{l=1}^{M_{i}} d_{kl}^{(i)} 
\sum_{\underset{l \neq k}{r=1}}^{L} \sum_{m=1}^{M_r}  J_{ir}(l,m)x_{m}^{(r)}.  
$$
By  definition of   $d_{kl}^{(i)} $  
$$
\sum_{l=1}^{M_{i}} d_{kl}^{(i)} h_i(l) = - \sum_{\underset{l \neq k}{l=1}} ^{M_{i}} x^{(i)} _{k}x^{(i)} _{l}  h_i(l) 
+ x^{(i)} _{k}(1-x^{(i)} _{k})  h_i(k) = - \sum_{ r =1} ^{M_{i}} x^{(i)} _{k}x^{(i)} _{r}  h_i(r)  +  x^{(i)} _{k}   h_i(k)
$$
\begin{equation}\label{del1} 
= x^{(i)} _{k} \left( h_i(k)  -  \sum_{ l =1} ^{M_{i} } x^{(i)} _{l}  h_i(l) \right),    
\end{equation}
which identifies the two first  terms  in  the right hand side (\ref{pairwise_selection3}). Secondly,  we have 
$$
\sum_{l=1}^{M_{i}} d_{kl}^{(i)} 
\sum_{\underset{r \neq i}{r=1}}^{L}\sum_{m=1}^{M_i}  J_{ir}(l,m)x_{m}^{(r)}  
$$ 
$$
=-\sum_{\underset{l \neq i}{l=1}}^{M_{i}} x^{(i)}_{k}x^{(i)}_{l}\sum_{\underset{r \neq i}{r=1}}^{L}\sum_{m=1}^{M_r}  J_{ir}(l,m)x_{m}^{(r)}   + x^{(i)} _{k}(1-x^{(i)} _{k}) \sum_{\underset{r \neq i}{r=1}}^{L}\sum_{m=1}^{M_r}  J_{ir}(k,m)x_{m}^{(r)}   
$$ 
$$
=   x^{(i)}_{k} \left( -\sum_{\underset{l \neq k}{l=1}}^{M_{i}} x^{(i)}_{l}\sum_{\underset{r \neq i}{r=1}}^{L}\sum_{m=1}^{M_r}  J_{ir}(l,m)x_{m}^{(r)}          
  + (1-x^{(i)}_{k}) \sum_{\underset{r \neq i}{r=1}}^{L}\sum_{m=1}^{M_r}  J_{ir}(k,m)x_{m}^{(r)}    \right)  
$$
$$
=  x^{(i)}_{k} \left(  -\sum_{\underset{r \neq i}{r=1}}^{L}\sum_{m=1}^{M_i} x_{m}^{(r)}\sum_{\underset{l \neq k}{l=1}}^{M_{r}} x^{(i)}_{l} J_{ir}(l,m)     +  \sum_{\underset{r \neq i}{r=1}}^{L}\sum_{m=1}^{M_r}  J_{ir}(k,m) x_{m}^{(r)} -  \sum_{\underset{r \neq i}{r=1}}^{L}\sum_{m=1}^{M_r}  x_{m}^{(r)}  x^{(i)}_{k} J_{ir}(k,m) \right)
$$ 
$$
=  x^{(i)}_{k} \left(  -\sum_{\underset{r \neq i}{r=1}}^{L}\sum_{m=1}^{M_r} x_{m}^{(r)}\sum_{   l=1 }^{M_{i}} x^{(i)}_{l} J_{ir}(l,m)  +   \sum_{\underset{r \neq i}{r=1}}^{L}\sum_{m=1}^{M_r}  J_{ir}(k,m) x_{m}^{(r)}  \right)    
$$
$$
=  x^{(i)}_{k} \left(   \sum_{\underset{r \neq i}{r=1}}^{L}\sum_{m=1}^{M_i}  J_{ir}(k,m) x_{m}^{(r)}   
-\sum_{\underset{r \neq i}{r=1}}^{L}\sum_{m=1}^{M_r} x_{m}^{(r)}\sum_{   l=1 }^{M_{i}} x^{(i)}_{l} J_{irl}(l,m)       \right)  
$$ 
\begin{equation}\label{del21} 
=  x^{(i)}_{k} \left(   \sum_{\underset{r \neq i}{r=1}}^{L}\sum_{m=1}^{M_r} \left[  J_{ir}(k,m)     
- \sum_{   l=1 }^{M_{i}} x^{(i)}_{l} J_{ir}(l,m)      \right] x_{m}^{(r)} \right).  
\end{equation} 
By adding    (\ref{del1}) and (\ref{del21}) give  the  expression  in (\ref{pairwise_selection3}) as desired.
\qed \\

\section{Appendix C: Inverse of  $D^{(i)}({\bf x})$  } \label{diffinverse}   
\renewcommand{\theequation}{\rm{C}.\arabic{equation}}
\setcounter{equation}{0}

 We establish  next  the inverse of the matrix $D^{(i)}({\bf x})$ explicitly.   It holds by  (\ref{driftlemma44}) 
 that  $D^{(i)}({\bf x})$ depends only on $x^{(i)}$. 
The result is as such known, see, e.g.,   \cite[p.1262]{papangelou2000large}, but not widely publicized, so   we prove  it for the sakes  of 
completeness and easy reference. The simple  important fact that emerges is that the inverse of  $D^{(i)}({\bf x})$
does not exist  at the facets of any ${\bf K}_{i}$.    
\begin{lemma}\label{ealemma} 
Assume that $ x_{k}^{(i)} >0$   for every $k \in \{ 1, \ldots, M_{i}  \}$. Then the  inverse matrix $D_{i}^{-1}(x^{(i)})$ 
is given  by  
\begin{equation}\label{eamatrixfact}
D_{i}^{-1}(x^{(i)}) = \left\{ \frac{1}{x^{(i)}_{l}} \delta_{l,k}  +   \frac{1}{x^{(i)}_{M_{i}}}     \right \}_{l=1, k=1}^{M_{i}-1, M_{i}-1}.   
\end{equation}
\end{lemma}  
{\em Proof}: We check that   $D_{i}^{-1}(x^{(i)})D_{i}(x^{(i)})=D_{i}(x^{(i)})D_{i}^{-1}(x^{(i)})= I$, where  $I$ is the  $  M_{i}-1   \times   M_{i}-1   $ unit matrix.     By (\ref{driftlemma44}) we have  
$$
\sum_{k=1}^{ M_{i}-1 }  \left(\frac{1}{x^{(i)}_{l}} \delta_{k,l}  +   \frac{1}{x^{(i)}_{M_{i}}} \right)  d_{km} = \sum_{k=1}^{ M_{i}-1 } \left(\frac{1}{x^{(i)}_{l}} \delta_{k,l}  +    \frac{1}{x^{(i)}_{M_{i}}} \right)\left(   x^{(i)}_{k} \delta_{k,m} -x^{(i)}_{k}x^{(i)}_{m} \right)   
$$  
$$
= \sum_{k=1}^{ M_{i}-1 }\frac{1}{x^{(i)}_{l}} \delta_{k,l}   x^{(i)}_{k} \delta_{k,m}  -  \sum_{k=1}^{M_{i-1}}
\frac{1}{x^{(i)}_{l}} \delta_{k,l} x^{(i)}_{k}x^{(i)}_{m} +   
\sum_{k=1}^{ M_{i}-1 }  \frac{1}{x^{(i)}_{M_{i}}}   x^{(i)}_{k} \delta_{k,m}  -   \sum_{k=1}^{M_{i}-1}    \frac{1}{x^{(i)}_{M_{i}}} 
x^{(i)}_{k}x^{(i)}_{m}     
$$
$$
=  \delta_{l,m} - x^{(i)}_{m} +   \frac{1}{x^{(i)}_{M_{i}}} x^{(i)}_{m}  -  \frac{1}{x^{(i)}_{M_{i}}} x^{(i)}_{m} \sum_{k=1}^{ M_{i}-1 } x_{k} 
$$
$$
=  \delta_{l,m} - x^{(i)}_{m} +   \frac{1}{x^{(i)}_{M_{i}}} x^{(i)}_{m}  -  \frac{1}{x^{(i)}_{M_{i}}} x^{(i)}_{m}  \left( 1-  x^{(i)}_{M_{i}} \right)   
$$
$$
=  \delta_{l,m} - x^{(i)}_{m} +   \frac{1}{x^{(i)}_{M_{i}}} x^{(i)}_{m}   -  \frac{1}{x^{(i)}_{M_{i}}} x^{(i)}_{m}   + x^{(i)}_{m}    
 = \delta_{l,m},    
$$
and the assertion in (\ref{eamatrixfact}) holds,  as claimed. \qed \\

\section*{Appendix D: Proof  of   lemma  \ref{bassats1} } \label{piilemma}   
\renewcommand{\theequation}{\rm{D}.\arabic{equation}}
\setcounter{equation}{0} 
 We start with a result  needed in the  proof of  the  lemma \ref{bassats1}, but which also  shows that $ \pi({\bf x})$ is the non-normalized  invariant density for 
 \begin{equation}
  dX(t) =  {\bf \mu}(X(t))dt +    {\bf D}^{1/2}\left(X(t)\right)d {\bf W}(t), 
 \end{equation} 
which  is (\ref{bigd55}) with   the  Svirezhev-Shahshahani gradient form removed. 
 \begin{lemma}\label{bassats}
 Assume    (\ref{indmut223}) and (\ref{dircond}).   
 Then,    for all  ${\bf x} \in \times_{i=1}^{L}{\bf K}_{i}$, $k=1, \ldots, M_{i}-1$,  and $i=1, \ldots, L$, 
 we have  with $g^{(i)}_{k} =u_{k}^{(i)}  -   \bar{u}^{(i)}   x_k^{(i)} $   that   
 \begin{equation}\label{flowlocusallellepart1} 
 \frac{1}{2}\sum_{l=1}^{M_i-1}\frac{\partial }{ \partial x_{l}^{(i)}}\left[d_{kl}^{(i)}({\bf x})\pi({\bf x}) \right] =  g_{k}^{(i)} \pi({\bf x}).
 \end{equation}  
  \end{lemma}   
\noindent {\em Proof}:  The proof is a straightforward but lengthy computation, but is recapitulated here for the sake of completeness.   For any $k=1, \ldots, M_{i}-1 $ we have   


\begin{eqnarray}\label{matrisfp2} 
\sum_{l=1}^{M_i-1}\frac{\partial }{ \partial x_{l}^{(i)}}  \left[d_{kl}^{(i)}({\bf x}) \pi({\bf x})\right]    
& = & \nonumber \\ 
&  &   \\ 
 \sum_{l=1}^{M_i-1}\frac{\partial }{ \partial x_{l}^{(i)}}  \left[d_{ kl}^{(i)}({\bf x})\right] \pi({\bf x}) & +  & 
 \sum_{l=1}^{M_i-1}d_{kl}^{(i)} ({\bf x})  \frac{\partial }{ \partial x_{l}^{(i)}}  \pi({\bf x}). \nonumber    
\end{eqnarray} 
 We evaluate first the second  term in the right hand side of  (\ref{matrisfp2}). For any  
 $x_{l}^{(i)}$ we get  by  straightforward  differentiation and  rearrangement that     
  \begin{equation}\label{rhsj1}
\frac{\partial }{ \partial x_{l}^{(i)}}   \pi({\bf x})   =  \left[  \frac{1}{ x_l^{(i)}}\left(2 u_{l}^{(i)} - 1 \right) 
 -   \frac{ 1 }{1-\sum_{k=1}^{M_i-1}x_{k}^{(i)}} \left (2 u_{M_{i}}^{(i)}  -1 \right)    \right ]  \pi({\bf x}).
 \end{equation} 
For this we note that $x_{l}^{(i)}$ is a variable in  one and only   one of the factors in  $\pi({\bf x})$.   
 Thus we get  by (\ref{driftlemma44}) that   
\begin{eqnarray}\label{pdel1}  
  \sum_{l=1}^{M_i-1}d_{kl}^{(i)} ({\bf x}) \frac{1}{ x_l^{(i)}}\left(2 u_{l}^{(i)} - 1 \right) & = &  \nonumber \\
-x_{k}^{(i)} \sum_{l=1, l\neq k}^{M_{i}-1} \left(2 u_{l}^{(i)} - 1 \right) & +&  
(1- x_{k}^{(i)}) \left(2 u_{k}^{(i)} - 1 \right) \nonumber  \\ 
= -x_{k}^{(i)} \sum_{l=1 }^{M_{i}-1} \left(2 u_{l}^{(i)} - 1 \right)   & + &  \left(2 u_{k}^{(i)} - 1 \right). 
\end{eqnarray}     
Next  (\ref{driftlemma44})  gives   
\begin{equation}\label{pdel2}  
  \sum_{l=1}^{M_i-1}d_{kl}^{(i)} ({\bf x})  \frac{ 1 }{1-\sum_{k=1}^{M_i-1}x_{k}^{(i)}} \left (2 u_{M_{i}}^{(i)}  -1 \right) =  \left (2 u_{M_{i}}^{(i)}  -1 \right)  \frac{ 1 }{1-\sum_{k=1}^{M_i-1}x_{k}^{(i)}}   \sum_{l=1}^{M_i-1}d_{kl}^{(i)} ({\bf x})  
\end{equation}   
Here 
$$  
 \sum_{l=1}^{M_i-1}d_{kl}^{(i)} ({\bf x}) =   -x^{(i)}_{k} \sum_{l=1, l \neq k }^{M_i-1}x^{(i)}_{l}  
+  \left( 1- x^{(i)}_{k} \right) x^{(i)}_{k}
$$ 
$$
= -x^{(i)}_{k} \sum_{l=1    }^{M_i-1}x^{(i)}_{l} + x^{(i)}_{k}  
$$
$$
=    x^{(i)}_{k} \left(  1-  \sum_{l=1    }^{M_i-1}x^{(i)}_{l} \right). 
$$
This means that   in the right hand  side of  (\ref{pdel2}) we get    
\begin{equation}\label{pdel3}  
    \left (2 u_{M_{i}}^{(i)}  -1 \right)  \frac{ 1 }{1-\sum_{k=1}^{M_i-1}x_{k}^{(i)}}   \sum_{l=1}^{M_i-1}d_{kl}^{(i)} ({\bf x})  =  x^{(i)}_{k}  \left(2 u_{M_{i}}^{(i)}  -1 \right).    
\end{equation}     
Hence we have in the right hand side of (\ref{rhsj1})  in view of (\ref{pdel1}) and  (\ref{pdel3})       
\begin{eqnarray}\label{delfin} 
  \sum_{l=1}^{M_i-1}d_{kl}^{(i)} ({\bf x})  \frac{\partial }{ \partial x_{l}^{(i)}}  \pi({\bf x}) & = &  \nonumber \\ 
 -x_{k}^{(i)} \sum_{l=1 }^{M_{i}-1} \left(2 u_{l}^{(i)} - 1 \right)   & + &  \left(2 u_{k}^{(i)} - 1 \right) -   x^{(i)}_{k}  \left (2 u_{M_{i}}^{(i)}  -1 \right)  \nonumber \\  
 & & \\ \nonumber 
 = -x_{k}^{(i)} \sum_{l=1 }^{M_{i}} \left(2 u_{l}^{(i)} - 1 \right) & + &  \left(2 u_{k}^{(i)} - 1 \right). 
\end{eqnarray} 
But  we now observe that   
$$
 -x_{k}^{(i)} \sum_{l=1 }^{M_{i}} \left(2 u_{l}^{(i)} - 1 \right) +   \left(2 u_{k}^{(i)} - 1 \right) 
$$
$$
=  2( u_{k}^{(i)} -   \bar{u}^{(i)}  x_{k}^{(i)})  + M_{i}x_{k}^{(i)}  -1,      
$$  
where we used the notation in  (\ref{eabar}). When we substitute this  in the right hand side of    
(\ref{delfin})  we obtain   
\begin{equation}\label{delfin2} 
  \sum_{l=1}^{M_i-1}d_{kl}^{(i)} ({\bf x})  \frac{\partial }{ \partial x_{l}^{(i)}}  \pi({\bf x}) = \left[  2( u_{k}^{(i)} -   \bar{u}^{(i)} x_{k}^{(i)})  + M_{i}x_{k}^{(i)}  -1 \right] \pi({\bf x}). 
\end{equation}  
Next we  compute  the first   term in the right hand side of  (\ref{matrisfp2}). By   (\ref{driftlemma44})   
$$ 
\sum_{l=1}^{M_i-1} \frac{\partial }{ \partial x_{l}^{(i)}} d_{kl}^{(i)} ({\bf x})  =   
- \sum_{l=1, l \neq k}^{M_{i} -1}x_{k}^{(i)}  +(1 -2x_{k}^{(i)})   
$$ 
$$
= -x_{k}^{(i)} (M_{i} -2)  +  1 -2x_{k}^{(i)} = -x_{k}^{(i)} M_{i} +1     
$$
or,  
\begin{equation}\label{delfin3}       
\sum_{l=1}^{M_i-1} \frac{\partial }{ \partial x_{l}^{(i)}} d_{kl}^{(i)} ({\bf x}) = -x_{k}^{(i)} M_{i} +1.
\end{equation}  
By (\ref{delfin2})  and  (\ref{delfin3}) we obtain in (\ref{matrisfp2}) that 
$$
\sum_{l=1}^{M_i-1}\frac{\partial }{ \partial x_{l}^{(i)}}  \left[d_{ kl}^{(i)}({\bf x})\right] \pi({\bf x})  +  
 \sum_{l=1}^{M_i-1}d_{kl}^{(i)} ({\bf x})  \frac{\partial }{ \partial x_{l}^{(i)}}  \pi({\bf x}) 
 $$
 $$
 = \left[  2( u_{k}^{(i)} -  \bar{u}^{(i)}  x_{k}^{(i)})  + M_{i}x_{k}^{(i)}  -1   -x_{k}^{(i)} M_{i} +1 \right]\pi({\bf x})  
 $$
 $$
= \left[  2( u_{k}^{(i)} -   \bar{u}^{(i)}  x_{k}^{(i)}) \right]\pi({\bf x}).  
$$  
In view of   lemma  \ref{ealemma4}  we get  the result as claimed  in the   lemma.  \qed  \\ 
An  inspection of the proof above shows that  it is strictly valid only in the interior of  
${\bf K}_{i}$. However, the final result can obviously be extended to the boundary by continuity.  
Next we prove lemma \ref{bassats1}. \\  
\noindent  {\em Proof of lemma \ref{bassats1}}: For any $k=1, \ldots, M_{i}-1 $ we have   
\begin{eqnarray}\label{matrisfbp} 
\sum_{l=1}^{M_i-1}\frac{\partial }{ \partial x_{l}^{(i)}}  \left[d_{kl}^{(i)}({\bf x}) P\left({\bf x}\right) \right]    
& = & \nonumber \\ 
&  &   \\ 
 \sum_{l=1}^{M_i-1}\frac{\partial }{ \partial x_{l}^{(i)}}  \left[d_{ kl}^{(i)}({\bf x})\right] P\left({\bf x}\right) & +  & 
 \sum_{l=1}^{M_i-1}d_{kl}^{(i)} ({\bf x})  \frac{\partial }{ \partial x_{l}^{(i)}} P\left({\bf x}\right). \nonumber    
\end{eqnarray}  
Here 
$$
 \sum_{l=1}^{M_i-1}d_{kl}^{(i)} ({\bf x})  \frac{\partial }{ \partial x_{l}^{(i)}} P\left({\bf x}\right)= 
   e^{2 V({\bf x})} \sum_{l=1}^{M_i-1}d_{kl}^{(i)} ({\bf x})  \frac{\partial }{ \partial x_{l}^{(i)}} \pi \left({\bf x}\right)  + P  \left({\bf x}\right)  \sum_{l=1}^{M_i-1}d_{kl}^{(i)} ({\bf x})    
 \left(2 V^{'}_{ x_{l}^{(i)}} \left({\bf x}\right)     \right).     
$$
Thus we have  
$$
 \sum_{l=1}^{M_i-1}\frac{\partial }{ \partial x_{l}^{(i)}}  \left[d_{ kl}^{(i)}({\bf x})\right] P\left({\bf x}\right) +  
 \sum_{l=1}^{M_i-1}d_{kl}^{(i)} ({\bf x})  \frac{\partial }{ \partial x_{l}^{(i)}} P\left({\bf x}\right) =
 $$
 $$
 =  e^{2 V({\bf x})} \left[ \sum_{l=1}^{M_i-1}\frac{\partial }{ \partial x_{l}^{(i)}}  \left[d_{ kl}^{(i)}({\bf x})\right] \pi\left({\bf x}\right) +  
   \sum_{l=1}^{M_i-1}d_{kl}^{(i)} ({\bf x})  \frac{\partial }{ \partial x_{l}^{(i)}} \pi \left({\bf x}\right) \right] 
 $$
$$                                                                                                                                                                                                                                                                                                                                                                                               
+ P  \left({\bf x}\right)  \sum_{l=1}^{M_i-1}d_{kl}^{(i)} ({\bf x})    
  2 V^{'}_{ x_{l}^{(i)}} \left({\bf x}\right).
$$
But by lemma   \ref{bassats}, or (\ref{flowlocusallellepart1}) in the preceding,    we  get above     
$$
= 2  e^{2 V({\bf x})}   g_{k}^{(i)} \pi({\bf x}) +  P  \left({\bf x}\right)  \sum_{l=1}^{M_i-1}d_{kl}^{(i)} ({\bf x})    
  2 V^{'}_{ x_{l}^{(i)}} \left({\bf x}\right).  
$$
$$
=  2  P  \left({\bf x}\right)   \left(  g_{k}^{(i)} +  \sum_{l=1}^{M_i-1}d_{kl}^{(i)} ({\bf x})    
   V^{'}_{ x_{l}^{(i)}} \left({\bf x}\right) \right).
$$
\qed

\section*{Appendix E: Diffusion approximation: technical steps   } \label{diffapprx}   
\renewcommand{\theequation}{\rm{E}.\arabic{equation}}
\setcounter{equation}{0} 
 This Appendix  contains the proofs of the   technical  conditions on limits of the sequences of  conditional incremental moments for the sequence of Markov chains  $X^{(N)}$   required in   the weak convergence statement 
 of  proposition \ref{diffapprthm}.   
 

     
\subsection{Conditional Expectation of the Increments}\label{expincfre} 
    Let us set  for $k=1, \ldots, M_{i}$   
 \begin{equation}\label{onelinterpol23}
 \mu_{i,k}^{(N)}\left( \frac{{\bf j}}{N} \right) \stackrel{\rm def}{=}  N  \sum_{r \in {\bf J}_{(i)}(N)} \left( \frac{r^{(i)}_{k}}{N} -\frac{j^{(i)}_{k}}{N} \right) P_{{\bf j}{\bf r}}    
 \end{equation}
which is  
$$
\mu_{i,k}^{(N)}\left( \frac{{\bf j}}{N} \right)=  N E \left[ \left(X_{i,k}^{(N)}(n+1) -X_{i,k}^{(N)}(n) \right)  \mid X^{(N)}(n)= \frac{{\bf j}}{N}  \right].  
$$ 
This is, of  course, the conditional expectation of the difference ratio    
$$
  \frac{\left(X_{i,k}^{(N)}(n+1) -X_{i,k}^{(N)}(n) \right)}{ \frac{1}{N} }.   
$$  
An  analogous  interpretation holds for the expressions studied  in the other subsections  of this Appendix E. 
We shall next evaluate  (\ref{onelinterpol23})   and then expand it  as a function  of $1/N$ evoking the rescalings  in the  assumption \ref{scaling}.   For ease of writing we drop for the moment  the subscript  for  locus in the computations that follow  
in this section. Hence the analysis holds  for any $ X_{k}^{(N)}(n+1) \equiv X_{i,k}^{(N)}(n)$ in (\ref{onelinterpol22}).      
$$
  N E \left[ \left(X_{k}^{(N)}(n+1) -X_{k}^{(N)}(n) \right)  \mid X^{(N)}(n)= \frac{{\bf j}}{N}  \right] 
$$
$$
=   E \left[ \left(Y_{k}^{(N)}\left( \frac{n+1}{N} \right) -Y_{k}^{(N)}\left( \frac{n }{N} \right) \right)  \mid X^{(N)}(n)= \frac{{\bf j}}{N}  \right]  
$$  
$$
=  E \left[ Y_{k}^{(N)}\left( \frac{n+1}{N} \right) \mid X^{(N)}(n)= \frac{{\bf j}}{N}  \right]   -  E \left[  Y_{k}^{(N)}\left( \frac{n }{N} \right)\mid X^{(N)}(n)= \frac{{\bf j}}{N}  \right]   
$$
$$
= E \left[ Y_{k}^{(N)}\left( \frac{n+1}{N} \right) \mid Y^{(N)}\left( \frac{n}{N} \right)=  {\bf j}   \right] - j_{k}  
$$
and by a property of  the multinomial distribution  (\ref{transitionmech}) 
$$
=  N p_k({\bf j}) - j_{k}.
$$
Thus when we return to  the full notations   
\begin{equation} 
\label{multinom_exp22}
 \mu_{i,k}^{(N)}\left( \frac{{\bf j}}{N} \right)=   N p_k({\bf j}) - j_{k}= N ( p_k({\bf j}) -x_k)  = N ( p^{(i)}_k({\bf j}) -x^{(i)}_k).  
 \end{equation}
 By  (\ref{not_approximated1})  we obtain   
\begin{equation}\label{not_approximated}
 N ( p^{(i)}_k({\bf j}) -x^{(i)}_k) =N \left\{\sum_{\underset{l \neq k}{l=1}}^{M } \left [ \upsilon^{(i)}_{lk}  x^{(i)}_l  \left ( \frac{\bar{v}_l^{(i)}}{\bar{v}} \right ) - \upsilon^{(i)}_{kl} x^{(i)}_k  \left ( \frac{\bar{v}_k^{(i)}}{\bar{v}} \right ) \right ] + x^{(i)}_k  \left (\frac{1}{\bar{v}}\right ) \left [ \bar{v}_k^{(i)}-\bar{v} \right ] \right\}.\\
\end{equation}  
When   (\ref{popscale1}),  
 (\ref{pairwise_selection11}) and (\ref{pairwise_selection12})  are inserted  in  (\ref{not_approximated})    
we obtain 
 $$
 \mu_{i,k}^{(N)}\left( \frac{{\bf j}}{N} \right)= N \cdot \left\{ p^{(i)}_k({\bf j}) -x^{(i)}_k \right \} = 
 $$
 $$ 
 N \left\{ \sum_{\underset{l \neq k}{l=1}}^{M } \left [\frac{u_{lk}^{(i)} }{N}  x^{(i)}_l  \left ( \frac{1 +  \frac{ \bar{m}_l^{(i)}} {N} }{ 1 +  \frac{\bar{m}}{N}} \right ) -\frac{u_{kl}^{(i)} }{N} x^{(i)}_k  \left ( \frac{1 +  \frac{ \bar{m}_k^{(i)}} {N} }{ 1 +  \frac{\bar{m}}{N}} \right ) \right ] \right. 
$$
$$  
\left.  + x^{(i)}_k  \left (\frac{1}{1 +  \frac{\bar{m}}{N}}\right )  \frac{ (\bar{m}_k^{(i)} -\bar{m}) }{N } \right\}. 
 $$   
If $N \rightarrow +\infty$, this expression clearly converges to  
$$
  p_{k}^{(i)}({\bf x}) =  \sum_{\underset{l \neq k}{l=1}}^{M_i} \left [ u_{lk}^{(i)}x_l^{(i)} - u_{kl}^{(i)}x_k^{(i)} \right ]+  x_k^{(i)} \left ( \bar{m}_k^{(i)}-\bar{m}\right ).
$$
\begin{lemma}\label{driftlemma1} 
\begin{equation}\label{driftlimit}  
\mu_{i,k}^{(N)}\left( \frac{{\bf j}}{N} \right)  \rightarrow p_{k}^{(i)}({\bf x} ), 
\end{equation} 
as   $ N  \rightarrow  + \infty$,  uniformly in ${\bf x} \in  \times_{i=1}^{L}{\bf  K}_{i}$,   where   for ${\bf x} \in  \times_{i=1}^{L}{\bf  K}_{i}$ 
\begin{equation}\label{driftlimit2}    
  p_{k}^{(i)}({\bf x}) =  \sum_{\underset{l \neq k}{l=1}}^{M_i} \left [ u_{lk}^{(i)}x_l^{(i)} - u_{kl}^{(i)}x_k^{(i)} \right ]+  x_k^{(i)} \left ( \bar{m}_k^{(i)}-\bar{m}\right )    
\end{equation} 
\end{lemma} 
{\em Proof:} 
 It  remains to prove that  the convergence in  (\ref{driftlimit}) is in fact uniform in $x^{(i)}$. 
To see this, let us check  
 $$
 | \mu_{i,k}^{(N)}\left( {\bf x} \right) -  p_{k}^{(i)}({\bf x} )|  
 = \left|  \sum_{\underset{l \neq k}{l=1}}^{M } \left [ u_{lk}^{(i)}  x^{(i)}_l  \left ( \frac{1 +  \frac{ \bar{m}_l^{(i)}} {N} }{ 1 +  \frac{\bar{m}}{N}} -1 \right ) - u_{kl}^{(i)}  x^{(i)}_k  \left ( \frac{1 +  \frac{ \bar{m}_k^{(i)}} {N} }{ 1 +  \frac{\bar{m}}{N}} -1 \right ) \right.  \right. 
$$
$$  
\left. \left.   + x^{(i)}_k  \left (\frac{1}{1 +  \frac{\bar{m}}{N}} -1 \right )    (\bar{m}_k^{(i)} -\bar{m}) \right ] \right | 
 $$   
$$
\leq    \sum_{\underset{l \neq k}{l=1}}^{M } \left|\left [ u_{lk}^{(i)}  x^{(i)}_l  \left ( \frac{\frac{1}{N}  \left( \bar{m}_l^{(i)}  -    \bar{m} \right) } {1 +  \frac{\bar{m}}{N} }   \right ) \right|  +   \left|  u_{kl}^{(i)}  x^{(i)}_k  \left ( \frac{\frac{1}{N} \left( \bar{m}_k^{(i)} -    \bar{m} \right)} {1 +  \frac{\bar{m}}{N} }   \right )  \right ] \right|  
$$
$$  
 + \left. \left|   x^{(i)}_k \left (-\frac{ \frac{\bar{m}}{N}}{1 +  \frac{\bar{m}}{N}}  \right )    (\bar{m}_k^{(i)} -\bar{m})  \right]\right|
$$
$$
\leq  \frac{1}{N}  \sup_{N} \sup_{x^{(i)} \in {\bf K}_{i}} \left\{ \sum_{\underset{l \neq k}{l=1}}^{M } \left|\left [ u_{lk}^{(i)}  x^{(i)}_l  \left ( \frac{   \left( \bar{m}_l^{(i)}  -    \bar{m} \right) } {1 +  \frac{\bar{m}}{N} }   \right ) \right|  +  \left|  u_{kl}^{(i)}  x^{(i)}_k  \left ( \frac{  \left( \bar{m}_k^{(i)} -    \bar{m} \right)} {1 +  \frac{\bar{m}}{N} }   \right )  \right ] \right|  \right.  
$$
$$  
 +\left.  \left. \left|   x^{(i)}_k \left (-\frac{\bar{m}}{1 +  \frac{\bar{m}}{N}}  \right )    (\bar{m}_k^{(i)} -\bar{m})  \right]\right| \right\} \leq  \frac{A}{N}.   
$$
I.e., 
\begin{equation}\label{unifconv} 
| \mu_{i,k}^{(N)}\left( {\bf x} \right) -  p_{k}^{(i)}({\bf x} )|   \leq  \frac{A}{N}.  
\end{equation}  
\qed \\ 
The result  in (\ref{unifconv}) agrees with the notion of uniform convergence of the incremental  conditional  moment characteristics of a sequence of  Markov chains  in both   
\cite[Theorem 7.1]{Durrett1996},   and      \cite[p. 642]{sato1976diffusion}.

 \subsection{Conditional Covariances of the Increments} 

Again we take   for any $ X_{k}^{(N)}(n+1) \equiv X_{i,k}^{(N)}(n)$ in (\ref{onelinterpol22}).      
 Next we set 
 \begin{equation}\label{onelinterpol24}
 d_{kl}^{(N)}\left( \frac{{\bf j}}{N} \right)\stackrel{\rm def}{=}  N  \sum_{r \in {\bf J}_{(i)}(N)} 
 \left( \frac{r_{k}}{N} -\frac{j_{k}}{N} \right)\left( \frac{r_{l}}{N} -\frac{j_{l}}{N} \right)  P_{\bf jr}, 
\end{equation}
which is  
$$
  d_{kl}^{(N)}\left( \frac{{\bf j}}{N} \right)=  N  E \left[ (X_{k}^{(N)}(n+1) -X_{k}^{(N)}(n)) 
 (X_{l}^{(N)}(n+1) -X_{l}^{(N)}(n))    \mid X^{(N)}(n)= \frac{{\bf j}}{N}  \right].   
$$ 
Then    
 $$
 d_{kl}^{(N)}\left( \frac{{\bf j}}{N} \right)    
= \frac{1}{N}  E \left[ (Y_{k}^{(N)}\left( \frac{n+1}{N} \right) -Y_{k}^{(N)}\left( \frac{n}{N} \right) \cdot 
 (Y_{l}^{(N)}\left( \frac{n+1}{N} \right) -Y_{l}^{(N)}\left( \frac{n}{N} \right))    \mid X^{(N)}(n)= \frac{{\bf j}}{N}  \right]. 
$$
 Take first  $k\neq l$. Then  
 $$
 E \left[ (Y_{k}^{(N)}\left( \frac{n+1}{N} \right) -Y_{k}^{(N)}\left( \frac{n}{N} \right)) \cdot 
 (Y_{l}^{(N)}\left( \frac{n+1}{N} \right) -Y_{l}^{(N)}\left( \frac{n}{N} \right))    \mid X^{(N)}(n)= \frac{{\bf j}}{N}  \right]= 
 $$ 
$$
= E \left[  Y_{k}^{(N)}\left( \frac{n+1}{N} \right)Y_{l}^{(N)}\left( \frac{n+1}{N} \right)   \mid X^{(N)}(n)= \frac{{\bf j}}{N}  \right]  $$
$$ -   
E \left[  Y_{k}^{(N)}\left( \frac{n+1}{N} \right)Y_{l}^{(N)}\left( \frac{n}{N} \right)   \mid X^{(N)}(n)= \frac{{\bf j}}{N}  \right]  
$$
$$
 -  E \left[  Y_{k}^{(N)}\left( \frac{n}{N} \right)Y_{l}^{(N)}\left( \frac{n+1}{N} \right)   \mid X^{(N)}(n)= \frac{{\bf j}}{N}  \right] 
 +  E \left[  Y_{k}^{(N)}\left( \frac{n }{N} \right)Y_{l}^{(N)}\left( \frac{n}{N} \right)   \mid X^{(N)}(n)= \frac{{\bf j}}{N}  \right]=     
$$
$$
=  E \left[ (Y_{k}^{(N)}\left( \frac{n+1}{N} \right)Y_{l}^{(N)}\left( \frac{n+1}{N} \right)   \mid X^{(N)}(n)= \frac{{\bf j}}{N}  \right] 
- N p_{k}({\bf j})\cdot j_{l}  -   N p_{l}({\bf j})\cdot j_{k} + j_{l}j_{k},
$$
where we invoked  the appropriate moments of  the multinomial distribution  in (\ref{transitionmech}).  By the  same token the first term in the right hand side of the inequality  above is evaluated as 
$$
E \left[ (Y_{k}^{(N)}\left( \frac{n+1}{N} \right)Y_{l}^{(N)}\left( \frac{n+1}{N} \right)   \mid X^{(N)}(n)= \frac{{\bf j}}{N}  \right]  = 
$$
$$
{\rm Cov}\left( Y_{k}^{(N)}\left( \frac{n+1}{N} \right), Y_{l}^{(N)}\left( \frac{n+1}{N} \right)\right) 
$$
$$ 
+ E \left[  (Y_{k}^{(N)}\left( \frac{n+1}{N} \right)  \mid X^{(N)}(n)= \frac{{\bf j}}{N}  \right] E \left[   Y_{l}^{(N)}\left( \frac{n+1}{N} \right)   \mid X^{(N)}(n)= \frac{{\bf j}}{N}  \right]    
$$  
$$
= -N p_k({\bf j})p_l({\bf j}) +   N^{2}p_{k}({\bf j})p_{l}({\bf j}).  
$$
Hence we have obtained 
$$
 d_{kl}^{(N)}\left( \frac{{\bf j}}{N} \right)= \frac{1}{N} \left[ -N p_k({\bf j})p_l({\bf j}) +   N^{2}p_{k}({\bf j})p_{l}({\bf j})  - N p_{k}({\bf j})\cdot j_{l}  -   N p_{l}({\bf j})\cdot j_{k} + j_{l}j_{k} \right]    
$$  
$$
= \frac{1}{N} \left[ -N p_k({\bf j})p_l({\bf j})      
+ N(  p_k({\bf j}) - x_{k})\cdot N  ( p_l({\bf j}) - x_{l} )     
 \right]    
$$ 
$$
=  -  p_k({\bf j})p_l({\bf j})       
+  (  p_k({\bf j}) - x_{k})\cdot N  ( p_l({\bf j}) - x_{l} ). 
$$
In view of (\ref{multinom_exp22}), (\ref{driftlimit}) and (\ref{unifconv})  we have 
$ p_l({\bf j}) = x_{l} + O\left( \frac{1}{N}\right)$ and $ p_k({\bf j}) = x_{k} + O\left( \frac{1}{N}\right)$ 
and we proved the following lemma. 
\begin{lemma}\label{driftlemma2} 
\begin{equation}\label{difflimit} 
  d_{kl}^{(N)}\left( \frac{{\bf j}}{N} \right) \rightarrow -  x_k x_l, \quad k \neq l,  
\end{equation} 
 as $N \rightarrow +\infty$,    the   convergence  is uniform in  ${\bf x}$.  
\end{lemma} \qed \\
The limiting genetic drift for $k=l$  is obtained readily, too.  We have 
$$
  d_{kk}^{(N)}\left( \frac{{\bf j}}{N} \right)=  N  E \left[ (X_{k}^{(N)}(n+1) -X_{k}^{(N)}(n))^{2} 
      \mid X^{(N)}(n)= \frac{{\bf j}}{N}  \right]
$$ 
$$
=   \frac{1}{N}  \left\{  E \left[  Y_{k}^{(N)}\left( \frac{n+1}{N} \right)^{2} 
      \mid X^{(N)}(n)= \frac{{\bf j}}{N}  \right] -2   E \left[ Y_{k}^{(N)}\left( \frac{n+1}{N} \right) Y_{k}^{(N)}\left( \frac{n}{N} \right)  
      \mid X^{(N)}(n)= \frac{{\bf j}}{N}  \right] \right. 
   $$ 
   $$
  \left. +     E \left[  Y_{k}^{(N)}\left( \frac{n+1}{N} \right)^{2} 
      \mid X^{(N)}(n)= \frac{{\bf j}}{N}  \right] \right \}
$$
$$
=     \frac{1}{N}  \left\{  {\rm Var} \left[  Y_{k}^{(N)}\left( \frac{n+1}{N} \right)   
      \mid X^{(N)}(n)= \frac{{\bf j}}{N}  \right] + \left(E \left[  Y_{k}^{(N)}\left( \frac{n+1}{N} \right) 
      \mid X^{(N)}(n)= \frac{{\bf j}}{N}  \right]\right)^{2} 
\right.
$$
$$
\left.      
         -2 N p_{k} ( {\bf j}) j_{k}     
   +     j_{k}^{2} \right\}
$$
and due to  (\ref{transitionmech}) it follows that  
$$
=  \frac{1}{N}  \left\{ N p_{k}({\bf j})( 1-  p_{k}({\bf j})) + N^{2} p_{k}({\bf j})^{2} -2  N p_{k} ( {\bf j}) j_{k}     
   +     j_{k}^{2}\right\}
$$
$$
= \frac{1}{N}  \left\{ N p_{k}({\bf j})( 1-  p_{k}({\bf j}))  + \left( N  p_{k}({\bf j}) -  j_{k}   \right)^{2}  \right\}
$$
$$
=   p_{k}({\bf j})( 1-  p_{k}({\bf j}))  +  \frac{1}{N} \left( N(  p_{k}({\bf j}) -  x_{k} )  \right)^{2}.  
$$
By (\ref{multinom_exp22}), (\ref{driftlimit}) and (\ref{unifconv})  we have 
$ p_l({\bf j}) = x_{l} + O\left( \frac{1}{N}\right)$ and $ p_k({\bf j}) = x_{k} + O\left( \frac{1}{N}\right)$.  
\begin{lemma}\label{driftlemma3}
\begin{equation}  \label{gdr_finished1}
 d_{kk}^{(N)}\left( \frac{{\bf j}}{N} \right) \rightarrow  x_{k}  ( 1-  x_{k}).
\end{equation} 
  uniformly,  as $N \rightarrow + \infty$. 
  \end{lemma}\qed \\ 
When the findings in the two lemmas above are collected to one statement  we have  for  $k=1, \ldots M_{i}$ and 
$l=1, \ldots M_{i}$ and for every locus $i$ the expression for $d_{kl}^{(i)}({\bf x})$. 
 
   Now we take  $  X_{i,k}^{(N)}(n)$ and $  X_{i,k}^{(N)}(n)$ in (\ref{onelinterpol22}) with two different loci.      
Let us set  
$$
  d_{kl, (i,j)}^{(N)}\left( \frac{{\bf j}}{N} \right)=  N  E \left[ (X_{i,k}^{(N)}(n+1) -X_{i,k}^{(N)}(n)) 
 (X_{j,l}^{(N)}(n+1) -X_{j,l}^{(N)}(n))    \mid X^{(N)}(n)= \frac{{\bf j}}{N}  \right].   
$$  
\begin{lemma}\label{driftlemma4} 
 \begin{equation}\label{mixed}  
N  E \left[ (X_{i,k}^{(N)}(n+1) -X_{i,k}^{(N)}(n)) 
 (X_{j,l}^{(N)}(n+1) -X_{j,l}^{(N)}(n))    \mid X^{(N)}(n)= \frac{{\bf j}}{N}  \right] \rightarrow 0. 
\end{equation}
as $ N \rightarrow +\infty$, uniformly. 
\end{lemma} 
{\em Proof:} Due to the assumption of conditional independence  (\ref{condindpe}) over loci  we get 
 $$
  =  N  E \left[ (X_{i,k}^{(N)}(n+1) -X_{i,k}^{(N)}(n))
   \mid X^{(N)}(n)= \frac{{\bf j}}{N}  \right] 
 E \left[   (X_{j,l}^{(N)}(n+1) -X_{j,l}^{(N)}(n))    \mid X^{(N)}(n)= \frac{{\bf j}}{N}  \right].
 $$  
 $$
=  \frac{1}{N} \left \{ N  E \left[ (X_{i,k}^{(N)}(n+1) -X_{i,k}^{(N)}(n))
   \mid X^{(N)}(n)= \frac{{\bf j}}{N}  \right] \right. 
   $$
   $$ 
   \left. 
 N E \left[   (X_{j,l}^{(N)}(n+1) -X_{j,l}^{(N)}(n))    \mid X^{(N)}(n)= \frac{{\bf j}}{N}  \right] \right \} 
 $$
 The computation in section \ref{expincfre} above   yields here 
 $$
 =   \frac{1}{N}  \mu_{i,k}^{(N)}\left( \frac{{\bf j}}{N} \right)  \mu_{j,l}^{(N)}\left( \frac{{\bf j}}{N} \right). 
 $$ 
 Since $\mu_{i,k}^{(N)}\left( \frac{{\bf j}}{N} \right)$ and   $\mu_{j,l}^{(N)}\left( \frac{{\bf j}}{N} \right)$  converge, as shown in section \ref{expincfre},    uniformly to finite limits as $N \rightarrow +\infty$, we get the lemma as asserted.  \qed  
 \subsection{No Jumps in the Limit}  
For the    diffusion approximation desired   we need to check the behaviour of (e.g.,) the fourth moment of the increments defined as 
\begin{equation}\label{onelinterpol25} 
 e_{{\bf j}{\bf r}, k}^{(N)}\left( \frac{{\bf j}}{N} \right) \stackrel{\rm def}{=}   N  \sum_{r \in {\bf J}_{(i)}(N)} \left( \frac{r_{k}}{N} -\frac{j_{k}}{N} \right)^{4}  P_{{\bf j}{\bf r}}.     
\end{equation}
 This is dependent on the locus $i$, i.e.,  $j_{k}= j^{(i)}_{k} $ and  $r_{k}= r^{(i)}_{k} $,but we omit once more this  for  reasons  of simplicity of  notation in the  calculations of this  subsection. 
 \begin{lemma}\label{driftlemma5}    
\begin{equation}\label{onelinterpol28}  
 e_{{\bf j}{\bf r}, k}^{(N)}\left( \frac{{\bf j}}{N} \right)  \rightarrow  0 
 \end{equation}  
uniformly, as $N  \rightarrow +\infty$.   
 \end{lemma} 
{\em Proof}: We have    
\begin{equation}\label{onelinterpol26} 
 e_{{\bf j}{\bf r}, k}^{(N)}\left( \frac{{\bf j}}{N} \right))  =   N  
E \left[  \left( X^{(N)}_{k}(n+1) -  X^{(N)}_{k}(n) \right)^{4} \mid   X (n) = \frac{{\bf j}}{N} \right].  
\end{equation}
$$
= \frac{1}{N^{3}}  
E \left[  \left( Y^{(N)}_{k}\left((\frac{n+1}{N}\right) -  Y^{(N)}_{k}\left((\frac{n}{N}\right) \right)^{4} \mid  Y^{(N)} \left((\frac{n}{N}\right) =  {\bf j}  \right].  
$$ 
First we insert  $N p_{k}( {\bf j})$, invoke  the  inequality $(a+b)^{4} \leq 2^{4}a^{4} +2^{4}b^{4}$  and obtain  the bound  
$$
\leq 2^{4} E \left[  \left( Y^{(N)}_{k}\left((\frac{n+1}{N}\right) - N p_{k}( {\bf j}) \right)^{4} \mid  Y^{(N)}(\left(\frac{n}{N}\right) =  {\bf j}  \right] $$
$$
 + 2^{4}  E \left[  \left(  N p_{k}( {\bf j}) - Y^{(N)}_{k}\left((\frac{n}{N}\right)    \right)^{4} \mid  Y^{(N)}(\left(\frac{n}{N}\right) =  {\bf j}  \right]. 
$$
 Here 
 $$
 E \left[  \left(  N p_{k}( {\bf j}) - Y_{k}\left(\frac{n}{N}\right)    \right)^{4} \mid  Y^{(N)}\left(\frac{n}{N}\right) =  {\bf j}  \right]  = E \left[  \left(  N p_{k}( {\bf j}) - j_{k}    \right)^{4} \mid  Y^{(N)}(\frac{n}{N}) =  {\bf j}  \right] 
 $$ 
$$
=   \left(  N p_{k}( {\bf j}) - j_{k}    \right)^{4}  =    \left(  N (p_{k}( {\bf j}) - x_{k}    \right)^{4}.   
$$
Hence, by (\ref{driftlimit}) and  (\ref{multinom_exp22}) as above  
$$
\frac{1}{N^{3}} \left(  N (p_{k}( {\bf j}) - x_{k})    \right)^{4}  \rightarrow 0, 
$$ 
uniformly in ${\bf x}$,  as $N \rightarrow \infty$. 
Next, we bound    
$$
E \left[  \left( Y^{(N)}_{k}\left(\frac{n+1}{N}\right) - N p_{k}( {\bf j}) \right)^{4} \mid  Y^{(N)}\left(\frac{n}{N}\right) =  {\bf j}  \right].
$$
 We know that  $ Y^{(N)}_{k}\left(\frac{n+1}{N}\right)$ conditioned on  $Y^{(N)}(\frac{n}{N}) =  {\bf j}$  has the   
multinomial distribution  (\ref{transitionmech}). This means that  $ Y^{(N)}_{k}(\frac{n+1}{N})$ is    in distribution equal  
$$
Y^{(N)}_{k}\left(\frac{n+1}{N} \right) \stackrel{d}{=}  \sum_{l=1}^{N} \xi_{l}, 
$$
where $ \xi_{k}= I\left( Z_{l} = k \right)$ is the indicator function of  the event,  $Z_{l}$, $l=1,2, \ldots N,$ are conditionally independent and identically distributed random variables  such that 
$
p_{k}({\bf j}) = P (Z_{l}= k), k=1, \ldots, M_{i}.     
$  
Here we apply  a technique from \cite[p.~308]{Durrett1996}.  
$$
E \left[  \left(Y^{(N)}_{k}(\frac{n+1}{N}) - N p_{k}( {\bf j}) \right)^{4} \mid  Y^{(N)}(\frac{n}{N}) =  {\bf j}  \right] =   
E \left[  \left(\sum_{l=1}^{N}  \xi_{k}   - N p_{k}( {\bf j}) \right)^{4} \mid  Y^{(N)}(\frac{n}{N}) =  {\bf j}  \right] 
$$ 
$$
= E \left[ \left( \sum_{l=1}^{N} ( \xi_{k}   -  p_{k}( {\bf j}))  \right)^{4} \mid  Y^{(N)}(\frac{n}{N}) =  {\bf j}  \right] 
$$   
$$
\leq   N E \left[ \left(  \xi_{k}   -  p_{k}( {\bf j})   \right)^{4} \mid  Y^{(N)}(\frac{n}{N}) =  {\bf j}  \right] + 6 \left( \begin{array}{ll} N \\ 2 \end{array}   \right) E \left[ \left(  \xi_{k}   -  p_{k}( {\bf j})   \right)^{2} \mid Y^{(N)}(\frac{n}{N}) =  {\bf j}  \right]  
$$
$$
\leq C N^{2}. 
$$
Above we evoked  the  inequality   $ E \left[ \left(  \xi_{k}   -  p_{k}( {\bf j})   \right)^{m} \mid  Y_{k}(\frac{n}{N}) =  {\bf j}  \right]   \leq 1$ for all $m \geq 1$. Hence we have that 
$$
\frac{1}{N^{3}} E \left[  \left( Y^{(N)}_{k}(\frac{n+1}{N}) - N p_{k}( {\bf j}) \right)^{4} \mid Y^{(N)}(\frac{n}{N}) =  {\bf j}  \right]   \leq  C N^{-1}.    
$$  
By the preceding  we have shown the asserted lemma. \qed

\end{document}